\documentclass{article}
\usepackage{graphicx} 
\usepackage{amsmath,amssymb,amsthm}
\usepackage{enumerate}
\usepackage{mathtools}
\usepackage{tikz}
\usepackage{tikz-3dplot}
\usepackage{comment}
\usepackage{url}

\newtheorem{theorem}{Theorem}[section]
\newtheorem*{theorem*}{Theorem}
  \newtheorem{corollary}[theorem]{Corollary}
  \newtheorem{definition}[theorem]{Definition}
   \newtheorem*{definition*}{Definition}
  \newtheorem{lemma}[theorem]{Lemma}
  \newtheorem{notation}[theorem]{Notation}
  \newtheorem{proposition}[theorem]{Proposition}
  \newtheorem{remark}[theorem]{Remark}
  \newtheorem{thmx}{Theorem}

\newcommand{\op}{\operatorname}
\newcommand{\R}{\mathbb{R}}
\newcommand{\C}{{\cal{C}}}
\newcommand{\N}{\mathbb{N}}
\newcommand{\ep}{\varepsilon}
\newcommand{\overbar}[1]{\mkern 1.5mu\overline{\mkern-1.5mu#1\mkern-1.5mu}\mkern 1.5mu}

\renewcommand{\(}{\left(}
\renewcommand{\)}{\right)}

\begin{document}

\title{On the locus of multiple maximizing geodesics on a globally hyperbolic spacetime}
\author{{\sc Alec Metsch$^{1}$} \\[2ex]
      $^{1}$ Universit\"at zu K\"oln, Institut f\"ur Mathematik, Weyertal 86-90, \\
      D\,-\,50931 K\"oln, Germany \\
      email: ametsch@math.uni-koeln.de \\[1ex]
      {\bf Key words:} 
      Lorentzian geometry, Lorentzian weak KAM, \\ 
     singularities of Lorentzian distance
      }
      

\maketitle

\begin{abstract}
\noindent
Extending the work of Cannarsa, Cheng and Fathi \cite{Cannarsa/Cheng/Fathi}, we investigate topological properties of the locus ${\cal NU}(M,g)$ of multiple maximizing geodesics on a globally hyperbolic spacetime $(M,g)$, i.e.\ the set of causally related pairs  $(x,y)$ for which there exists more than one maximizing geodesic (up to reparametrization) from $x$ to $y$. We will prove that this set is locally contractible. We will also define the notion of a Lorentzian Aubry set ${\cal A}$ and prove that the inclusions ${\cal NU}(M,g)\hookrightarrow \op{Cut}_M\hookrightarrow J^+\backslash {\cal A}$ are homotopy equivalences.
\end{abstract}

\section{Introduction}

In the recent work \cite{Cannarsa/Cheng/Fathi}, the authors Cannarsa, Cheng and Fathi established the following result:
\begin{theorem*}
    Consider the closed subset $C\subseteq N$ of the complete and connected Riemannian manifold $(N,h)$. Then the set of singularities $x\in N\backslash C$ of the function $d_C(x):=\inf_{c\in C}d(c,x)$ is locally contractible.
\end{theorem*}
Let us recall the definition of local contractibility.
\begin{definition}\rm\label{defa}
    Let $X$ be a topological space. A subset $A\subseteq X$ is called \emph{locally contractible} if for each $x_0\in A$ and each open neighbourhood $V$ of $x_0$, there exists an open neighbourhood $U$ of $x_0$ and a homotopy $h:(U\cap A)\times [0,1]\to V\cap A$ such that $h(x,0)=x$ for all $x\in U\cap A$ and $h(x,1)=h(y,1)$ for all $x,y\in U\cap A$.
\end{definition}
Let $(N,h)$ be as in the theorem. Following \cite{Cannarsa/Cheng/Fathi}, let us denote by ${\cal NU}(N,h)$ the set of pairs $(x,y)\in N\times N$ for which there exist two distinct minimizing geodesics connecting $x$ to $y$. It is well-known (\cite{FathiHJ}, Corollary 4.24) that $(x,y)\in {\cal NU}(N,h)$ if and only if $x\neq y$ and the distance function $d$ is not differentiable at $(x,y)$. This is also equivalent to $x\neq y$ and the non-differentiability of $d(x,\cdot)$ at $y$. 

Applying the above theorem to a singleton $C=\{x\}$, $x\in N$, and the diagonal $C=\Delta:=\{(x,x)\mid x\in N\}$ in the product manifold $N\times N$, the following result is an easy consequence of the above theorem (see also Theorem 1.3 in \cite{Cannarsa/Cheng/Fathi}):
\begin{thmx}\label{theorema}
Let $(N,h)$ be a complete and connected Riemannian manifold.
\begin{enumerate}[(a)]
    \item The set ${\cal NU}(N,h)$ is locally contractible. 
    \item For any $x\in N$, the set $\{y\in N\mid (x,y)\in {\cal NU}(N,h)\}$ is locally contractible.
\end{enumerate}
\end{thmx}

Using similar methods, the authors also derived a global topological result. To state it, we introduce the following notion.

\begin{definition*}
    For a closed set $C\subseteq N$ of a complete and connected Riemannian manifold $(N,h)$, the Aubry set ${\cal A}(C)$ is defined as the set of points $x\in N$ such that there exists a geodesic $\gamma:[0,\infty) \to N$, parametrized by arc length, with $\gamma(t_0)=x$ for some $t_0>0$ and $d_C(\gamma(t))=t$ for all $t\geq 0$. 
\end{definition*}

Note that the corresponding Definition 1.4 in \cite{Cannarsa/Cheng/Fathi} of the Aubry set is slightly different, as only points $x\in N\backslash C$ are considered. However, due to Theorem 1.6 in \cite{Cannarsa/Cheng/Fathi}, this is only a matter of taste. In our notation, this theorem reads as follows:
\begin{theorem*}
    If $C$ is a closed subset of the complete and connected Riemannian manifold $(N,h)$, the inclusion $\{x\in N\backslash C\mid d_C \text{ is not diff.\ at } x\}\hookrightarrow N\backslash {\cal A}(C)$ is a homotopy equivalence.
\end{theorem*}

Applying this result again to the special cases $C=\{x\}$ and $C=\Delta$, one obtains:

\begin{thmx}\label{theoremb}
    Let $(N,h)$ be a complete and connected Riemannian manifold. 
\begin{enumerate}[(a)]
    \item The inclusion ${\cal NU}(N,h)\hookrightarrow (N\times N)\backslash {\cal A}(\Delta)$ is a homotopy equivalence.
    \item For any $x\in N$, the inclusion 
    \[
    \{y\in N\mid (x,y)\in {\cal NU}(N,h)\}\hookrightarrow N\backslash {\cal A}(\{x\})
    \]
    is a homotopy equivalence.
\end{enumerate}
\end{thmx}

The aim of the present work is to extend Theorem \ref{theorema} and Theorem \ref{theoremb} to the setting of a globally hyperbolic spacetime. If $(M,g)$ is a globally hyperbolic spacetime, let us denote by ${\cal NU}(M,g)$ (resp.\ ${\cal NU}^t(M,g)$) the set of causally (resp.\ chronologically) related points for which there exist two distinct (up to reparametrization) maximizing geodesics connecting them. Our first main result is the extension of Theorem \ref{theorema}.
\begin{theorem}\label{aa}
Let $(M,g)$ be a globally hyperbolic spacetime.
\begin{enumerate}[(a)]
    \item The set ${\cal NU}(M,g)$ is locally contractible. 
    \item For any $x\in M$, the set $\{y\in M\mid (x,y)\in {\cal NU}(M,g)\}$ is locally contractible. 
    \item[(b')] For any $y\in M$, the set $\{x\in M\mid (x,y)\in {\cal NU}(M,g)\}$ is locally contractible. 
\end{enumerate}
\end{theorem}
Note that (b') follows from (b) by reversing the time orientation on $M\!$.
\medskip 

In contrast to the Riemannian case, parts (a) and (b) cannot be deduced from a more general result, since the product of two spacetimes, equipped with the product metric, is not a spacetime. 

To extend Theorem \ref{theoremb}, we will introduce the Lorentzian Aubry set(s) in a way that is the natural extension of the Riemannian case. 

In this paper, $J^+(x)$ (resp.\ $I^+(x)$) denotes the causal (resp.\ chronological) future of $x$, while $J^+$ (resp.\ $I^+$) denotes the set of points $(x,y)$ with $y\in J^+(x)$ (resp.\ $y\in I^+(x)$). 

\begin{definition}\rm
    Given $x\in M$, we define the \emph{future Aubry set} ${\cal A}(x)\subseteq M$ as the set of all points $y\in J^+(x)$ such that there exists a future ray through $y$ emerging from $x$, i.e.\ a future inextendible maximizing\footnote{For rays and lines (defined next), maximizing of course means that any restriction to a compact interval is maximizing.} geodesic $\gamma:[0,a)\to M$, $a\in (0,\infty]$, with $\gamma(0)=x$ and $\gamma(t)=y$ for some $t\in (0,a)$. We define the \emph{Aubry set} ${\cal A}\subseteq M\times M$ as the set of points $(x,y)\in J^+$ such that there exists a line through $x$ and $y$, i.e.\ a future and past inextendible maximizing geodesic $\gamma:I\to M$, with $\gamma(t_1)=x$ and $\gamma(t_2)=y$ for some $t_1,t_2\in I$, $t_1<t_2$. 
\end{definition}

\begin{remark}\rm\label{sdsxdfgh}
    Another suitable definition of the Aubry set ${\cal A}$,  which also extends the Riemannian definition ${\cal A}(\Delta)$ to the Lorenzian setting, is to consider only future or past inextendible maximizing geodesics $\gamma:(-a,a)\to M$, $a\in (0,\infty]$, such that there exists $t\in (0,a)$ with $\gamma(t)=y$ and $\gamma(-t)=x$. We denote this set by $\tilde {\cal A}$.

    In the Riemannian case, the geodesic flow is complete, so both definitions are equivalent. In the Lorentzian case, however, geodesics can be incomplete, and the two definitions may differ.  Still, ${\cal A}\subseteq \tilde {\cal A}$. The following theorem remains valid for $\tilde {\cal A}$, and the proof becomes even simpler in that case. However, our version based on lines is more natural and interesting in the Lorentzian context. Note that, under the assumption of causal geodesic completeness, the definitions become equivalent.
\end{remark}
Our second main result is the extension of Theorem \ref{theoremb}.

\begin{theorem}\label{b}
Let $(M,g)$ be a globally hyperbolic spacetime.
\begin{enumerate}[(a)]
    \item The inclusions
        ${\cal NU}(M,g)\hookrightarrow J^+\backslash {\cal A}$
    and ${\cal NU}^t(M,g)\hookrightarrow I^+\backslash {\cal A}$ are homotopy equivalences.
    \item For any $x\in M$, the inclusions 
        \[
        \{y\in M\mid (x,y)\in {\cal NU}(M,g)\}\hookrightarrow J^+(x)\backslash {\cal A}(x)
        \]
    and
    \[
    \{y\in M\mid (x,y)\in {\cal NU}^t(M,g)\}\hookrightarrow I^+(x)\backslash {\cal A}(x)
    \]
    are homotopy equivalences.
\end{enumerate}
\end{theorem}
Obviously, by reversing the time orientability of $M$ as above, one also obtains a corresponding result (b') for the "past" Aubry set.
\medskip 

In \cite{Cannarsa/Cheng/Fathi}, the authors also investigated the topological structure of singularities of continuous viscosity solutions to the evolutionary Hamilton-Jacobi equation
\begin{align*}
    \partial_tU+H(x,d_xU)=0
\end{align*}
    for a general Tonelli Hamiltonian $H$ (see Theorem 1.8 and 1.10 in \cite{Cannarsa/Cheng/Fathi}). This framework includes, in particular, the function $U:(0,\infty)\times N\to \R, \ (t,x)\mapsto \frac{d_C(x)^2}{2t},$ on a complete and connected Riemannian manifold $(N,h)$, as it is a viscosity solution of the evolutionary Hamilton-Jacobi equation for the Hamiltonian
\[
H:T^*N\to \R,\ H(x,p):=\frac 12 |p|_h^2.
\]
Here, $|p|_h$ denotes the dual norm of $p\in T_x^*N$. Except for some refinements in the arguments for the homotopy equivalence, the results in the Riemannian case can be proved very similarly to those general ones concerning continuous viscosity solutions of the evolutionary Hamilton-Jacobi equation.
The proofs rely on the representation formula of continuous viscosity solutions via the Lax-Oleinik evolution of some lower semicontinuous function $u$ (cf.\ Theorem 1.2 in \cite{FathiHJ}). For the Lax-Oleinik evolution, it is then well-known that, locally, we find $s\ll t$ such that the Lasry-Lions-type regularization $\hat T_sT_tu$ is $C^1$. Here, $\hat T_s$ and $T_t$ denote the backward and forward Lax-Oleinik semigroup, respectively, i.e.\ 
\begin{align*}
    \hat T_tf(x):=\sup_{y\in N}\{f(y)-h_t(x,y)\}  \text{ and } T_tf(y):=\inf_{x\in N}\{f(x)+h_t(x,y)\},
\end{align*}
with $h_t$ being the minimal action to go from $x$ to $y$ in time $t$ for the Lagrangian associated to $H$, i.e.\ 
\begin{align*}
    L:TN\to \R,\ L(x,v):=\frac 12 |v|_h^2.
\end{align*}
The differentiability of the function $\hat T_sT_tu$ plays the key role in the proof of the homotopy properties (see Claim 4.9 in \cite{Cannarsa/Cheng/Fathi}).
\medskip

In the Lorentzian setting, we investigate the Hamiltonian 
\begin{align}
  H:T^*M\to \R\cup\{+\infty\},\ H(x,p)=
  \begin{cases}
      \frac{1}{4|p|_g},\ &p\in \op{int}(\C_x^*),
      \\[10pt]
      +\infty,\ &\text{otherwise,} 
  \end{cases}
  \label{Ham}
\end{align}
where $\C_x^*\subseteq T_x^*M$ denotes the cone of future-directed causal covectors, cf.\ Definition \ref{causalco}. It lacks the regularity and superlinearity properties required by the classical theory of Tonelli Hamiltonians, on which the proof of the differentiability of $\hat T_sT_tu$ heavily relies. This prevents us from establishing general regularity results for $\hat T_sT_tu$, where $\hat T_s$ and $T_t$ now denote the Lax-Oleinik semigroups w.r.t.\ to the Hamiltonian \eqref{Ham} and its associated Lagrangian \eqref{Lag}; that is, $h_t$ is being replaced by $c_t$, the minimal (Lorentzian) action to go from $x$ to $y$ in time $t$ for the Lagrangian
\begin{align}
     L:TM\to \R\cup\{+\infty\},\ L(x,v):=
    \begin{cases}
        -|v|^{\frac 12}_g, & \text{ if } v\in \C_x,
        \\[10pt]
        +\infty,& \text{otherwise,}
    \end{cases}
    \label{Lag}
\end{align}
$\C_x$ denoting the cone of future-directed causal vectors.
 \medskip

The core of this paper is therefore devoted to the study of regularity properties of $\hat T_sT_tu$, where $u$ takes the specific form $u=\chi_x$, see \eqref{chi}. We establish local $C^1$-regularity on $I^+(x)$ for $s\ll t$, and also examine the regularity w.r.t.\ $s,t$ and $x$.

As in the Riemannian case, this $C^1$-regularity follows from proving that $\hat T_sT_t\chi_x$ is both locally semiconvex and semiconcave. In the Riemannian setting, local semiconvexity holds for all $s<t$, as shown in \cite{FathiHJ}, Theorem 6.2. The proof relies on the  superlinearity of Tonelli Lagrangians. A slight modification of the argument shows that for any fixed $t>0$ and any point $y_0\in M$, the (nearly) optimal points $z$ in the definition of $\hat T_sT_t\chi_x$ remain uniformly bounded as $s$ is small and $y$ stays close to $y_0$. More precisely, for some uniform constant $C>0$, 
\begin{align*}
    \hat T_sT_t\chi_x(y)>\sup\{T_t\chi_x(z)-h_s(y,z)\mid d_h(y,z)> Cs\}, 
\end{align*}
and consequently, 
\begin{align}
    \hat T_sT_t\chi_x(y)=\sup\{T_t\chi_x(z)-h_s(y,z)\mid d_h(y,z)\leq Cs\}. \label{eqyy}
\end{align}
 The idea here is that $h_s(y,z)$ dominates $T_t\chi_x(z)$ and it gets too large when $s\to 0$ and $d_h(y,z)$ tends slower to $0$ than $s$. Due to the local semiconcavity of $h_s$ on $M\times M$, the representation \eqref{eqyy} is sufficient to ensure local semiconvexity of $\hat T_sT_t\chi_x$ near $y_0$, as shown in \cite{Fathi/Figalli}. 
 
 In the Lorentzian context, fixing a complete Riemannian metric $h$ on $M$ for reference, things become more subtle, as superlinearity and regularity of the Lagrangian fail. To get a similar control over the (nearly) optimal points $z$, we instead show that $T_t\chi_x(z)$ dominates $c_s(y,z)$ when $s\to 0$ and $d_h(y,z)$ tends slower to $0$ than $\sqrt{s}$; in the sense that
 \[
 T_t\chi_x(z)-c_s(y,z) > T_t\chi_x(y)-c_s(y,y)+O(\sqrt{s}).
 \]
This results in a analogous bound \eqref{eqyy}, with $s$ replaced by $\sqrt{s}$. However, unlike in the Riemannian setting, this condition alone is not sufficient to guarantee local semiconvexity. This is due to the fact that, in contrast to $h_s$, the minimal Lorentzian action $c_s$ is not locally semiconcave on all of $M\times M$, but only on the set of chronologically related points $I^+$ \cite{McCann2}. As a result, we must also ensure that the optimal points $z$ stay uniformly bounded away from the boundary of the causal future. These subtleties introduce additional difficulties in establishing local semiconvexity, see Section 3. 

To prove local semiconcavity, we use an approach similar to the Riemannian case: we approximate the locally semiconcave function $T_t\chi_x$ from above by a family $f_i$ of smooth functions that is uniformly locally semiconcave (cf.\ \cite{Bernard}). One then shows that, for small $s>0$, the family $(\hat T_sf_i)_i$ remains uniformly locally semiconcave and approximates $\hat T_sT_t\chi_x$. This kind of approximation result is known in the Riemannian setting for complete manifolds \cite{Fathi/Figalli/Rifford}. In the Lorentzian case, however, similar difficulties as for semiconvexity appear. In particular, the construction of the family $f_i$ requires additional care to ensure the propagation of both the approximation property and the regularity (Section 4). 

Theorem \ref{main2} collects the main results from Sections 3 and 4 (in particular, its proof shows the local $C^1$-regularity of $\hat T_sT_t\chi_x$) and implies Theorem \ref{ab} (compare to the important Lemma 3.6 in \cite{Cannarsa/Cheng/Fathi}), which in turn leads to the key Corollary \ref{cor4} – a result that is not needed in \cite{Cannarsa/Cheng/Fathi}. These results play a central role in the proofs of our main theorems. 

For instance, much like in the Riemannian setting (cf.\ the proof of Proposition 3.1 in \cite{Cannarsa/Cheng/Fathi}), Theorem \ref{ab} allows us to show local contractibility of the set ${\cal NU}^t(M,g)={\cal NU}(M,g)\cap I^+$, as well as its analogue for a fixed $x\in M$. The argument is carried out in Subsection \ref{ppoiouhgs}, where we also show how the proof can be adapted to include the case of lightlike geodesics. This issue, of course, does not arise in the Riemannian setting and stems from from the fact that the conclusion from Theorem \ref{ab} holds only for chronologically related points. As a consequence, it not only makes Corollary \ref{cor4} necessary, but also requires us – independently of the theory developed in Sections 3 and 4 – to establish the local contractibility of $\operatorname{Cut}_M$ (and $\operatorname{Cut}_M(x)$), with the homotopy required to satisfy certain additional conditions (cf.\ Lemma \ref{lo}).

In Subsections \ref{ppoiouhgs1} and \ref{ppoiouhgs2}, we prove Theorem \ref{b}. Note that the proof of Theorem \ref{theoremb} from \cite{Cannarsa/Cheng/Fathi} does not carry over to our setting, as there is no analogue to Proposition 7.1 from that work. Nevertheless, Corollary \ref{cor4} still implies that the inclusion ${\cal NU}^t(M,g)\hookrightarrow \op{Cut}_M^t$ is a homotopy equivalence. This also holds for the versions with fixed $x$, and can be generalized to include lightlike geodesics.

To complete the proof of Theorem \ref{b}, we combine this with the fact that $\op{Cut}_M^t$ (resp.\ $\op{Cut}_M$) is a strong deformation retract of $I^+\backslash {\cal A}$ (resp.\ $J^+\backslash {\cal A}$), including the corresponding versions for fixed $x$. These latter results do not rely on the theory developed in Sections 3 and 4 - hence neither on Theorem \ref{ab} nor on Corollary \ref{cor4} - but rather on classical results about maximizing geodesics in globally hyperbolic spacetimes (in particular, Theorem \ref{thme} and Lemma \ref{lm}). The main difficulty lies in proving that $\op{Cut}_M$ is a strong deformation retraction of $J^+\backslash {\cal A}$ (the version for fixed $x$ is considerably simpler). The geometric intuition is to move the two points $(x,y)\in J^+\backslash {\cal A}$ in opposite directions along the maximizing geodesic connecting them, which, by definition, is not globally maximizing. The challenge arises since the geodesic flow may be incomplete, and we lack information about where exactly the geodesic ceases to be maximizing. Therefore, the points $x$ and $y$ must be transported with individual chosen speeds that depend continuously on $(x,y)$ in order to construct a valid homotopy. 

Actually, this is precisely the reason why this difficulty does not appear when working with the set $\tilde {\cal A}$ instead, as introduced in Remark \ref{sdsxdfgh}. In particular, this problem also does not arise when $(M,g)$ is assumed to be causal geodesically complete, and likewise it's not an issue in the Riemannian setting. In fact, our approach to Theorem \ref{b} carries over to the Riemannian case with far fewer complications as here and offers an alternative (and, in my view, simpler) proof of Theorem 1.6 in \cite{Cannarsa/Cheng/Fathi}.

\section{The Lagrangian}
In this chapter, and throughout the following ones, let $(M,g)$ be a globally hyperbolic spacetime, where the metric $g$ is taken to have signature $(-,+,...,+)$. We refer to future-directed causal (resp.\ timelike) vectors simply as causal (resp.\ timelike), and explicitly specify past-directed causal when needed. We understand $0$ to be a causal vector. A curve is always assumed to be piecewise smooth if not otherwise said. In particular, a curve is referred to as causal (timelike) if it is piecewise smooth and future-directed causal (timelike). The Lorentzian distance function is denoted by $d$.
For $x\in M$, we denote by $\C_x\subseteq T_xM$ the cone of causal vectors. Note that $\C_x$ is closed. We also set $\C:=\{(x,v)\in TM\mid v\in \C_x\}$. We can equip $M$ with a complete Riemannian metric, which will be fixed and denoted by $h$. All balls $B_r(x)$, $x\in M$, $r>0$, are understood to be taken w.r.t.\ the metric $h$.

\begin{definition}\rm
    We define the Lagrangian $L:TM\to \R\cup\{+\infty\}$ as
    \begin{align*}
        L(x,v):=
        \begin{cases}
          -|v|^{\frac 12}_g, & \text{ if } v\in \C_x,
          \\[10pt]
          +\infty,& \text{otherwise.}
        \end{cases}
    \end{align*}
    Here, $|v|_g:=\sqrt{|g_x(v,v)|}$. See also \cite{Mondino/Suhr}, Section 2, and \cite{McCann2}, Section 3.
\end{definition}

\begin{definition}\rm
\begin{enumerate}[(a)]
    \item 
    The \emph{action} of a curve $\gamma:[a,b]\to M$ is defined by
    \begin{align*}
        \mathbb{L}(\gamma):=\int_a^b L(\gamma(t),\dot \gamma(t))\, dt\in (-\infty,\infty].
    \end{align*}
    Note that $\mathbb{L}(\gamma)$ is finite if and only if $\gamma$ is causal.
    \item A curve $\gamma:[a,b]\to M$ is called an $L$-\emph{minimizer} if for any other curve $\tilde \gamma:[a,b]\to M$ with the same endpoints, we have $\mathbb{L}(\gamma)\leq \mathbb{L}(\tilde \gamma)$.
\end{enumerate}
\end{definition}

\begin{notation} \rm
    We reserve the term \emph{maximize} to refer specifically to the Lorentzian length functional. That is, a maximizing curve $\gamma:[a,b]\to M$ is a causal curve that satisfies $\ell_g(\gamma):=\int_a^b |\dot \gamma(t)|_g\, dt=d(\gamma(a),\gamma(b))$. The term maximal, on the other hand, is reserved for geodesics and refers to a geodesic defined on its maximal existence interval. In particular, a maximal causal geodesic is not necessarily maximizing. The following result is well-known.
\end{notation}

\begin{theorem}\label{asdfghj}
    For any two points $(x,y)\in J^+$, there exists a maximizing geodesic connecting $x$ to $y$. Moreover, any maximizing curve must be a pregeodesic.
\end{theorem}
\begin{proof}
    See \cite{ONeill}, Chapter 14, Proposition 19 and \cite{Minguzzi}, Theorem 2.9.
\end{proof}

\begin{lemma}\label{lt}
    Let $x\in M$, $y\in J^+(x)$ and $t>0$. Let $\gamma:[0,t]\to M$ be a curve connecting $x$ to $y$. Then:
    \begin{align*}
        \gamma \text{  is a maximizing geodesic} \Rightarrow \gamma \text{ is $L$-minimizing }  
    \end{align*}
    If $y\in I^+(x)$, the implication becomes an equivalence. 
\end{lemma}
\begin{proof}
    This follows immediately from the Cauchy-Schwarz inequality. Indeed, for any  causal curve $\gamma:[0,t]\to M$, we have
    \begin{align*}
        \mathbb{L}(\gamma) \geq -t^\frac 12 \ell_g(\gamma)^\frac 12
    \end{align*}
    with equality if and only if $|\dot \gamma|_g$ is constant. This easily implies that a maximizing geodesic is $L$-minimizing. Conversely, let $y\in I^+(x)$ and suppose $\gamma$ is $L$-minimizing. Being $(M,g)$ globally hyperbolic, the above theorem guarantees the existence of a maximizing geodesic $\tilde \gamma:[0,t]\to M$ connecting $x$ to $y$. Then
    \begin{align*}
        \ell_g(\gamma)^\frac 12 \geq -t^{-\frac 12} \mathbb{L}(\gamma) \geq -t^{-\frac 12} \mathbb{L}(\tilde \gamma)= \ell_g(\tilde \gamma)^\frac 12 \geq \ell_g(\gamma)^\frac 12.
    \end{align*}
    Thus, $\gamma$ is maximizing as well, and equality must hold in each of the above steps, implying that $|\dot \gamma|_g$ is constant. Since $y\in I^+(x)$, we must have $|\dot \gamma|_g=cons.\neq 0$. Combining with the fact that $\gamma$ is a pregeodesic by the theorem above, we conclude that $\gamma$ is in fact a geodesic.
\end{proof}

\begin{remark}\rm
    It is not difficult to verify that the second derivative along the fibres, $\frac{\partial^2 L}{\partial v^2}$, is positive definite at every point $(x,v)\in \op{int}(\C)$ (\cite{Mondino/Suhr}, Lemma 2.1). As a consequence, it is well-known that there exists a smooth local flow $\phi_t$ on $\op{int}(\C)$ whose orbits are precisely the speed curves of extremals for $L$; that is, the curves of the form $(\gamma(t),\dot \gamma(t))$, where $\gamma:I\to M$ is a $C^2$-curve satisfying, in local coordinates, the Euler-Lagrange equation
    \begin{align*}
        \frac{d}{dt}\bigg(\frac{\partial L}{\partial v}(\gamma(t),\dot \gamma(t))\bigg)=\frac{\partial L}{\partial x}(\gamma(t),\dot \gamma(t)).
    \end{align*}
    Moreover, every timelike $L$-minimizing curve $\gamma:[a,b]\to M$ solves the Euler-Lagrange equation. Since the timelike $L$-minimizing curves are exactly the timelike maximizing geodesics, and since every timelike geodesic is locally maximizing (Proposition 7.3 in \cite{Metsch1}), it follows that the Euler-Lagrange flow coincides with the geodesic flow restricted to the invariant set $\op{int}(\C)$.
\end{remark}

\begin{definition}\rm
The \emph{minimal time-$t$-action} to go from $x$ to $y$ is defined by
\begin{align*}
    c_t(x,y) :=\inf\{\mathbb{L}(\gamma)\mid \gamma:[0,t] \to M\text{ is a curve connecting } x \text{ to } y\},
\end{align*}
where $\inf(\emptyset):=\infty$.
\end{definition}

\begin{corollary}\label{a}
    We have the following identity:
    \begin{align*}
        c_t(x,y)=
        \begin{cases}
        0, &\text{ if } t=0 \text{ and } x=y,
        \\[10pt]
            -t^{\frac 12} d(x,y)^{\frac 12}, &\text{ if } t>0 \text{ and }y\in J^+(x),
            \\[10pt]
            +\infty, & \text{ otherwise.}
        \end{cases}
    \end{align*}
    Moreover, for any $x\in M$, $y\in J^+(x)$ and $t>0$, there exists a smooth $L$-minimizing geodesic $\gamma:[0,t]\to M$ connecting $x$ to $y$.
\end{corollary}
\begin{proof}
    This follows immediately from Theorem \ref{asdfghj} and the proof of Lemma \ref{lt}.
\end{proof}

\begin{definition}\rm\label{causalco}
    For each $x\in M$, we denote the canonical isomorphism by
    \begin{align*}
        T_xM\to T_x^*M,\ v\mapsto v^\flat:=g(v,\cdot).
    \end{align*}
    We define $\C_x^*$ as the image of $\C_x$ under this isomorphism, and $\C^*:=\{(x,p)\in T^*M\mid p\in \C_x^*\}.$
\end{definition}

In the following lemma, we state semiconcavity and related properties for the time-$t$-action. For defintions, see Appendix A in \cite{Fathi/Figalli}.

\begin{lemma}\label{l4}
\begin{enumerate}[(a)]
    \item The function 
    \begin{align*}
        \C:(0,\infty)\times M\times M\to \R\cup\{+\infty\},\ (t,x,y)\mapsto c_t(x,y),
    \end{align*}
    is real-valued and continuous on $(0,\infty)\times J^+$, and locally semiconcave on $(0,\infty)\times I^+$.
    \item If $x\in M$, $y\in I^+(x)$ and $t>0$, then the set of super-differentials of $\C$ at the point $(t,x,y)$ is given by
    \begin{align*}
        \partial^+ \C(t,x,y)
        =
        \op{conv}\left(\bigg\{\left(\partial_t c_t(x,y),-\frac{\partial L}{\partial v}(x,\dot \gamma(0)), \frac{\partial L}{\partial v}(y,\dot \gamma(t))\right)\bigg\}\right),
    \end{align*}
    where the set runs over all maximizing geodesics $\gamma:[0,t]\to M$ connecting $x$ to $y$.

    In particular, $\C$ is differentiable at $(t,x,y)$ if and only if there is a unique maximizing geodesic connecting $x$ to $y$ in time $t$ (equivalently, in time $1$).
\end{enumerate}
\end{lemma}
\begin{proof}
Part (a) follows from the well-known continuity (\cite{ONeill}, Chapter 14, Lemma 21) and semiconvexity (\cite{McCann2}, Proposition 3.4) properties of the Lorentzian distance function on $J^+$ and $I^+$, respectively, combined with standard properties of semiconvex/semiconcave functions (\cite{Cannarsa}, Proposition 2.1.12). For part (b), note that the super-differential can be computed in terms of the sub-differential of $d$ (\cite{McCann}, Lemma 5).
For the convenience of the reader, we provide a proof in the appendix.
\end{proof}

\begin{definition}\rm
    The \emph{Legendre transform} of $L$ is the map
    \begin{align*}
        {\cal L}:\op{int}(\C)\to T^*M,\ 
        (x,v)\mapsto \left(x,\frac{\partial L}{\partial v}(x,v)\right).
    \end{align*}
\end{definition}

\begin{proposition}
    The Legendre transform is a diffeomorphism onto its image
        $\op{im}({\cal L})=\op{int}(\C^{*})$.
\end{proposition}
\begin{proof}
We have
    \begin{align}
        \bigg(x,\frac{\partial L}{\partial v}(x,v)\bigg)
        =
        \bigg(x,\frac 12 |v|_g^{-\frac 32}v^\flat\bigg),
        \label{eqk}
    \end{align}
    which is clearly a diffeomorphism from $\op{int}(\C)$ to $\op{int}(\C^*)$.
\end{proof}

\begin{corollary}\label{lc}
    Let $x\in M$ and $y\in I^+(x)$. Then the following are equivalent:
    \begin{align*}
        (x,y)\in \op{sing}(d) 
        \Leftrightarrow 
        y\in \op{sing}(d(x,\cdot))
        \Leftrightarrow
        (x,y)\in {\cal NU}(M,g)
    \end{align*}
    Here, $\op{sing}(f)$ denotes the set of points where a function $f:M\to \overbar \R$ fails to be differentiable.
\end{corollary}
\begin{proof}
    Since $d$ is positive on $I^+$, it suffices to prove the lemma with $c_1=-d^\frac 12$ instead of $d$. It is clear that $y\in \op{sing}(c_1(x,\cdot))\Rightarrow (x,y)\in \op{sing}(c_1)$. 
    
    To prove $(x,y)\in \op{sing}(c_1)
        \Rightarrow
        (x,y)\in {\cal NU}(M,g)$, we argue by contraposition. 
    Since unique super-differentiability of locally semiconcave functions is known to imply differentiability (\cite{Villani}, Theorem 10.8), Lemma \ref{l4} implies that if there exists a unique maximizing geodesic connecting $x$ to $y$, then $(x,y)\notin \op{sing}(c_1)$.
    
    Since differentiability implies unique super-diferentiability, Lemma \ref{l4} and the fact that the Legendre transform is a diffeomorphism imply that if there exist two distinct maximizing geodesics connecting $x$ to $y$, then the super-differential $\partial^+(c_1(x,\cdot))(y)$ is not reduced to a singleton. Hence, $y\in \op{sing}(d(x,\cdot))$.
\end{proof}

\begin{definition}\rm
The \emph{forward Lax-Oleinik semigroup} is the family of maps $(T_t)_{t\geq 0}$, defined on the space of functions $f:M\to \overbar \R$, given by
\begin{align*}
    T_tf:M\to \overbar \R,\ T_tf(y):=\inf\{f(x)+c_t(x,y)\mid x\in M\}.
\end{align*}
Here, we use the convention $-\infty+\infty:=\infty$.

The \emph{backward Lax-Oleinik semigroup} is the family of maps $(\hat T_t)_{t\geq 0}$, defined on the space of functions $f:M\to \overbar \R$, given by
\begin{align*}
    \hat T_tf:M\to \overbar \R,\ \hat T_tf(x):=\sup\{f(y)-c_t(x,y)\mid y\in M\}.
\end{align*}
Here, we use the convention $+\infty-\infty:=-\infty$.
\end{definition}

\begin{lemma}\label{lf}
\begin{enumerate}[(a)]
\item 
 Let $f:M\to \overbar \R$ be any function, and let $x\in M$, $y\in I^+(x)$ and $t>0$. Suppose that $f(x)$ and $T_{t}f(y)$ are finite. Additionally, assume that 
    $T_{t}f(y)=f(x)+c_{t}(x,y)$.
    Then the function $Tf$ is super-differentiable at $(t,y)$, and
    \begin{align}
    \left(\partial_t c_{t}(x,y),\frac{\partial L}{\partial v}(y,\dot \gamma(t))\right) \in \partial^+ Tf(t,y), \label{zz}
    \end{align}
    where $\gamma:[0,t]\to M$ is a maximizing geodesic connecting $x$ to $y$.
    Moreover, 
    \begin{align}
    \frac{\partial L}{\partial v}(x,\dot \gamma(0))\in \partial^- f(x). \label{zzz}
    \end{align}

    \item Let $f:M\to \overbar \R$ be any function, and let $x\in M$, $y\in I^+(x)$ and $t>0$. Suppose that $f(y)$ and $\hat T_{t}f(x)$ are finite. Additionally, assume that 
    $\hat T_{t}f(x)=f(y)-c_{t}(x,y)$.
    Then the function $\hat Tf$ is sub-differentiable at $(t,x)$, and
    \[
    \left(-\partial_t c_{t}(x,y),\frac{\partial L}{\partial v}(x,\dot \gamma(0))\right)\in \partial^- \hat Tf(t,x),
    \]
    where $\gamma:[0,t]\to M$ is a maximizing geodesic connecting $x$ to $y$.
     Moreover,
    \[
    \frac{\partial L}{\partial v}(y,\dot \gamma(t))\in \partial^+ f(y).
    \]
    \end{enumerate}
\end{lemma}
\begin{proof}
    This is a simple consequence of Lemma \ref{l4}.
\end{proof}

\begin{definition}\rm
The \emph{Hamiltonian} associated with $L$ is the function
\begin{align*}
    &H:T^*M\to \R\cup\{+\infty\},\ 
    \\
    &H(x,p):=\sup\{pv-L(x,v)\mid v\in T_xM\}.
\end{align*}
\end{definition}

\begin{lemma}[See \cite{Mondino/Suhr}, Section 2, and \cite{McCann2}, Lemma 3.1]
    We have
    \begin{align*}
        H({\cal L}(x,v))=\frac{\partial L}{\partial v}(x,v)(v)-L(x,v)
        =-\frac 12 L(x,v)
    \end{align*}
    for all $(x,v)\in \op{int}(\C)$. 
\end{lemma}
\begin{proof}
    For $p\in \op{int}(\C_x^{*})$, observe that 
    \begin{align*}
        pv-L(x,v)\xrightarrow{|v|_h\to \infty } -\infty.
    \end{align*}
    Indeed, if $v\notin \C_x$, we have $pv-L(x,v)=-\infty$. Othwerwise, since $p\in \op{int}(\C_x^{*})$, we can define
    \begin{align*}
        &\alpha:=\sup\{pv\mid v\in \C_x,|v|_{h}=1\}<0, \text{ and }
        \\[10pt]
        &\beta:=\inf\{L(x,v)\mid v\in \C_x,|v|_{h}=1\}>-\infty.
    \end{align*}
    Therefore, on $\C_x\backslash \{0\}$,
    \begin{align*}
        pv-L(x,v)\leq \alpha |v|_{h}-\beta |v|_{h}^\frac 12 \xrightarrow{|v|_{h}\to \infty}-\infty.
    \end{align*}
    Thus, since continuous functions defined on compact sets attain their supremum, there is $v\in \C_x$ with
    \begin{align}
        H(x,p)=pv-L(x,v). \label{eqo}
    \end{align}
    We claim that $v\notin \partial \C_x$. First $v\neq 0$: If $v=0$, then we have, for any nonzero $w\in \C_x$,
    \begin{align*}
        p(\lambda w)-L(x,\lambda w)
        =
        \lambda pw+(\lambda|w|_g)^\frac 12
        > 0=pv-L(x,v)
    \end{align*}
    for sufficiently small $\lambda>0$, contradicting \eqref{eqo}.
    Now suppose $v\in \partial \C_x\backslash \{0\}$. Let $(e_0,...,e_n)$ be a generalized orthonormal frame in $T_{x}M$ with $e_0$ (future-directed) timelike. Then
    \begin{align*}
        v=\sum_{i=0}^n \lambda_i e_i, \text{ with } \lambda_0>0 \text{ and }\lambda_0^2-\sum_{i=1}^n \lambda_i^2=0.
    \end{align*}
    However, if we define, for small $\ep>0$, 
    \begin{align*}
        v(\ep):=(\lambda_0+\ep)e_0+\sum_{i=1}^n \lambda_i e_i,
    \end{align*}
    then $v(\ep)$ is causal and
    \begin{align*}
        |v(\ep)|_g^\frac 12= \left((\lambda_0+\ep)^2-\sum_{i=1}^n \lambda_i^2\right)^\frac 14 = (2\lambda_0 \ep+\ep^2)^\frac 14
        \geq (2\lambda_0)^\frac 14 \ep^\frac 14.
    \end{align*}
    Since $\lambda_0>0$, we conclude that
    \begin{align*}
       pv(\ep)-L(x,v(\ep))= pv(\ep)+|v(\ep)|_g^\frac 12
        \geq pv +\ep pe_0+(2\lambda_0)^\frac 14\ep^\frac 14>pv =pv-L(x,v).
    \end{align*}
    for small $\ep$, meaning that $v$ cannot be optimal in \eqref{eqo}.
   Hence, $v\in \op{int}(\C_x)$.
   
   Since $L$ is smooth on $\op{int}(\C)$, we can differentiate $w\mapsto pw-L(x,w)$ at its maximum point $w=v$ yielding
    \begin{align*}
        p=\frac{\partial L}{\partial v}(x,v).
    \end{align*}
    Therefore, using \eqref{eqk}, we get
    \begin{align*}
        H(x,p)=\frac{\partial L}{\partial v}(x,v)(v)-L(x,v)
        =
        \frac 12 |v|_g^\frac 12
        =
        -\frac 12 L(x,v).
    \end{align*}
\end{proof}

\begin{remark}\rm
\begin{enumerate}[(a)]
    \item The above lemma and the identity $p=\frac{\partial L}{\partial v}(x,v)$ imply that $H$ is given explicitly by \eqref{Ham}.
    \item From the above lemma, we conclude that $H$ is smooth on the open set $\op{int}(\C^{*})$, and satisfies the identity $H({\cal L}(x,v))={\cal L}(x,v)(v)-L(x,v)$. It is therefore well-known \cite{Fathi} that the Euler-Lagrange flow $\phi_t$ (i.e.\ the geodesic flow on $\op{int}(\C)$) is conjugate, via the diffeomorphism ${\cal L}:\op{int}(\C)\to \op{int}(\C^*)$, to the Hamiltonian flow $\psi_t$ on $\op{int}(\C^{*})$. The latter is understood with respect to its canonical symplectic structure as an open subset of the cotangent bundle.
\end{enumerate}

\end{remark}

\section{Local semiconvexity of $\hat T_sT_{t}\chi_x$}
For $x\in M$, we define the characteristic function
\begin{align}
    \chi_x:M\to \overbar \R,\ \chi_x(y):=
    \begin{cases}
        0,& \text{ if } x=y,
        \\[10pt]
        +\infty, &\text{ otherwise}.
    \end{cases} \label{chi}
\end{align}

Our main result in this section is the following theorem:

\begin{theorem}\label{thmb}
    Let $x_0\in M$, $y_0\in I^+(x_0)$ and $t_0>0$. Then there exist two open neighbourhoods $U$ and $V$ of $x_0$ and $y_0$, respectively, with $U\times V\subseteq I^+$, some number $s_0>0$, a constant $C_0>0$, and a non-decreasing sequence $(\delta_s)_{0< s \leq s_0}$ of positive numbers such that, for all $s\in (0,s_0]$, $t\in [t_0/2,3t_0/2]$ and $(x,y) \in U\times V$, we have
    \begin{align}
        \hat T_sT_{t}\chi_x(y)>\sup\{T_{t}\chi_x(z)-c_s(y,z)\mid d_h(y,z)\geq C_0\sqrt{s} \text{ or } d(y,z)\leq \delta_s\}. \label{eqae}
    \end{align}
    In particular, for all $s\in [0,s_0]$, $t\in [t_0/2,3t_0/2]$ and $(x,y) \in U\times V$, there exists $z\in M$ with $\hat T_sT_t\chi_x(y)=T_t\chi_x(z)-c_s(y,z)$. If, in addition, $s>0$, then necessarily $z\in I^+(y)$. 
    
    Moreover, the mapping $[0,s_0]\times [t_0/2,3t_0/2]\times U\times V\to \R,\ (s,t,x,y)\mapsto \hat T_sT_t\chi_x(y)$, is continuous, and for any $s_1\in (0,s_0]$, the family $\{\hat T_sT_t\chi_x\mid s\in [s_1,s_0],t\in [t_0/2,3t_0/2],x\in U\}$ is uniformly locally semiconvex on $V$.
\end{theorem}

We will prove this theorem in several steps.

\begin{lemma}\label{lh}
    For $x,y\in M$ and $t>0$, we have
    \begin{align*}
        T_{t}\chi_x(y)=c_t(x,y)=\C(t,x,y).
    \end{align*}
\end{lemma}
\begin{proof}
    This follows directly from the definition.
\end{proof}

\begin{lemma}\label{l1}
    Let $x_0\in M$, $y_0\in I^+(x_0)$ and $t_0>0$. Then there exist two open neighbourhoods $U$ and $V$ of $x_0$ and $y_0$, respectively, with $U\times V\subseteq I^+$, some number $s_0>0$, and a constant $C_0>0$ such that, for all $s\in (0,s_0]$, $t\in [t_0/2,3t_0/2]$ and $(x,y)\in U\times V$, we have
    \begin{align}
        \hat T_sT_{t}\chi_x(y)>\sup\{T_{t}\chi_x(z)-c_s(y,z)\mid
				d_{h}(y,z)\geq C_0\sqrt{s}\}. \label{eqaa}
    \end{align}
\end{lemma}
\begin{proof}
        The idea of proof comes from Lemma 6.3 in \cite{FathiHJ}.
        Let $r>0$ such that the compact set $K:=\overbar B_r(x_0)\times \overbar B_r(y_0)$ is contained in $I^+$. By Lemma \ref{lh} and Lemma \ref{l4}\footnote{To apply Lemma \ref{l4}, we use the fact that whenever $(t,x,y)\in \op{int}(\C)$, then $\partial^+ \C(t,x,\cdot)(y)=\pi_3\circ \partial^+ \C(t,x,y)$, where $\pi_3:T_t^*\R\times T_x^*M\times T_y^*M\to T_y^*M$ denotes the projection onto the third factor.} and standard compactness arguments together with the well-known continuity of the super-differential of locally semiconcave functions as multivalued maps, we find
        \begin{align*}
            \{(y,p)\in T^*M\mid p\in \partial^+T_t\chi_x(y),(x,y)\in K,t\in [t_0/4,3t_0/2]\}\Subset \op{int}(\C^{*}).
        \end{align*}
        Another compactness argument yields the existence of a constant $C>0$ such that
        \begin{align}
            \forall (t,x,y)\in [t_0/4,3t_0/2]&\times K,\  p\in \partial^+ T_t\chi_x(y): 
            pv\leq -C|v|_{h} \ \forall v\in \C_y. \label{eqp}
        \end{align}
        Moreover, since $\C$ is locally semiconcave on $(0,\infty)\times I^+$, Lemma \ref{lh} gives us a constant $M>0$ such that
        \begin{align}
            \op{Lip}(T_{\cdot}\chi_x(y)) \leq  M\text{ on } [t_0/4,3t_0/2]\label{eqs} 
        \end{align}
        for all $(x,y)\in K$, and also
        \begin{align}
            |c_t(y,z)|\leq \sqrt{t}M \text{ on } (0,\infty)\times ((\overbar B_r(y_0)\times \overbar B_r(y_0))\cap J^+). \label{eqw} 
        \end{align}

        Now, fix $r'\in (0,r)$ and define $U:=B_{r}(x_0)$ and $V:=B_{r'}(y_0)$. Choose $0<s_0\leq \min\{t_0/4,1\}$ with $4\frac{M}{C}\sqrt{s_0}\leq r-r'$.
        
        We claim that \eqref{eqaa} holds with $C_0:=4\frac{M}{C}$. To see this, let $s\in (0,s_0]$, $t\in [t_0/2,3t_0/2]$, $(x,y)\in U\times V$ be given, and suppose there exists a sequence $z_k\in M$ with $d_h(y,z_k)\geq C_0\sqrt{s}$ and
        \begin{align}
            T_{t}\chi_x(z_k)-c_s(y,z_k)
            \geq 
            \hat T_sT_{t}\chi_x(y)-\frac 1k
            \geq T_t\chi_x(y)-\frac 1k,\ k\in \N. \label{eqq}
        \end{align}
    
        Let $\gamma_k:[0,s]\to M$ be a maximizing geodesic connecting $y$ to $z_k$, and let $\tau_k\in [0,s]$ be the first time for which $\gamma_k(\tau_k)\in \partial B_{C_0\sqrt{s}}(y)\subseteq \overbar B_r(y_0)$. Set $\tilde z_k:=\gamma_k(\tau_k)$. Using the semigroup property and the fact that $\gamma_k$ is maximizing, we obtain
        \begin{align*}
            T_{t}\chi_{x}(z_k)-c_s(y, z_k)
            \leq
            T_{t-s+\tau_k}\chi_{x}(\tilde z_k)-c_{\tau_k}(y,\tilde z_k).
        \end{align*}
        Combining this with \eqref{eqq}, we get
        \begin{align*}
           &T_{t-s+\tau_k}\chi_{x}(\tilde z_k)-c_{\tau_k}(y,\tilde z_k)
           \geq T_{t}\chi_{x}(y)-\frac 1k
           \\
           \Leftrightarrow &T_{t-s+\tau_k}\chi_{x}(\tilde z_k)- T_{t}\chi_{x}(y)
           \geq 
           c_{\tau_k}(y,\tilde z_k)-\frac 1k.
           \end{align*}
           We now show that this is impossible. 
           
           Since $s\leq t_0/4$, it follows that $t-s+\tau_k\in [t_0/4,3t_0/2]$, so we can apply \eqref{eqs}. Then, the above inequality implies
           \begin{align}
           T_{t}\chi_{x}(\tilde z_k)- T_{t}\chi_{x}(y)
           \geq 
           c_{\tau_k}(y,\tilde z_k)-\frac 1k-Ms.  \label{eqqq}
           \end{align}
           Moreover, since $T_{t}\chi_{x}$ is locally semiconcave on $I^+(x)$, the mapping $(0,\tau_k)\ni \tau \mapsto T_{t}\chi_{x}(\gamma_k(\tau))$ is locally semiconcave as well as the composition of a locally semiconcave with a smooth function (\cite{Fathi/Figalli}, Lemma A.9.). Thus, it is almost everywhere differentiable with
           \begin{align*}
               \frac{d}{d\tau} (T_{t}\chi_{x}\circ \gamma_k)(\tau)= p_\tau(\dot \gamma_k(\tau)) 
           \end{align*}
           for some (or any) $p_\tau\in  \partial^+ T_{t}\chi_{x}(\gamma_k(\tau))$. In particular, since ${\gamma_k}_{|[0,\tau_k]}$ maps to $\overbar B_r(y_0)$, we can apply \eqref{eqp} and obtain
           \begin{align*}
               T_{t}\chi_{x}(\tilde z_k)- T_{t}\chi_{x}(y)
               &=
               \int_0^{\tau_k} \frac{d}{d\tau} (T_{t}\chi_{x}\circ \gamma_k)(\tau)\, d\tau= \int_0^{\tau_k} p_\tau(\dot \gamma_k(\tau)) \, d\tau
               \\[10pt]
               &\leq 
               -C \int_0^{\tau_k} |\dot \gamma_k(\tau)|_{ h}\, d\tau
               \\[10pt]
               &\leq -Cd_{h}(y,\tilde z_k)
               =-CC_0\sqrt{s}.
           \end{align*}
           We now estimate the right-hand side of \eqref{eqq}. Thanks to \eqref{eqw}, we get
           \begin{align*}
              c_{\tau_k}(y,\tilde z_k)-\frac 1k-Ms
              \geq 
              -\sqrt{s}M-\frac 1k -Ms \geq -2M\sqrt{s}-\frac 1k.
           \end{align*}
           Putting everything together, \eqref{eqqq} implies
           \begin{align}
               -CC_0\sqrt{s} 
               \geq 
               -2M\sqrt{s} -\frac 1k, \label{eqr}
           \end{align}
           but this contradicts the definition of $C_0$ if $k$ is large.
\end{proof}

\begin{lemma}\label{l2}
     Let $x_0\in M$, $y_0\in I^+(x_0)$ and $t_0>0$. Then there exist two open neighbourhoods $U$ and $V$ of $x_0$ and $y_0$, respectively, with $U\times V\subseteq I^+$, and some number $s_0>0$ such that, for all $s\in (0,s_0]$, $t\in [t_0/2,3t_0/2]$ and $(x,y)\in U\times V$, we have
    \begin{align}
        \hat T_sT_{t}\chi_x(y)>\ T_t\chi_x(y). \label{eqac}
    \end{align}
\end{lemma}
\begin{proof}
    Let $U, V$, $s_0$ and $C_0$ be as in Lemma \ref{l1}, and let $s,t,x,y$ be as in the statement. Let $\gamma:[0,1]\to M$ be a maximizing geodesic connecting $x$ to $y$. 
    
    For small $\ep>0$, we can extend $\gamma$ to a geodesic defined on the interval $[0,1+\ep]$. We will prove that, for $\ep>0$ sufficiently small,
    \begin{align}
        T_t\chi_x(\gamma(1+\ep))-c_s(y,\gamma(1+\ep))>T_t\chi_x(y). \label{eqab}
    \end{align}
    As a composition of a locally semiconcave (hence locally Lipschitz) with a smooth function, $T_t\chi_x\circ \gamma$ is locally Lipschitz, hence there is $C>0$ such that    \begin{align}
        |T_t\chi_x(\gamma(1+\ep))-T_t\chi_x(y)|
        \leq C\ep \label{eqg}
    \end{align}
    for small $\ep>0$.
    
    Moreover, since $\gamma$ is a maximizing geodesic, its "speed" $|\dot \gamma(\tau)|_g$ is constant and equal to $d(x,y)$.
    Thus, $d(y,\gamma(1+\ep))\geq \ell_g({\gamma_{|[1,1+\ep]}})=\ep d(x,y)$. Hence
    \begin{align}
       -c_s(y,\gamma(1+\ep))
        \geq \sqrt{s}(\ep d(x,y))^\frac 12. \label{eqh}
    \end{align}
    Since $d(x,y)>0$, we can combine inequalities \eqref{eqg} and \eqref{eqh} to obtain \eqref{eqab} for sufficiently small $\ep$.
    \end{proof}

\begin{corollary} \label{c1}
    Let $x_0\in M$, $y_0\in I^+(x_0)$ and $t_0>0$. Then there exist two open neighbourhoods $U$ and $V$ of $x_0$ and $y_0$, respectively, with $U\times V\subseteq I^+$, and some number $s_0>0$ such that, for all $s\in (0,s_0]$, $t\in [t_0/2,3t_0/2]$ and $(x,y)\in U\times V$, the supremum in the definition of $\hat T_sT_{t}\chi_x(y)$ is attained at some $z\in M$, and necessarily $z\in I^+(y)$. In particular, the mapping $(0,s_0]\times [t_0/2,3t_0/2]\times U\times V\ni (s,t,x,y)\mapsto \hat T_sT_{t}\chi_x(y)$ is lower semicontinuous.
\end{corollary}
\begin{proof}
     Let $U, V$, $s_0$ and $C_0$ be such that both Lemma \ref{l1} and Lemma \ref{l2} apply. 
     
     Fix $s\in (0,s_0]$, $t\in [t_0/2,3t_0/2]$ and $(x,y)\in U\times V$, and let $z_k$ be a maximizing sequence in the definition of $\hat T_sT_{t}\chi_x(y)$. Then, as $U\times V\subseteq I^+$, we have $z_k\in J^+(y)\subseteq I^+(x)$. By Lemma \ref{l1} and the completeness of the metric $h$, it follows that, up to a subsequence, $z_k\to z\in J^+(y)\subseteq I^+(x)$. The continuity of $\C$ implies
    \begin{align*}
        \hat T_sT_{t}\chi_x(y)
        =\lim_{k\to \infty} (T_{t}\chi_x(z_k)-c_s(y,z_k))
        =
        T_{t}\chi_x(z)-c_s(y,z).
    \end{align*}
    This shows that the supremum is attained at some $z$. Assuming momentarily $z\in I^+(y)$, let us show lower semicontinuity (at the point $(s,t,x,y)$). 
    
    Let $(s_k,t_k,x_k,y_k) \in (0,s_0]\times [t_0/2,3t_0/2]\times U\times V$ with $(s_k,t_k,x_k,y_k)\to (s,t,x,y)$. Since $y\in I^-(z)$ and $I^-(z)$ is open, also $y_k\in I^-(z)$ for large $k$, i.e.\ $z\in I^+(y_k)$. Hence, by continuity of $\C$,
    \begin{align*}
        \liminf_{k\to \infty} \hat T_{s_k}T_{t_k}\chi_{x_k}(y_k)
        \geq \liminf_{k\to \infty}\,  [T_{t_k}\chi_{x_k}(z)-c_{s_k}(y_k,z)] &= T_{t}\chi_x(z)-c_{s}(y,z)
        \\
        &=\hat T_sT_{t}\chi_x(y).
    \end{align*}
    This proves the lower semicontinuity and we are left to show $z\in I^+(y)$ for any possible optimal $z$.
    
    Let $z\in M$ be optimal. Then $z\in J^+(y)$. Suppose, for contradiction, that $d(y,z)=0$, and let $\gamma_0:[0,s]\to M$ be a causal null geodesic connecting $y$ with $z$. By Lemma \ref{l2}, $z\neq y$, so $\gamma_0$ is non-constant. Hence, we can choose $\xi_0\in T_{z}M$ with $g_z(\dot \gamma_0(s),\xi_0)<0$. Let $\xi:[0,s]\to TM$ be the parallel transport (w.r.t.\ the Levi-Civita connection of $g$) of $\xi_0$ along $\gamma_0$, and consider a smooth variation $\gamma:(-\ep,\ep)\times [0,s]\to M$ of $\gamma_0$ with variational vector field $\tilde \xi(\tau)=\tau \xi(\tau)$, fixing $\gamma(r,0)=\gamma_0(0)=y$ for all $r\in (-\ep,\ep)$. Our goal is to show that, for $r>0$ small enough,
    \begin{align}
        T_{t}\chi_x(\gamma(r,s))-c_s(y,\gamma(r,s))>T_{t}\chi_x(z)-c_s(y,z), \label{eqx}
    \end{align}
    contradicting the optimality of $z$, proving the claim. 

    By the compatibility of the connection with the metric, we compute
    \begin{align*}
        \frac{d}{dr}\Big\vert_{r=0} g(\partial_\tau \gamma(r,\tau),\partial_\tau \gamma(r,\tau))
        &= 2 g\bigg(\frac{D}{dr}\Big\vert_{r=0}\frac{\partial \gamma}{\partial \tau}(r,\tau),\frac{\partial \gamma}{\partial \tau}(0,\tau)\bigg)
        \\[10pt]
        &=
        2 g\bigg(\frac{D}{d\tau}\frac{\partial \gamma}{\partial r}(0,\tau),\dot \gamma_0(\tau)\bigg)
        \\[10pt]
        &=
        2 g\bigg(\frac{D}{d\tau}(\tau\xi(\tau)),\dot \gamma_0(\tau)\bigg)
        \\[10pt]
        &= 2 g(\xi(\tau),\dot \gamma_0(\tau))
        \\[10pt]
        &= 2 g(\xi(s),\dot \gamma_0(s))=:-2a<0,
    \end{align*}
    where in the step to the last line we used the fact that $\xi$ and $\dot \gamma_0$ are parallel along $\gamma_0$.
    Thus, Taylor expansion and the fact that $\gamma_0$ is a null geodesic yield
    \begin{align}
        g(\partial_\tau \gamma(r,\tau),\partial_\tau \gamma(r,\tau))\leq -|\dot \gamma_0(\tau)|^2_g -2ar +O(r^2)\leq -ar \label{hhhh}
    \end{align}
    for small values $r>0$ and all $\tau\in [0,s]$. In particular, $\gamma(r,\cdot)$ is either (future-directed) timelike or past-directed timelike for small $r$.
    Since $y\neq z$ and $z\in J^+(y)$ we must have $z\notin J^-(y)$. Thus $z\approx \gamma(r,s)\notin J^-(z)$ for small $r$, so that $\gamma(r,\cdot)$ is in fact (future-directed) timelike for small $r$.
This, together with
\begin{align*}
    \ell_g(\gamma(r,\cdot))\geq s \sqrt{ar},
\end{align*}
as follows from \eqref{hhhh}, gives
\begin{align*}
    d(y,\gamma(r,s))\geq  s\sqrt{ar}
\end{align*}
and therefore
\begin{align}
    -c_s(y,\gamma(r,s))+c_s(y,z)= -c_s(y,\gamma(r,s))=s^\frac 12 d(x,\gamma(r,s))^\frac 12 
    \geq s^\frac 12 (ar)^\frac 14 \label{hwdauidjaias}
\end{align}
for small $r$.
On the other hand, since $\partial_r \gamma(0,s)=s\xi_0$, it holds $d_h(\gamma(r,s),z)\leq 2s|\xi_0|_hr$, if $r$ is small. Thus, denoting by $L$ a local Lipschitz constant of $T_t\chi_x$ near $z\in I^+(x)$, it follows for small $r$ that
\begin{align*}
    |T_{t}\chi_x(\gamma(r,s))-T_{t}\chi_x(z)| \leq 2Ls|\xi_0|_hr.
\end{align*}
Combining this inequality with \eqref{hwdauidjaias}, we see that that \eqref{eqx} holds for small $r$.
\end{proof}

\begin{corollary}\label{c3}
     Let $x_0\in M$, $y_0\in I^+(x_0)$ and $t_0>0$. Then there exist two open neighbourhoods $U$ and $V$ of $x_0$ and $y_0$, respectively, with $U\times V\subseteq I^+$, some number $s_0>0$, and a non-decreasing sequence $(\delta_s)_{0<s\leq s_0}$ of positive numbers such that, for all $s\in (0,s_0]$, $t\in [t_0/2,3t_0/2]$ and $(x,y)\in U\times V$, we have
    \begin{align}
        \hat T_sT_{t}\chi_x(y)>\sup\{T_{t}\chi_x(z)-c_s(y,z)\mid d(y,z)\leq \delta_s\}. \label{eqad}
    \end{align}
\end{corollary}
\begin{proof}
     Let $\tilde U,\tilde V$, $s_0$ and $C_0$ be given by Lemma \ref{l1} and Corollary \ref{c1}. Choose two open neighbourhoods $U\Subset\tilde U$ and $V\Subset \tilde V$ of $x$ and $y$, respectively. It suffices\footnote{
     Indeed, we can set $\delta_s:=\delta_{s_1}'$ for $s\in [s_1,s_0]$, $\delta_s:=\min\{\delta_{s_1/2}',\delta_{s_1}'\}$ for $s\in [s_1/2,s_1]$ and so on.
     } 
     to show that, if $s_1\in (0,s_0]$, there exists $\delta'_{s_1}$ such that \eqref{eqad} holds for $s\in [s_1,s_0]$, $t\in [t_0/2,3t_0/2]$ and $(x,y)\in U\times V$. 
     
     Suppose the contrary. Then we can find sequences $s_k\in [s_1,s_0]$, $t_k\in [t_0/2,3t_0/2]$ and $(x_k,y_k)\in U\times V$ such that
     \begin{align*}
        \hat T_{s_k}T_{t_k}\chi_{x_k}(y_k)=\sup\{T_{t_k}\chi_{x_k}(z_k)-c_{s_k}(y_k,z_k)\mid d(y_k,z_k)\leq 1/k\}. 
    \end{align*}
     This, together with Lemma \ref{l1}, implies that for each $k\in \N$ we can find  $z_k\in J^+(y_k)$ with $d(y_k,z_k)\leq 1/k$, $d_h(y_k,z_k)\leq C_0\sqrt{s_0}$ and
    \begin{align*}
        T_{t_k}\chi_{x_k}(z_k)-c_{s_k}(y_k,z_k)
        \geq
        \hat T_{s_k}T_{t_k}\chi_{x_k}(y_k)-\frac 1k.
    \end{align*}
    Up to subsequences, $s_k\to s\in [s_1,s_0]$, $t_k\to [t_0/2,3t_0/2]$, $(x_k,y_k)\to (x,y)\in \overbar U\times \overbar V\subseteq \tilde U\times \tilde V$ and $z_k\to z\in J^+(y)$ with $d(y,z)=0$. Moreover, by the lower semicontinuity of $(s,t,x,y)\mapsto \hat T_sT_{t}\chi_x(y)$ on $(0,s_0]\times [t_0/2,3t_0/2]\times U\times V$ and the continuity of $\C$,
    \begin{align*}
        \hat T_sT_{t}\chi_x(y) 
        \leq \liminf_{k\to \infty} 
        \hat T_{s_k}T_{t_k}\chi_{x_k}(y_k)
        &\leq \liminf_{k\to \infty} 
        \big[T_{t_k}\chi_{x_k}(z_k)-c_{s_k}(y_k,z_k)\big]
        \\
        &= T_t\chi_x(z)-c_s(y,z).
    \end{align*}
   By definition of $\hat T_sT_t\chi_x(y)$, this inequality must actually hold as an equality. However, Corollary \ref{c1} shows that this is impossible since $d(y,z)=0$.
\end{proof}

\begin{proof}[Proof of Theorem \ref{thmb}]
    Let $U,V,s_0,C_0$ and $(\delta_s)_{0<s<s_0}$ be given by Lemma \ref{l1} and the above corollary. \eqref{eqae} follows from the above corollary and Lemma \ref{l1}. The existence of a maximizer $z\in M$ in the definition of $\hat T_sT_t\chi_x(y)$ is trivial when $s=0$ (namely $z=y$), and follows from \eqref{eqae} for $s>0$. Moreover, \eqref{eqae} guarantees that $z\in I^+(y)$ whenever $s>0$. 
    
    Next, fix $s_1\in (0,s_0]$ and $y\in V$. By \eqref{eqae} there exists an open neighbourhood $V'\subseteq V$ of $y$ and a compact set $K$ with $V'\times K\subseteq I^+$ and such that, for any $s\in [s_1,s_0]$, $t\in [t_0/2,3t_0/2]$, $x\in U$ and $y'\in V'$, we have
    \begin{align*}
        \hat T_sT_{t}\chi_x(y')=\sup\{T_{t}\chi_x(z)-c_s(y',z)\mid z\in K\}.
    \end{align*}
    Since also $U\times K\subseteq I^+$, and $\C$ is locally semiconcave on $(0,\infty)\times I^+$, the family of functions 
    \[
    ([s_1,s_0]\times [t_0/2,3t_0/2]\times U\times V'\ni (s,t,x,y')\mapsto T_t\chi_x(z)-c_s(y',z))_{z\in K},
    \]
    is locally equi-continuous. Hence, their pointwise finite supremum, namely $\hat T_sT_t\chi_x(y)$, is continuous.
   Since $y\in V$ and $s_1$ were arbitrary, it follows that the mapping
   \begin{align*}
       (0,s_0]\times [t_0/2,3t_0/2]\times U\times V\ni (s,t,x,y)\mapsto \hat T_sT_t\chi_x(y),
   \end{align*}
   is continuous. To establish continuity at points of the form $(0,t,x,y)$, note first that $\hat T_0T_t\chi_x(y)=T_t\chi_x(y)$. Additionally, as follows from the definition,
   \begin{align*}
       T_t\chi_x(y)\leq \hat T_sT_t\chi_x(y) \leq T_{t-s}\chi_x(y)
   \end{align*}
   as soon as $s\leq t$ (cf.\ Lemmas 3.3 and 3.4 in \cite{Metsch2}). In particular, for any sequence $(s_k,t_k,x_k,y_k)\in [0,s_0]\times [t_0/2,3t_0/2]\times U\times V$ converging to $(0,t,x,y)$, the continuity of $\C$ implies
   \begin{align*}
       T_t\chi_x(y) \leq \lim_{k\to \infty} \hat T_{s_k}T_{t_k}\chi_{x_k}(y_k)\leq T_t\chi_x(y),
   \end{align*}
   which yields the desired continuity.
    
    Finally, fix again $s_1\in (0,s_0]$ and $y\in V$. Define $V'$ and $K$ as in the beginning of the proof. By the same reasoning as above ($\C$ is locally semiconcave on $I^+$), the family $(-c_s(\cdot,z))_{s\in [s_1,s_0],z\in K}$ is uniformly locally semiconvex on $V'$ (\cite{Fathi/Figalli}, Proposition A.17 and Proposition A.4). Thus, the family of functions $(T_{t}\chi_x(z)-c_s(\cdot,z))_{s\in [s_1,s_0],t\in [t_0/2,3t_0/2],x\in U,z\in K}$ is uniformly locally semiconvex on $V'$. Therefore, also the family $(\sup\{T_{t}\chi_x(z)-c_s(\cdot,z)\mid z\in K\})_{s\in [s_1,s_0],t\in [t_0/2,3t_0/2],x\in U}$ is uniformly locally semiconvex on $V'$, provided the suprema are everywhere finite (\cite{Fathi/Figalli}, Theorem A.11). The suprema are, however, precisely $\hat T_sT_t\chi_x$. Since $y\in V$ was arbitrary, this concludes the rest of the proof.
\end{proof}

\section{Local semiconcavity of $\hat T_sT_{t}\chi_x$}

Our main theorem in this section is the following, which is somehow the converse to Theorem \ref{thmb}.
\begin{theorem}\label{thmc}
Let $x_0\in M$, $y_0\in I^+(x_0)$ and $t_0>0$. Then there exist two open neighbourhoods $U$, $V$ of $x_0$ and $y_0$, respectively, with $U\times V\subseteq I^+$, an open interval $I$ containing $t_0$, and some number $s_0>0$ such that the family of functions $\{\hat T_sT_{t}\chi_x\mid s\in [0,s_0],t\in I,x\in U\}$ is uniformly locally semiconcave on $V$. 
\end{theorem}

We prove the theorem in several steps. We start with the following general lemma.

\begin{lemma}
    Let $N$ be a pseudo-Riemannian manifold and let $x_0\in N$. Then there exists a chart $(\phi,U)$ centered at $x_0$ (i.e.\ $\phi(x_0)=0)$ such that the following holds: Whenever $r>0$ is such that $B_r(0)\subseteq \phi(U)\subseteq \R^n$ and $\gamma:[a,b]\to U$ is a geodesic with $\phi(\gamma(a))\in B_r(0)$ and $\phi(\gamma(t_0))\notin B_r(0)$ for some $t_0\in [a,b]$, then the map 
    \begin{align*}
        [t_0,b]\ni t\mapsto |\phi(\gamma(t))| \text{ is non-decreasing.}
    \end{align*}
		Here, $|\cdot|$ denotes the Euclidean norm in $\R^n$.
\end{lemma}
\begin{proof}
    Since the statement is local, we may assume that $N$ is an open subset of $\R^n$ equipped with some pseudo-Riemannian metric and that $x_0=0$. 
		
    Let $\ep>0$ and choose an open precompact neighbourhood $V$ of $0$ such that, for any $x\in \overbar V$ and any unit vector $v\in S^{n-1}$, the geodesic $\gamma_{x,v}(t)=\exp_x(tv)$ exists and is injective on $[-\ep,\ep]$. Define
    \begin{align*}
        &C:=\sup\{|\ddot \gamma_{x,v}(t)|\mid x\in V,\ v\in S^{n-1},\ t\in [-\ep,\ep]\}<\infty, 
        \\[10pt]
        &\delta_1:=\inf\{|\dot \gamma_{x,v}(t)|\mid x\in V,\ v\in S^{n-1},\ t\in [-\ep,\ep]\}>0 \text{ and } 
        \\[10pt]
        &\delta_2:=\inf\{|\gamma_{x,v}(\ep)-x|\mid x\in V,\ v\in S^{n-1}\}>0.
    \end{align*}
    Now set
    \[
    R:= \min\bigg\{\frac{\delta_2}{2},\frac{\delta_1^2}{C}\bigg\}
    \]
    and $U:=B_R(0)$. 
		
    Assume for contradiction that we find $r<R$, a geodesic $\gamma:[a,b]\to U$ with $\gamma(a)\in B_r(0)$, $t_0\in [a,b]$ with $\gamma(t_0)\notin B_r(0)$, and $b\geq t_2\geq t_1\geq t_0$ with $|\gamma(t_2)|< |\gamma(t_1)|$. Reparametrizing, we may assume that $a=0$ and $|\dot \gamma(a)|=1$. By continuity, there must exist $t_3\in (a,t_2)$ such that $|\gamma(t)|$ attains a local maximum at $t=t_3$. At this point, the second derivative is non-positive:
    \begin{align*}
        |\dot \gamma(t_3)|^2+\langle \gamma(t_3),\ddot \gamma(t_3)\rangle =\frac 12\frac{d^2|\gamma|^2}{dt^2}(t_3)\leq 0.
    \end{align*}
		Here, $\langle\cdot,\cdot\rangle$ denotes the Euclidean scalar product. 
    Observe that, if $b\geq \ep$, then $|\gamma(\ep)|\geq |\gamma(\ep)-\gamma(0)|-|\gamma(0)|\geq \delta_2-R\geq  R$, so that $\gamma(\ep)\notin U$, which contradicts the definition of $\gamma$. Thus, $b\leq \ep$. Therefore, we also have $t_3\in [0,\ep]$ and
    \begin{align*}
        |\dot \gamma(t_3)|^2+\langle \gamma(t_3),\ddot \gamma(t_3)\rangle
        \geq 
        \delta_1^2 -|\gamma(t_3)||\ddot \gamma(t_3)|
        > \delta_1^2-RC
        \geq 0.
    \end{align*}
    This is a contradiction and, hence, concludes the proof.
\end{proof}

In the following proposition, we use the notion of $C^2$-boundedness. For definitions, see Subsection 6.1 in \cite{Metsch2}.

\begin{proposition}\label{pro1}
 Let $N$ be a smooth manifold and $f:N\times M\to  \overbar \R$ be a function that is locally semiconcave on an open neighbourhood of $(q_0,y_0)$. Suppose further that
 \begin{align*}
    \partial^+ f_{q_0}(y_0)\Subset \op{int}(\C_{y_0}),
 \end{align*}
 where $f_{q_0}(y):=f(q_0,y)$.
    
    Then there exists an open neighbourhood $U_1$ of $q_0$ and three open neighbourhoods $V_1\Subset V_2\Subset V_3$ of $y_0$ as well as a chart $\phi:V_3\to W_3\subseteq \R^n$ and a family of smooth functions $f_{i,q}:V_3\to \R$, $i\in I_{q}$, $q\in U_1$ such that 
    \begin{enumerate}[(i)]
        \item 
        The family $(f_{i,q}\circ \phi^{-1})_{i\in I_q,q\in U_1}$ is bounded in $C^2(W_3)$.
        
        \item 
        The set
        \[
        \{(y,d_yf_{i,q})\mid y\in  V_1,\ q\in U_1,i\in I_{q}\}
        \]
        is relatively compact in $\op{int}(\C^{*})$.
        
        \item 
         $\forall q\in U_1$: $f(q,\cdot)= \inf_{i\in I_{q}} f_{i,q} \text{ on } V_1$.

        \item 
         $\forall q\in U_1,\ y\in V_1,\ p\in \partial^+f_{q}(y)$: $\exists i\in I_q$: $f_{i,q}(y)=f(q,y)$ and $d_yf_{i,q}=p$.
         
         \item 
        $\forall q\in U_1,\ i\in I_q$: $f(q,y)> f_{i,q}(y)$ for all $y\notin V_2$.
        
         \item Any causal curve with start- and endpoint in $V_2$ lies entirely within $V_3$.
        \item 
         There is $C>0$ such that, for all $q\in U_1$, $i\in I_q$, and any causal geodesic $\gamma:[a,b]\to V_3$ with $\gamma(a)\in V_1$, it holds
         \begin{align*}
             &\frac{d}{dt}(f_{i,q}\circ \gamma)(t)\leq -C|\dot \gamma(t)|_h \text{ whenever } \gamma(t)\in V_1 \text{ and }
             \\[10pt]
             &\frac{d}{dt}(f_{i,q}\circ \gamma)(t)\leq 0 \text{ for all } t\in [a,b].
         \end{align*}
    \end{enumerate}
\end{proposition}

\begin{proof}
    Let $W\subseteq N\times M$ be an open neighbourhood of $(q_0,y_0)$ where $f$ is locally semiconcave. Let $(\tilde \phi,\tilde V_3)$ be a chart around $y_0$ as in the preceding lemma. Since locally semiconcave functions are locally Lipschitz (\cite{Fathi/Figalli}, Lemma A.5), there exists a relatively compact open neighbourhood $V_3\Subset\tilde V_3$ of $y_0$, and an open neighbourhood $U_2\subseteq U$ of $q_0$ with $U_2\times V_3\subseteq W$ such that:
     \begin{align}
        \forall q\in U_2:\ f_q\circ \phi^{-1} \text{ is $K$-concave, $K$-Lipschitz and bounded by $K$}. \label{eqn}
     \end{align}
     Here, $f_q:=f(q,\cdot)$. Set $\phi:=\tilde \phi_{|V_3}$ and $W_3:=\phi(V_3)\subseteq \R^n$.
     
     For each $(q,y)\in U_2\times V_3$ and $p\in \partial^+f_q(y)$, we define a smooth function $\tilde f_{y,p,q}:V_3\to \R$ by
     \begin{align*}
         &\tilde f_{y,p,q}\circ \phi^{-1}(x):= f_q(y)+p\circ d_{\phi(y)}\phi^{-1}(x-\phi(y))+K|x-\phi(y)|^2.
     \end{align*}
     Note that, by \eqref{eqn}, the functions $f_q$, $q\in U_2$, are bounded by $K$ on $V_3$, and the linear maps $p\circ d_{\phi(y)}\phi^{-1}$, $(q,y)\in U_2\times V_3$, $p\in \partial^+ f_q(y)$, have operator norm bounded by $K$. Hence, being $W_3$ precompact, there exists a constant $\tilde K>0$ such that all functions $\tilde f_{y,p,q}\circ \phi^{-1}$, $(q,y)\in U_2\times V_3$, $p\in \partial^+f_q(y)$, are equi-Lipschitz and bounded by $\tilde K$.
     
     Moreover, the continuity of the super-differential of a locally semiconcave function as multivalued map implies that $d_{y'}\tilde f_{y,p,q}\to\partial^+f_{q_0}(y_0)$ in the topology of the tangent bundle whenever $y',y\to y_0$, $q\to q_0$ and $p\in \partial^+f_q(y)$. By assumption, $\partial^+ f_{q_0}(y_0)\Subset \op{int}(\C_{y_0}^*)$, so there exist $V_2'\subseteq V_3$ and $U_1\subseteq U_2$, two open neighbourhoods of $y_0$ and $q_0$, respectively, with
     \begin{align}
        \{(y',d_{y'}\tilde f_{y,p,q})\mid y',y\in V_2',\ q\in U_1,p\in \partial^+f_q(y)\} \Subset \op{int}(\C^{*}). \label{equ} 
    \end{align}
    In particular, there exists a constant $C>0$ such that, for any $y,y'\in  V_2'$, $q\in U_1$ and $p\in \partial^+f_q(y)$, we have
        \begin{align}
          d_{y'}\tilde f_{y,p,q}(v)\leq -C|v|_h \text{ for all } v\in \C_{y'}. \label{equu}
        \end{align}
      
      Now let $V_2\Subset  V_2'$ be an open neighbourhood of $y_0$ such that any causal curve with start- and endpoint in $V_2$ lies entirely in $V_3$. Such a neighbourhood exists since $M$ is globally hyperbolic, hence strongly causal. Let $V_1\Subset V_2$ be any coordinate ball as in the foregoing lemma (w.r.t.\ the chart $(\phi,V_3)$).
			
	Let us define, for $(q,y)\in U_1\times V_1$ and $p\in \partial^+f_q(y)$, the family of functions
     \begin{align}
        f_{y,p,q}\circ \phi^{-1}:W_3\to \R,\ 
        x\mapsto \rho(x) (\tilde f_{y,p,q}\circ \phi^{-1})(x)
        - 2 \tilde K(1-\rho(x)).
        \label{rrrr}
     \end{align}
    Here, $\rho:\R^n\to [0,1]$ is a smooth function satisfying $\rho\equiv 1$ on $\phi(V_1)$, $\rho\equiv 0$ on $W_3\backslash \phi(V_2)$, and $\rho(x)\leq \rho(x')$ for any $|x|\geq |x'|$.
    
    We finally define, for each $q\in U_1$, the index set
		\begin{align*}
		I_q:=\{(y,p)\mid y\in V_1,\ p\in \partial^+f_q(y)\},
		\end{align*}
		and we consider the family of smooth functions
		\begin{align*}
		(f_{i,q}:V_3\to \R)_{q\in U_1, i\in I_q}.
		\end{align*}
    \begin{enumerate}[(i)]
        \item Thanks to the fact that $\rho\equiv 0$ outside $\phi(V_2)$ and that the family $\tilde f_{i,q}\circ \phi^{-1}$, $i\in I_q$, $q\in U_1$, is bounded by $\tilde K$ and equi-Lipschitz, it suffices to prove that the family of maps $\tilde f_{i,q}\circ \phi^{-1}$, $i\in I_q$, $q\in U_1$, is bounded in $C^2(W_3)$. This, however, follows from the definition of $\tilde f_{i,q}$.
        
        \item Since $\rho\equiv 1$ on $\phi(V_1)$, this follows immediately from \eqref{equ}.
        
        \item By the $K$-concavity of $f_q\circ \phi^{-1}$, $q\in U_1$, it is clear that $\tilde f_{y,p,q}\geq f_{q}$ on $V_3$, implying that $f_{i,q}\geq f_{q}$ on $V_1$ for all $i\in I_q$. The fact that $f_{q}=\inf_{i\in I_q} f_{i,q}$ follows from (iv) (note that, by the local semiconcavity of $f$, $\partial^+ f_{q}(y)$ is never empty).
        
        \item Given $(q,y)\in U_1\times V_1$ and $p\in \partial^+ f_{q}(y)$, consider $i=(y,p)\in I_q$. Then $f_{i,q}(y)=f_q(y)$ and $d_yf_{i,q}=d_y\tilde f_{y,p,q}=p$.

        \item We have $f_{i,q}\equiv -2\tilde K<-K$ on $V_3\backslash V_2$, while $|f|$ is bounded by $K$ on $U_1\times V_3$.
        
        \item This was the definition of $V_2$.

        \item Given $q\in U_1$ and $i=(y,p)\in I_q$, let $\gamma:[a,b]\to V_3$ be a causal geodesic. Suppose first that $\gamma(t)\in V_1$. Since $f_{i,q}=\tilde f_{y,p,q}$ on $V_1$, \eqref{equu} gives
        \begin{align*}
            \frac{d}{dt}(f_{i,q}\circ \gamma)(t)
            =
            d_{\gamma(t)}f_{i,q}(\dot \gamma(t))
            \leq 
            -C|\dot \gamma(t)|_h.
        \end{align*}
         On the other hand, if $\gamma(t)\notin \overbar V_2$, then we have $d_{\gamma(t)}f_{i,q}=0$ as $f_{i,q}\equiv -2\tilde K$ on $V_3\backslash V_2$. Finally, if $\gamma(t)\in \overbar V_2\backslash V_1$,  we have
        \begin{align*}
           \frac{d}{dt}( \tilde f_{i,s}\circ \gamma)(t)
           = (\tilde f_{y,p,q}(\gamma(t))+2\tilde K)\cdot \frac{d}{dt} (\rho\circ \gamma)(t)
           +\rho(\gamma(t))\cdot  \frac{d}{dt}( \tilde f_{y,p,q}\circ \gamma)(t).
        \end{align*}
        Now, $(\tilde f_{y,p,q}(\gamma(t))+2\tilde K)\geq 0$ by definition of $\tilde K$, and the foregoing lemma states that $|\gamma|$ is increasing on $[t,b]$, hence $\rho\circ \gamma$ in decreasing on $[t,b]$. In total, the first term on the right-hand side in the above inequality is non-positive. For the second term, note that $\rho(\gamma(t))\geq 0$, and since $V_2\subseteq V_2'$, \eqref{equu} implies
        \begin{align*}
            \frac{d}{dt}(\tilde f_{y,p,q}\circ \gamma)(t)
            =
            d_{\gamma(t)}\tilde f_{y,p,q}(\dot \gamma(t))
            \leq 
            -C|\dot \gamma(t)|_h\leq 0.
        \end{align*}
    \end{enumerate}
    This proves all required properties.
\end{proof}

\begin{remark}\rm\label{rema}
    In this remark, we fix our notation which will be used throughout \textbf{the rest of the chapter}. 
    
    Let $t_0>0$ and $(x_0,y_0)\in I^+$ be fixed, let $N:=(0,\infty)\times M$ and let $f:=\C:(0,\infty)\times M^2\to \R\cup\{+\infty\}$ (here, $q=(t,x)$). Then $f$ is locally semiconcave on $(0,\infty)\times I^+$. Using that the super-differential of a locally semiconcave function at any point is compact (in the cotangent space), Lemma \ref{l4} implies that
    \begin{align*}
        \partial^+f_{t_0,x_0}(y_0)\Subset \op{int}(\C_{y_0}).
    \end{align*}
    Thus, we can apply the preceding lemma. Let $I$ be an open interval containing $t_0$ and $U$ be an open neighbourhood of $x_0$ such that $I\times U\subseteq U_1$. Let $V_1,V_2,V_3$, $\phi:V_3\to W_3$ and the constant $C$ be as in the lemma, and let $f_{i,t,x}$, $(t,x)\in I\times U$, $i\in I_{t,x}$, be the associated family of functions. Moreover, fix an open neighbourhood $V\Subset V_1$ of $y_0$ with $\phi(V)$ being a convex set in $\R^n$. By possibly shrinking $U$ and $V$ if necessary, we may assume that Theorem \ref{thmb} is applicable with $U\times V\subseteq I^+$ and constants $s_0',C_0>0$. Finally, let \begin{align*}
        s_1\in (0,s_0'] \text{ such that } C_0\sqrt{s_1}\leq d_{h}(V,\partial V_1).
    \end{align*}
\end{remark}

\begin{definition}\rm
    For $s\geq 0$ and a map $g:V_3\to \R$, we define the function
    \begin{align*}
        \hat T_s^{loc}g:V_3\to \R\cup\{+\infty\},\
        \hat T_s^{loc}g(y):=\sup\{g(z)-c_s(y,z)\mid z \in V_3\}.
    \end{align*}
\end{definition}

\begin{lemma}\label{li}
    Let $s\in [0,s_1]$, $t\in I$ and $x\in U$. Then 
    \begin{align*}
        \inf_{i\in I_{t,x}} \hat T^{loc}_sf_{i,t,x}\geq \hat T_sT_{t}\chi_x \text{ on } V.
    \end{align*}
\end{lemma}
\begin{proof}
    Let $y\in V$. By the choice of $s_1$, Theorem \ref{thmb} states that there exists $z\in V_1$ with
    \begin{align*}
        \hat T_sT_{t}\chi_x(y)=T_{t}\chi_x(z)-c_s(y,z).
    \end{align*}
    Given $i\in I_{t,x}$, since $f_{i,t,x}\geq T_{t}\chi_x$ on $V_1$ by Proposition \ref{pro1}(iii), it follows that
    \begin{align*}
        \hat T^{loc}_sf_{i,t,x}(y)
        \geq
        f_{i,t,x}(z)-c_s(y,z)
        \geq
        T_{t}\chi_x(z)-c_s(y,z)
        =
        \hat T_sT_{t}\chi_x(y).
    \end{align*}
    This concludes the proof.
\end{proof}

\begin{corollary}\label{cor3}
    Let $s\in [0,s_1]$, $t\in I$, $x\in U$ and $i\in I_{t,x}$. Then
    \begin{align*}
        \hat T^{loc}_sf_{i,t,x}(y)>\sup\{f_{i,t,x}(z)-c_s(y,z)\mid z\in V_3\backslash V_2\}
    \end{align*}
    for all $y\in V$.
\end{corollary}
\begin{proof}
    If $s=0$, this is trivial. Thus, suppose $s> 0$. Proposition \ref{pro1}(v), $C_0\sqrt{s_1}\leq d_h(V,\partial V_1)$, Theorem \ref{thmb} and the above lemma give for $y\in V$
    \begin{align*}
       \sup\{f_{i,t,x}(z)-c_s(y,z)\mid z\notin V_2\} 
       &\leq
       \sup\{T_{t}\chi_x(z)-c_s(y,z)\mid z\notin V_2\}
       \\[10pt]
       &\leq
        \sup\{T_{t}\chi_x(z)-c_s(y,z)\mid d_{h}(y,z)\geq C_0\sqrt{s}\}
        \\[10pt]
        &<
        \hat T_sT_{t}\chi_x(y)
        \\[10pt]
        &\leq \hat T^{loc}_sf_{i,t,x}(y).
    \end{align*}
    This proves the lemma.
\end{proof}

\begin{lemma} \label{lj}
Let $s\in (0,s_1]$, $t\in I$, $x\in U$ and $i\in I_{t,x}$.
\begin{enumerate}[(a)]
    \item If $y\in V$, there is $z\in V_3$ with $\hat T^{loc}_sf_{i,t,x}(y)=f_{i,t,x}(z)-c_s(y,z)$, and necessarily  $z\in I^+(y)$. In particular, $\hat T^{loc}_sf_{i,t,x}$ is lower semicontinuous on $V$.
    
    \item If $y\in V$ and $\hat T_s^{loc}f_{i,t,x}$ is differentiable at $y$, then there is a unique $z\in V_3$ with $\hat T^{loc}_sf_{i,t,x}(y)=f_{i,t,x}(z)-c_s(y,z)$, and it holds $z\in I^+(y)$ and
    \begin{align*}
        d_y(\hat T^{loc}_sf_{i,t,x})=\frac{\partial L}{\partial v}(y,\dot \gamma(0)),\ d_zf_{i,t,x}=\frac{\partial L}{\partial v}(z,\dot \gamma(s)),
    \end{align*}
    where $\gamma:[0,s]\to M$ is the unique(!) maximizing geodesic connecting $y$ to $z$.
    \item  The function $\hat T^{loc}_sf_{i,t,x}$ is continuous and even locally semiconvex on $V$.
\end{enumerate}   
\end{lemma}

\begin{proof}
    \begin{enumerate}[(a)]
    \item 
    Let $y\in V$. The existence of a point $z\in V_3$ with $\hat T^{loc}_sf_{i,t,x}(y)=f_{i,t,x}(z)-c_s(y,z)$ follows easily from the preceding corollary, $\overbar V_2\Subset V_3$, and the continuity of both $f_{i,t,x}$ on $V_3$ and $c_s(y,\cdot)$ on the closed set $J^+(y)$. To prove that necessarily $z\in I^+(y)$ for any optimal $z$, we first observe that $z\neq y$. Indeed, let $\gamma:[0,1]\to M$ be a maximizing geodesic connecting $x$ to $y$. Then $\gamma(1+\ep)$ is defined for small $\ep>0$, and it suffices to show that 
    \begin{align}
        f_{i,t,x}(\gamma(1+\ep))-c_s(y,\gamma(1+\ep))> f_{i,t,x}(y). \label{eqab'}
    \end{align} 
    To this aim, note that, by the smoothness of $f_{i,t,x}$, there is $C'>0$ such that, for small $\ep$,
    \begin{align}
        |f_{i,t,x}(\gamma(1+\ep))-f_{i,t,x}(y)| \leq C'\ep. \label{eqag'}
    \end{align}
    Moreover, using that $d(y,\gamma(1+\ep))\geq \ell_g(\gamma_{|[1,1+\ep]})\geq \ep d(x,y)$, we obtain
    \begin{align}
        -c_s(y,\gamma(1+\ep)) \geq \sqrt{s}(\ep d(x,y))^\frac 12. \label{eqah'}
    \end{align}
    Combining \eqref{eqag'} and \eqref{eqah'} yields \eqref{eqab'} for small $\ep$ (note that $d(x,y)>0$ since $U\times V\subseteq I^+$).
    
    Thus, $z\neq y$. If $z\in \partial J^+(y)\backslash \{y\}\subseteq I^+(x)$, one can adapt the argument in Corollary \ref{c1}, using a variation of a null geodesic from
    $y$ to $z$, to show that $z$ cannot attain the maximum in the definition of $\hat T^{loc}_sf_{i,t,x}(y)$. The lower semicontinuity follows as in Corollary \ref{c1}.
    
    \item 
    Let $y\in V$. If $z\in V_3$ is such that $\hat T^{loc}_sf_{i,t,x}(y)=f_{i,t,x}(z)-c_s(y,z)$, then $z\in I^+(y)$ by part (a), and it follows as\footnote{The only difference is that we consider a local version of the semigroup.} as in the proof of Lemma \ref{lf} that
    \begin{align*}
        d_y(\hat T^{loc}_sf_{i,t,x})=\frac{\partial L}{\partial v}(y,\dot \gamma(0)),\ d_zf_{i,t,x}=\frac{\partial L}{\partial v}(z,\dot \gamma(s)),
    \end{align*}
    where $\gamma:[0,s]\to M$ is a maximizing connecting $y$ to $z$. If there exists another curve $\tilde \gamma:[0,s]\to M$ such that $\hat T^{loc}_sf_{i,t,x}(y)=f_{i,t,x}(\tilde \gamma(s))-c_s(y,\tilde \gamma(s))$, then $\tilde \gamma$ must also be a maximizing timelike geodesic by part (a), and again it follows that 
    \begin{align*}
        d_y(\hat T^{loc}_sf_{i,t,x})=\frac{\partial L}{\partial v}(y,\dot {\tilde \gamma}(0)),\ d_{\tilde \gamma(s)}f_{i,t,x}=\frac{\partial L}{\partial v}(z,\dot {\tilde \gamma}(s)).
    \end{align*}
    However, since the Legendre transform is a diffeomorphism on $\op{int}(\C)$, it follows that $(\gamma(0),\dot \gamma(0))=(\tilde \gamma(0),\dot {\tilde \gamma}(0))$, which implies $\gamma=\tilde \gamma$. This proves both uniqueness of a curve and of $z$.
    
    \item 
    Let $y\in  V$ be arbitrary. For a sequence $y_k\to y$, we can find $z_k\in V_3\cap I^+(y_k)$ such that
    \begin{align*}
        f_{i,t,x}(z_k)-c_s(y_k,z_k)= \hat T^{loc}_sf_{i,t,x}(y_k).
    \end{align*}
    By Corollary \ref{cor3}, we have $z_k\in V_2$ for large $k$.
    Since $\overbar V_2$ is compact, it follows that, up to a subsequence, $z_k\to z\in \overbar V_2$ with $z\in J^+(y)$. Therefore,
    \begin{align*}
        \hat T^{loc}_sf_{i,t,x}(y)
        \geq
        f_{i,t,x}(z)-c_s(y,z)
        &=\lim_{k\to \infty}
        f_{i,t,x}(z_k)-c_s(y_k,z_k)
        \\
        &= \lim_{k\to \infty}\hat T^{loc}_sf_{i,t,x}(y_k).
    \end{align*}
    Hence, $\hat T^{loc}_sf_{i,t,x}$ is upper semicontinuous on $V$. Combining with the lower semicontinuity from part (a), this shows the desired continuity. To show local semiconvexity, using Corollary \ref{cor3}, the compactness of $\overbar V_2$, the established continuity and part (a), an easy compactness argument yields the existence of an open neighbourhood $V'\subseteq V$ of $y$ and $\delta>0$ such that, for any $y'\in V'$,
    \begin{align*}
        \hat T^{loc}_sf_{i,t,x}(y')
        =
        \sup\{f_{i,t,x}(z')-c_s(y',z')\mid z'\in V_2,\ d(y',z')\geq \delta \}.
    \end{align*}
    Using the precompactness of $V_2$, it follows as in the proof of Theorem \ref{thmb} that $\hat T^{loc}_sf_{i,t,x}$ is locally semiconvex on $V'$, and hence on $V$.
    \end{enumerate}
\end{proof}

\begin{lemma}\label{lk}
    There exists a constant $C_1>0$ and $s_2\in (0,s_1]$ such that, for all $s\in (0,s_2]$, $t\in I$, $x\in U$, $i\in I_{t,x}$ and $y\in V$, we have
    \begin{align*}
        \hat T^{loc}_sf_{i,t,x}(y)
        > \sup\{f_{i,t,x}(z)-c_s(y,z)\mid z\in V_3,\ d_{h}(y,z)\geq C_1\sqrt{s}\}.
    \end{align*}
\end{lemma}
\begin{proof}
    Set $C_1:=\frac{2M}{C}$, where 
    \begin{align*}
        M:=\sup\{|c_1(y,z)|\mid y,z\in \overbar V_2, z\in J^+( y)\}
    \end{align*}
    is finite thanks to the compactness of $\overbar V_2$. Now choose $s_2\in (0,s_1]$ with $C_1\sqrt{s_2}\leq d_h(V,\partial V_1)$. Fix $s,t,x,i$ and $y$ as in the lemma, and suppose by contradiction that $z_k\in V_3\cap J^+(y)$ is a sequence with $d_h(y,z_k)\geq C_1\sqrt{s}$ and
    \begin{align}
        f_{i,t,x}(z_k)-c_s(y,z_k) \geq \hat T^{loc}_sf_{i,t,x}(y)-\frac 1k\geq f_{i,t,x}(y)-\frac 1k. \label{eqaf}
    \end{align}
    By Corollary \ref{cor3}, $z_k\in V_2$ for all sufficiently large $k$. For these values of $k$, let $\gamma_k:[0,s]\to M$ be a maximizing geodesic connecting $y$ to $z_k$. Then by Proposition \ref{pro1}(vi), $\gamma(\tau)\in V_3$ for all $\tau\in [0,s]$. In particular, Proposition \ref{pro1}(vii) implies $\frac{d}{d\tau} (f_{i,t,x}\circ \gamma_k)(\tau)\leq 0$ for all $\tau$ and $\frac{d}{d\tau} (f_{i,t,x}\circ \gamma_k)(\tau)\leq -C|\dot \gamma_k(\tau)|_h$ as long as $\gamma_k(\tau)\in V_1$. In particular, since $C_1\sqrt{s}\leq C_1\sqrt{s_2}\leq d_h(V,\partial V_1)$, we have
    \begin{align*}
        f_{i,t,x}(z_k)-f_{i,t,x}(y)=\int_0^s \frac{d}{d\tau} (f_{i,t,x}\circ \gamma)(\tau)\, d\tau
        \leq
        -CC_1\sqrt{s}.
    \end{align*}
    Thus,
    \begin{align*}
        f_{i,t,x}(z_k)-c_s(y,z_k)
        \leq 
        f_{i,t,x}(y)-CC_1\sqrt{s}+\sqrt{s}M.
    \end{align*}
    For large $k$, this is a contradiction to \eqref{eqaf} by the definition of $C_1$.
\end{proof}

The next theorem is well known in the context of Tonelli-Hamiltonian systems and first appeared in \cite{Bernard} in the compact setting. An adaption of the standard proof in the non-compact case can be found in Appendix B of \cite{Fathi/Figalli/Rifford}. The extension of this result to the Lorentzian setting follows the same approach.

\begin{theorem}\label{thmd}
    Let $s_2$ be as in the previous lemma. There exists $s_0\in (0,s_2]$ such that, for any $s\in [0,s_0]$, $t\in I$, $x\in U$, and $i\in I_{t,x}$, the map 
    \[
    \psi_{i,t,x,s}:V_1\to M,\ \psi_{i,t,x,s}(z):=\pi^*\circ \psi_{-s}(z,d_zf_{i,t,x}),
    \]
    where $\pi^*:T^*M\to M$ is the canonical projection, is well-defined, a homeomorphism onto its image which contains $V$, and satisfies
\begin{align*}
    \hat T_s^{loc}f_{i,t,x}(\psi_{i,t,x,s}(z))=f_{i,t,x}(z)-
    c_s(\psi_{i,t,x,s}(z),z),
\end{align*}
 whenever $\psi_{i,t,x,s}(z)\in V$. Moreover, the family of maps 
 \[
 \{\hat T_s^{loc}f_{i,t,x}\mid s\in [0,s_0], t\in I, x\in U, i\in I_{t,x}\}
 \]
 is uniformly locally semiconcave on $V$.
\end{theorem}
\begin{proof}
    In proving the theorem, it suffices to consider the cases $s\in (0,s_0]$ (with $s_0$  to be defined), since by Proposition \ref{pro1}(i), all statements also hold when adding $s=0$ (note that $C^2$-boundedness implies uniform local semiconcavity \cite{Fathi}).

    We will work both in local coordinates w.r.t.\ the chart $\phi$, and intrinsically on the manifold. To distinguish between these settings, we will often use a tilde to indicate the local coordinate version of an object.

    Define the index set $J:=\{(i,t,x)\mid t\in I,x\in U,i\in I_{t,x}\}$. By Proposition \ref{pro1}(ii), the set
    \[
    C:={\{(z,d_zf_j)\mid z\in V_1, j\in J\}}
    \]
    is relatively compact in $\op{int}(\C^{*})$. Therefore, there exists a precompact open neighbourhood $O\Subset T^*V_3$ of $\overbar C$, and a time $T>0$  such that the Hamiltonian flow $\psi$ is well-defined on $(-2T,2T)\times O$ and takes values in $T^*V_3$. Let $\tilde O:=T^*\phi(O)$, where $T^*\phi$ denotes the contangent bundle chart associated with $\phi$.
    Then we can consider the map
    \begin{align*}
         F:(-2T,2T)\times \tilde O\to \R^n,\ (s,\tilde z,\tilde p)\mapsto(\phi\circ \pi^*\circ \psi_{-s})((T^*\phi)^{-1}(\tilde z,\tilde p))-\tilde z.
    \end{align*}
     This map $F$ is smooth and  satisfies $F(0,\tilde z,\tilde p)=0$. Consequently, there exists a modulus of continuity $\omega_F:[0,\infty)\to [0,\infty)$ such that, for any $s\in [-T,T]$, the map $F(s,\cdot,\cdot)$ has Lipschitz constant $\omega_F(s)$ on the compact set $\tilde C:=T^*\phi(C)\subseteq  \tilde O$. 
     
     We now define $0<s_0\leq \min\{s_2,T\}$ such that
     \begin{align*}
        1-(K+1)\omega_F(s_0)>0 \text{ and } C_1\sqrt{s_0}\leq d_{h}(V,\partial V_1)/2,
     \end{align*}
      where $K$ bounds the family $(f_j\circ \phi^{-1})_j$ in $C^2(W_3)$ (Proposition \ref{pro1}(i)), and $C_1$ is the constant from the previous lemma.
    \medskip 

    Now, fix $s\in (0,s_0]$ and $j\in J$. Note that $\psi_{j,s}$ is well-defined on $V_1$ by definition of $s_0, T$ and $C$. Since $s_0\leq s_2\leq s_1$, Lemma \ref{lj} guarantees that $\hat T^{loc}_sf_j$ is locally semiconvex on $V$, and thus differentiable almost everywhere on $V$. 
    Furthermore, Lemma \ref{lj}(b) states that for every differentiability point $y\in V$ of $\hat T^{loc}_sf_j$, there exists a unique $z\in V_3$ such that
    \begin{align}
        \hat T^{loc}_sf_j(y)=f_j(z)-c_s(y,z), \label{j}
    \end{align}
    and necessarily $z\in I^+(y)$. Moreover, by Lemma \ref{lk} and the choice of $s_0$, we have $z\in V_1$ and $d_{h}(z,\partial V_1)\geq d_h(V,\partial V_1)/2$. We claim that $y=\psi_{j,s}(z)$.
 
     Indeed, Lemma \ref{lj}(b) gives
    \begin{align*}
        d_y(\hat T^{loc}_sf_j)=\frac{\partial L}{\partial v}(y,\dot \gamma(0)) \text{ and } d_zf_j=\frac{\partial L}{\partial v}(z,\dot \gamma(s)),
    \end{align*}
     where $\gamma:[0,s]\to M$ is the unique maximizing geodesic connecting $y$ to $z$. In particular, 
    \begin{align}
        (y,d_y(\hat T^{loc}_sf_j))=\bigg(\gamma(0),\frac{\partial L}{\partial v}(\gamma(0),\dot \gamma(0))\bigg) 
        &={\cal L}(\gamma(0),\dot \gamma(0)) \nonumber
        \\
        &={\cal L}(\phi_{-s}(\gamma(s),\dot \gamma(s))) \nonumber
        \\
        &=\psi_{-s}\bigg(z,\frac{\partial L}{\partial v}(z,\dot \gamma(s))\bigg) \nonumber
        \\
        &=\psi_{-s}(z,d_zf_j). \label{eqy}
    \end{align}
    It follows that $y=\psi_{j,s}(z)$, as claimed. 
 
    Now let $\tilde \psi_{j,s}:=\phi\circ \psi_{j,s}\circ (\phi_{|V_1})^{-1}:W_1:=\phi(V_1)\to W_3$ be the map $\psi_{j,s}$ in local coordinates. 
    Since the map $W_3\ni z\mapsto (z,D(f_j\circ \phi^{-1})(z))$ has Lipschitz constant $K+1$ thanks to the convexity of $V$ and the uniform $K$-boundedness of the second derivatives, and since $s\leq T$, the mapping
    \begin{align*}
         W_1\to \R^n,\ \tilde z\mapsto\tilde \psi_{j,s}(\tilde z)-\tilde z=F(s,\tilde z,D(f_j\circ \phi^{-1})(\tilde z)),
    \end{align*}
    has Lipschitz constant $(K+1)\omega_F(s)$. Thus, for any $\tilde z,\tilde z'\in W_1$, we have
    \begin{align}
        |\tilde \psi_{j,s}(\tilde z)-\tilde \psi_{j,s}(\tilde z')|
        \geq
        |\tilde z-\tilde z'|-(K+1)\omega_F(s)|\tilde z-\tilde z'|. \label{eqv}
    \end{align}
    By definition of $s_0$, this shows that $\tilde \psi_{j,s}$, and hence $\psi_{j,s}$, is a homeomorphism onto its image with locally Lipschitz inverse.
    
    Now, if $\hat T^{loc}_sf_j$ is differentiable at $y\in V$ and $z$ is such that \eqref{j} holds, we already saw that $z\in V_1\cap I^+(y)$, $d_h(z,\partial V_1)\geq d_h(V,\partial V_1)/2$ and $y=\psi_{j,s}(z)$. In particular, $\psi_{j,s}(V_1)$ contains a set of full measure in $V$, namely all the differentiability points. It follows that $\psi_{j,s}(V_1)\supseteq V$: Indeed, if $y_k\in V$ is a sequence of differentiability points converging to $y\in V$, let $z_k$ be such that \eqref{j} holds (with $y,z$ replaced by $y_k,z_k$). Then, thanks to the local Lipschitz continuity of $\psi_{j,s}^{-1}$ and $d_h(z_k,\partial V_1)\geq d_h(V,\partial V_1)/2$, $z_k$ converges to some $z\in V_1$, and the continuity of $\psi_{j,s}$ implies $\psi_{j,s}(z)=y$. Thus, we indeed have $\psi_{j,s}(V_1)\supseteq V$.

    Moreover, denoting by $\Gamma_f(y):=(y,d_yf)$ the graph of the derivative of a smooth function $f$, \eqref{eqy} shows for any differentiability point $y\in V$ that
    \begin{align*}
        \Gamma_{\hat T^{loc}_sf_j}(y)=\psi_{-s}\circ \Gamma_{f_j}\circ \psi_{j,s}^{-1}(y)
    \end{align*}
    Hence, if $\tilde \psi$ denotes the Hamiltonian flow in local coordinates, we have for $\tilde y=\phi(y)$, the differentiability point in local coordinates,
    \begin{align*}
        \Gamma_{(\hat T^{loc}_sf_j)\circ \phi^{-1}}(\tilde y)=
        \tilde \psi_{-s} \circ \Gamma_{f_j\circ \phi^{-1}}\circ \tilde \psi_{j,s}^{-1}(\tilde y).
    \end{align*}
    The right-hand side of the above equation is well-defined on $W:=\phi(V)$ and Lipschitz continuous with Lipschitz constant 
    \begin{align*}
        \op{Lip}(\tilde \psi_{|[-T,T]\times \tilde C}) (K+1) (1-(K+1)\omega_F(s_0))^{-1}
    \end{align*}
    which is independent of $s$ and $j$. In particular,  $\hat T^{loc}_sf_j\circ \phi^{-1}$ is locally semiconvex on $W$ and its a.e.\ derivative admits a Lipschitz extension on $V$. Then  $\hat T_s^{loc}f_j\circ \phi^{-1}$ must be $C^{1,1}$ on $W$ and the Lipschitz constant of the derivative is independent of $j$ and $s\in (0,s_2]$. 
    In particular, being $W$ convex, it follows that the family $\{\hat T^{loc}_sf_{j}\circ \phi^{-1}\mid s\in (0,s_2],j\in J\}$ is uniformly semiconcave on $W$ \cite{Fathi}. As a consequence, by definition, $\{\hat T^{loc}_sf_{j}\mid s\in (0,s_2],j\in J\}$ is uniformly locally semiconcave on $V$. This shows the last part of the lemma.

    Finally, if $s\in (0,s_0]$, $j\in J$ and $z\in V_1$ are such that $\psi_{j,s}(z)\in V$, then $\psi_{j,s}(z)\in V$ is a differentiability point of $\hat T^{loc}_sf_j$, and our first claim in the proof shows that
    \begin{align*}
        \hat T^{loc}_sf(\psi_{j,s}(z)) = f(z)-c_s(\psi_{j,s}(z),z).
    \end{align*}
    This concludes the proof.
\end{proof}

\begin{proposition}
    Let $s_0$ be as in the previous theorem, and let $s\in [0,s_0]$, $t\in I$ and $x\in U$. We have
    \begin{align*}
       \inf_{i\in I_{t,x}} \hat T^{loc}_sf_{i,t,x}=\hat T_sT_{t}\chi_x \text{ on } V. 
    \end{align*}
\end{proposition}
\begin{proof}
    For $s=0$ the result follows immediately from Proposition \ref{pro1}(iii). Thus, suppose $s\neq 0$, and let $y\in V$. By Theorem \ref{thmb} and $C_0\sqrt{s}\leq C_0\sqrt{s_0'}\leq d_h(V,\partial V_1)$ (see Remark \ref{rema}), there exists $z\in V_1$ with 
    \begin{align*}
        \hat T_sT_{t}\chi_x(y)=T_{t}\chi_x(z)-c_s(y,z),
    \end{align*}
    and necessarily $z\in I^+(y)$.
    If $\gamma:[0,s]\to M$ is a necessarily timelike maximizing geodesic with $\gamma(0)=y$ and $\gamma(s)=z$, Lemma \ref{l4} implies that
    \begin{align}
        p:=\frac{\partial L}{\partial v}(z,\dot \gamma(s))\in \partial^+(T_{t}\chi_x)(z). \label{eqag}
    \end{align}
    Therefore, Proposition \ref{pro1}(iv) ensures the existence of $i\in I_{t,x}$ with $f_{i,t,x}(z)=T_{t}\chi_x(z)$ and $d_zf_{i,t,x}=p$. Since $\gamma$ is a geodesic with ${\cal L}(\gamma(s),\dot \gamma(s))=(z,p)$, we have $\gamma(\tau)=\pi^*\circ \psi_{-s+\tau}(z,p)=\psi_{i,t,x,s-\tau}(z)$ for $\tau\in [0,s]$, and in particular, $\psi_{i,t,x,s}(z)=\gamma(0)=y\in V$.
       Hence, by the previous theorem, we obtain
    \begin{align*}
        \hat T^{loc}_sf_{i,t,x}(y)
        =
        \hat T^{loc}_sf_{i,t,x}(\psi_{i,t,x,s}(z))
        =
        f_{i,t,x}(z)-c_s(y,z)
        &=
        T_{t}\chi_x(z)-c_s(y,z)
        \\
        &=
        \hat T_sT_{t}\chi_x(y).
    \end{align*}
    This shows that
    \begin{align*}
       \inf_{i\in I_{t,x}} \hat T_sf_{i,t,x}(y)
       \leq \hat T_sT_{t}\chi_x(y). 
    \end{align*}
    However, Lemma \ref{li} shows that equality holds. Since $y$ was arbitrary, this proves the proposition.
\end{proof}

\begin{proof}[Proof of Theorem \ref{thmc}]
Let $U,V$ and $I$ be as in Remark \ref{rema}, and let $s_0$ be as in Theorem \ref{thmd}. For the first part, note that, by Theorem \ref{thmd}, the family $\{\hat T^{loc}_sf_{i,t,x}\mid s\in [0,s_0],t\in I,x\in U,i\in I_{t,x}\}$ is uniformly locally semiconcave on $V$. Moreover, from the previous proposition, we have
\begin{align*}
       \inf_{i\in I_{t,x}} \hat T^{loc}_sf_{i,t,x}=\hat T_sT_{t}\chi_x \text{ on } V.
    \end{align*}
    It follows that also the family $\{\hat T_sT_{t}\chi_x\mid s\in [0,s_0],t\in I,x\in U\}$ is uniformly locally semiconcave on $V$ (\cite{Fathi/Figalli}, Theorem A.11).
\end{proof}

\section{The main result}

In this section we will prove Theorems \ref{aa} and \ref{b}. The following theorem (compare with Claims 4.7, 4.9 and 4.10 in \cite{Cannarsa/Cheng/Fathi}) summarizes all the important results from the last two sections that will be needed in the proofs. Its proof is an easy consequence from Theorems \ref{thmb} and \ref{thmc}.

\begin{theorem}\label{main2}
    Let $x_0\in M$ and $y_0\in I^+(x_0)$. Then there exist two open neighbourhoods $U$, $V$ of $x_0$ and $y_0$, respectively, with $U\times V\subseteq I^+$, some number $s_0>0$ and a constant $C_0>0$ such that:
   \begin{enumerate}[(a)]
   \item The mapping $[0,s_0]\times U\times V\to \R,\ (s,x,y)\mapsto \hat T_sT_{1+s}\chi_x(y)$, is continuous.
   \item If $x\in U$, $y\in V$ and $s\in [0,s_0]$, there exists a unique $z\in M$ with
   \begin{align*}
        \hat T_sT_{1+s}\chi_x(y)
        =
        T_{1+s}\chi_x(z)-c_s(y,z).
    \end{align*}
    Moreover, $d_{h}(y,z)\leq C_0\sqrt{s}$, and if $s>0$, then necessarily $z\in I^+(y)$. 
    \end{enumerate}
\end{theorem}
\begin{proof}
    We choose $U,V$, $s_0\leq 1/2$ and $C_0$ such that Theorem \ref{thmb} holds (with $t_0=1$), and such that $\hat T_sT_{1+s}\chi_x$ is locally semiconcave on $V$ for all $x\in U$ and $s\in [0,s_0]$ (see Theorem \ref{thmc}). Then part (a) follows immediately from Theorem \ref{thmb}, using the fact that $1+s\leq 3/2=3t_0/2$. Part (b) follows from Theorem \ref{thmb}, except for the uniqueness of $z$. However, uniqueness is obvious when $s=0$. Now suppose $s>0$ and that there exist two such $z,\tilde z\in M$. Then necessarily $z,\tilde z\in I^+(y)$. Let $\gamma,\tilde \gamma:[0,s]\to M$ be two maximizing geodesics connecting $y$ to $z$ and $\tilde z$, respectively. From Lemma \ref{lf},
    \[
    \frac{\partial L}{\partial v}(y,\dot \gamma(0)),\frac{\partial L}{\partial v}(y,\dot {\tilde \gamma}(0))\in \partial^- (\hat T_sT_{1+s}\chi_x)(y).
    \]
    Since $\hat T_sT_{1+s}\chi_x$ is both locally semiconcave and locally semiconvex on $V$ thanks to Theorem \ref{thmb} and Theorem \ref{thmc}, it is differentiable (and even $C^1$) on $V$ \cite{Fathi}. It follows that 
    \[
    \frac{\partial L}{\partial v}(y,\dot \gamma(0))=\frac{\partial L}{\partial v}(y,\dot {\tilde \gamma}(0)).
    \]
    In particular, $\gamma(0)=\tilde \gamma(0)$ and since the Legendre transform is a diffeomorphism, also $\dot \gamma(0)=\dot {\tilde \gamma}(0)$. Hence, $\gamma=\tilde \gamma$, implying $z=\tilde z$. This concludes the proof.
\end{proof}

For the proofs of both Theorems \ref{aa} and \ref{b}, the cut locus plays a crucial role. For instance, we will prove Theorem \ref{b}(a) by first showing that $\op{Cut}_M$ (see below) is a strong deformation retract of $J^+\backslash {\cal A}$, and that the inclusion ${\cal NU}(M,g)\hookrightarrow \op{Cut}_M$ is a homotopy equivalence.

Let us recall the definition of the cut locus, along with a powerful characterization. We will also revisit some important results concerning the compactness of maximizing geodesics. Recall that we are assuming $M$ to be globally hyperbolic.

\begin{definition}\rm\label{defb}
\begin{enumerate}[(a)]
    \item Let $x\in M$ and let $\gamma:[0,a)\to M$, $a\in (0,\infty]$, be a future inextendible causal geodesic starting at $x$. Set
    \begin{align*}
        t_0:=\sup\{t\in [0,a)\mid d(x,\gamma(t))=\ell_g(\gamma_{|[0,t]})\}\in [0,a].
    \end{align*}
    Then $t_0>0$\footnote{This is not true in general, but it is in our case since $M$ is globally hyperbolic} and if $t_0<a$, the point $\gamma(t_0)$ is called the \emph{cut point of $x$ along $\gamma$}.

    \item A point $y$ is called \emph{causal/timelike/null cut point of $x$} if $y$ is the cut point of $x$ along $\gamma$ for some causal/timelike/null geodesic $\gamma:[0,a)\to M$ emerging from $x$. 

    \item The \emph{causal (resp.\ timelike/null) cut locus} $\op{Cut}_M(x)$ (resp.\ $\op{Cut}_M^t(x),\ \op{Cut}_M^n(x)$) is defined as the set of all causal (resp.\ timelike, null) cut points of $x$.

    \item The set $\op{Cut}_M\subseteq M\times M$ (resp.\ $\op{Cut}^t_M,\op{Cut}_M^n$) is defined as the set of all $(x,y)\in M\times M$ such that $y\in \op{Cut}_M(x)$ (resp.\ $y\in \op{Cut}_M^t(x)$, $y\in \op{Cut}_M^n(x)$).
\end{enumerate}
\end{definition}

\begin{remark}\rm
    Let $\gamma:I\to M$ be a future inextendible causal geodesic and $t_0 \in I$. Set 
    \begin{align*}
        t_1:=\sup\{t\in I\cap [t_0,\infty)\mid d(\gamma(t_0),\gamma(t))=\ell_g(\gamma_{|[t_0,t_1]})\}\in [t_0,\sup I].
    \end{align*}
    If $t_1<\sup I$, we also say that $\gamma(t_1)$ is the cut point of $\gamma(t_0)$ along $\gamma$. Note that this definition of cut points coincides with the definition introduced above.
\end{remark}

\begin{definition}\rm
Let $\gamma:I\to M$ be a geodesic.
    \begin{enumerate}[(a)]
        \item A \emph{Jacobi field along $\gamma$} is a smooth vector field $J$ along $\gamma$ such that
        \begin{align*}
            \frac{D^2J}{dt^2}+R(J,\dot \gamma)\dot \gamma=0,
        \end{align*}
        where 
        \[
        R(X,Y)Z:=\nabla_X \nabla_Y Z-\nabla_Y \nabla_X Z-\nabla_{[X,Y]}Z,
        \]
        is the Riemannian curvature tensor.
        \item 
        A point $\gamma(t_1)$, $t_1\in I$, is said to be a \emph{conjugate point} of $\gamma(t_0)$, $t_0\in I$, along $\gamma$ if there exists a non-zero Jacobi field $J$ along $\gamma$ such that $J(t_0)=J(t_1)=0$.
    \end{enumerate}
\end{definition}

\begin{remark}\rm \label{hdsdoa}
    Let $\gamma:[a,b]\to M$ be a geodesic. Then $\gamma(b)$ is a conjugate point of $\gamma(a)$ along $\gamma$ if and only if $(b-a)\dot \gamma(a)$ is a critical point of the exponential map $\exp_{\gamma(a)}$ (Proposition 3.5 in \cite{DoCarmo}).
\end{remark}

There is the following characterization of cut points. It shows ${\cal NU}(M,g)\subseteq \op{Cut}_M$.
\begin{theorem}\label{thme}
    Let $\gamma:[0,a)\to M$ be a future inextendible causal geodesic emerging from $x$. If $y=\gamma(t_0)\in J^+(x)$ is the cut point of $x$ along $\gamma$, then at least one of the following hold:
    \begin{enumerate}[(i)]
        \item $y$ is the first conjugate point of $x$ along $\gamma$.
        \item There exists another distinct maximizing geodesic connecting $x$ to $y$.
    \end{enumerate}
    Conversely, if $y=\gamma(t_0)$ is a conjugate point of $x$ along $\gamma$ or (ii) holds, then $\gamma$ ceases to be maximizing beyond $t=t_0$.
\end{theorem}
\begin{proof}
    The implication $\Rightarrow$ is established by Theorems 9.12 and 9.15 in \cite{Ehrlich}. The fact that (ii) implies $\gamma$ is not maximizing beyond $t_0$ follows from Corollaries 9.4 and 9.11 in \cite{Ehrlich}; although the proofs are omitted there, this conclusion essentially follows from \cite{Minguzzi}, Theorem 2.9, which states that any maximizing causal curve must be a pregeodesic.
    The fact that $\gamma$ is not maximizing beyond a conjugate point follows from Theorem 9.10 in \cite{Ehrlich} for the timelike case (see Proposition 10.12 in \cite{Ehrlich} for a proof) and Theorem 10.72 in \cite{Ehrlich} for the lightlike case. For the lightlike case, see also Proposition 2.2.8 in \cite{Baer}.
\end{proof}

\begin{definition}\rm
    We define the function
    \begin{align*}
        \alpha:\C\to [0,\infty],\ \alpha(x,v):=\sup\{t\geq 0\mid d(x,\exp_x(tv))=t|v|_g\}.
    \end{align*}
    Note that, if $\exp_x(tv)$ is only defined for $t\in [0,a)$, then $\alpha(x,v)\leq a$.
\end{definition}

\begin{lemma}\label{lq}
      $\alpha$ is continuous at $(x,v)$ unless $\alpha(x,v)$ is finite and $\exp_x(\alpha(x,v)v)$ does not exist.
\end{lemma}
\begin{proof}
    See \cite{Ehrlich}, Proposition 9.33.
\end{proof}

The following lemma is well-known.

\begin{lemma}\label{lm}
    Let $x_k\in M$ and $v_k\in \C_{x_k}$ such that the curves $[0,1]\ni t\mapsto \exp_{x_k}(tv_k)$ are maximizing geodesics connecting $x_k$ to $y_k:=\exp_{x_k}(v_k)$. Suppose that $x_k\to x$ and $y_k\to y$. Then $(x,y)\in J^+$, and $(x_k,v_k)$ converges, along a subsequence, to some $(x,v)\in \C$. Moreover, $[0,1]\ni t\mapsto \exp_x(tv)$ is well-defined and a maximizing geodesic connecting $x$ to $y$.
\end{lemma}
\begin{proof}
    Since $J^+$ is closed, $(x,y)\in J^+$. The case $x\ll y$ follows from Lemma 9.6 in \cite{Ehrlich}, together with the continuity of $d$ and a scaling argument. The case $x=y$ is trivial: for large $k$, global hyperbolicity implies $(x_k,v_k)=(x_k,\exp_{x_k}^{-1}(y_k))\xrightarrow{k\to \infty} (x,0)$. If $x\neq y$ and $d(x,y)=0$, the claim follows from Lemma 9.14 and the proof of Lemma 9.25 in \cite{Ehrlich}.
\end{proof}

\begin{corollary}\label{cor1}
    Let $(x,y)\in \op{Cut}_M$. Then for $\ep>0$ sufficiently small there exist two open neighbourhoods $U$ of $x$ and $V$ of $y$ such that the following holds:
    \begin{itemize}
        \item Every maximizing causal geodesic $\gamma:[0,1]\to M$ with $\gamma(0)\in U$ and $\gamma(1)\in V$ can be extended to a geodesic parametrized over $[0,1+\ep]$. This extension, however, is never maximizing.
    \end{itemize}

\end{corollary}
\begin{proof}
    Let $U'$ and $V'$ be two arbitrary precompact neighbourhoods of $x_0$ and $y_0$, respectively. By the closedness of $J^+$ and the preceding lemma, the set of maximizing causal geodesics $\gamma:[0,1]\to M$ with $\gamma(0)\in \overbar U'$ and $\gamma(1)\in \overbar V'$ is compact in the $C^1$-topology. This compactness ensures the existence of $\ep_0>0$ such that any such geodesic can be extended to $[0,1+\ep_0]$. 
    
    Let $0< \ep\leq \ep_0$ and suppose, for contradiction, that there exists a sequence $(x_k,y_k)\in J^+$ converging to $(x,y)$, as well as maximizing geodesics $\gamma_k:[0,1]\to M$ connecting $x_k$ to $y_k$, each extendable to a maximizing geodesic on $[0,1+\ep]$. Write $\gamma_k(t)=\exp_{x_k}(tv_k)$ for some $v_k\in \C_{x_k}$. Lemma \ref{lm} guarantees, along a subsequence, $(x_k,v_k)\to (x,v)\in \C$, and that $[0,1]\ni t\mapsto \exp_x(tv)$ is a maximizing geodesic connecting $x$ to $y$. Since $(x,y)\in \op{Cut}_M$, we have $\alpha(x,v)=1$. Hence, by continuity of $\alpha$, $\alpha(x_k,v_k)<1+\ep$ for large $k$, contradicting the assumption.
\end{proof}

From now on, we will abbreviate ${\cal NU}:={\cal NU}(M,g)$ and ${\cal NU}^t:={\cal NU}^t(M,g)$.

\begin{theorem}\label{ab}
    Let $(x_0,y_0)\in \op{Cut}_M^t$ and $V_0$ be an open neighbourhood of $y_0$. Then there exist two open neighbourhoods $U,V$ of $x_0$ and $y_0$, respectively, with $U\times V\subseteq I^+$, and $s_0>0$ such that the map
    \begin{align*}
        F:[0,s_0]\times U\times V\to V_0, \ F(s,x,y):=z,
    \end{align*}
    where $z\in M$ is the unique(!) point satiyfying $\hat T_sT_{1+s}\chi_x(y)=T_{1+s}\chi_x(z)-c_s(y,z)$, is well-defined, continuous and satisfies
        \begin{enumerate}[(i)]
        \item $F(s,x,y)\in I^+(x)$ for all $(s,x,y)$,
            \item $F(0,x,y)=y$ for all $(x,y)$,
            \item $(x,F(s,x,y))\in {\cal NU}$ for all $s>0$ and $(x,y)\in (U\times V)\cap \op{Cut}_M^t$,
            \item $(x,F(s_0,x,y))\in  {\cal NU}$ for all $(x,y)$.
        \end{enumerate}
\end{theorem}
\begin{proof}
      Let $U$, $V$, $s_0$ and $C_0$ be given as in Theorem \ref{main2}. Without loss of generality, we may assume that $C_0\sqrt{s_0}<d_{h}(V,\partial V_0)$. Since $(x_0,y_0)\in \op{Cut}_M^t$, by Corollary \ref{cor1} we may also assume that, if $(x,y)\in U\times V$, no maximizing geodesic $\gamma:[0,1]\to M$ connecting $x$ to $y$ can be extended to a maximizing geodesic on $[0,1+s_0]$.

    By Theorem \ref{main2}, for $(x,y)\in U\times V$ and $s\in [0,s_0]$, there exists a unique $z\in M$ with 
     \begin{align*}
        \hat T_sT_{1+s}\chi_x(y)
        =
        T_{1+s}\chi_x(z)-c_s(y,z).
    \end{align*}
    Moreover, $z\in J^+(y)$ and $d_{h}(y,z)\leq C_0\sqrt{s_0}$, so $z\in V_0$. Therefore, the map $F$ is well-defined. 
    \bigskip
    
    \noindent \textbf{Claim:} $F$ is continuous. 
    \medskip 
    
    \noindent \textbf{Proof of claim:} 
    Let $(s_k,x_k,y_k)\in [0,s_0]\times U\times V$ be any sequence converging to $(s,x,y)\in [0,s_0]\times U\times V$. Thanks to Theorem \ref{main2}, the sequence $z_k:=F(s_k,x_k,y_k)$ is precompact and converges, along a subsequence $z_{k_l}$, to some $z\in J^+(y)\subseteq J^+(I^+(x))=I^+(x)$. Using part (a) of Theorem \ref{main2} and the continuity of $\C$ on $(0,\infty)\times J^+$ (if $s>0$) and at $(0,y,y)$ (if $s=0$. Note that $z_k\to y$ in this case), we have
    \begin{align}
        \hat T_{s}T_{1+s}\chi_x(y)
        =
        T_{1+s}\chi_x(z)-c_s(y,z). \label{eqm}
    \end{align}
    The uniqueness part in Theorem \ref{main2} gives $z=F(s,x,y)$. This shows that $F(s_{k_l},x_{k_l},y_{k_l})\to F(s,x,y)$. However, the sequence $(s_k,x_k,y_k)$ was arbitrary, implying that $F(s_k,x_k,y_k)\to F(s,x,y)$.
    \hfill \checkmark 
    \bigskip

    Part (i) follows from $F(s,x,y)\in J^+(y)\subseteq I^+(x)$. Property (ii) is immediate. 
    
    To prove (iii) and (iv), let $(x,y)\in U\times V$ and $s\in (0,s_0]$. Setting $z:=F(s,x,y)$, by definition we have
    \[
    \hat T_{s}T_{1+s}\chi_x(y)=T_{1+s}\chi_x(z)-c_{s}(y,z).
    \]
    Let $\gamma:[0,s]\to M$ be a (necessarily timelike, by Theorem \ref{main2}) maximizing geodesic connecting $y$ to $z$. From Lemma \ref{lf} we deduce 
    \begin{align}
        \frac{\partial L}{\partial v}(z,\dot \gamma(s))\in \partial^+ (T_{1+s}\chi_x)(z)
        \text{ and } 
        \frac{\partial L}{\partial v}(y,\dot \gamma(0))\in \partial^- (\hat T_{s}T_{1+s}\chi_x)(y). \label{eqww}
    \end{align}

     Suppose there exists a unique maximizing geodesic ${\tilde \gamma}:[0,1+s]\to M$ connecting $x$ to $z$. Since $z\in  I^+(x)$, and since unique super-differentiability implies differentiability for locally semiconcave functions (\cite{Villani}, Theorem 10.8), Lemma \ref{l4} implies that $T_{1+s}\chi_x$ is differentiable at $z$ with
    \begin{align}
        d_zT_{1+s}\chi_x=\frac{\partial L}{\partial v}(z,\dot {\tilde \gamma}(1+s)). \label{eqxx}
    \end{align}
    This derivative must also be the unique superdifferential. Hence, since the Legendre transform is a diffeomorphism, we obtain from \eqref{eqww} and \eqref{eqxx} that
    \begin{align*}
       (\gamma(s),\dot \gamma(s)) =({\tilde \gamma}(1+s),\dot {\tilde \gamma}(1+s)). 
    \end{align*}
    Since both $\gamma$ and ${\tilde \gamma}$ are maximizing geodesics, it follows that $\gamma(\tau)={\tilde \gamma}(1+\tau)$ for any $\tau\in [0,s]$. In particular, ${\tilde \gamma}(1)=\gamma(0)=y$, hence ${\tilde \gamma}$ is a maximizing geodesic defined on $[0,1+s]$ with ${\tilde \gamma}(0)=x$ and ${\tilde \gamma}(1)=y$.

    To prove (iii), let $(x,y)\in \op{Cut}_t(M)$. By the preceding lemma, no maximizing geodesic connecting $x$ to $y$ can be extended to a maximizing geodesic beyond $[0,1]$. Clearly, this contradicts the fact that ${\tilde \gamma}$ is maximizing on $[0,1+s]$
    
    To prove (iv), let $s=s_0$. Then ${\tilde \gamma}$ is maximizing on $[0,1+s_0]$, which is a contradiction to the construction of $U$ and $V$.
\end{proof}

\begin{corollary}\label{cor4}
    Let $\ep:\op{Cut}_M^t\to (0,\infty)$ be a continuous function. Then there exists a continuous function $s:\op{Cut}_M^t\to (0,\infty)$ such that the map
    \begin{align*}
        \bar F:[0,1]\times \op{Cut}_M^t\to M,\ \bar F(t,x,y):= z,
    \end{align*}
    where $z\in M$ is the unique point satiyfying $\hat T_{ts(x,y)}T_{1+ts(x,y)}\chi_x(y)=T_{1+ts(x,y)}\chi_x(z)-c_{ts(x,y)}(y,z)$, is well-defined, continuous and satisfies
    \begin{enumerate}[(i)]
        \item $\bar F(t,x,y)\in I^+(x)$ for all $(t,x,y)$,
        \item $\bar F(0,x,y)=y$ for all $(x,y)$
        \item $(x,\bar F(t,x,y))\in {\cal NU}$ for all $t>0$ and $(x,y)$,
        \item $d_h(\bar F(t,x,y),y)\leq \ep(x,y)$ for all $(t,x,y)$.
    \end{enumerate}
\end{corollary}
\begin{proof}
    For each $(x,y)\in \op{Cut}_M^t$, let $U(x,y)\subseteq M$ and $V(x,y)\subseteq B_{\ep(x,y)/4}(y)\subseteq M$ be open neighbourhoods of $x$ and $y$, respectively, and let $s_0(x,y)>0$ such that the statement of the above theorem holds with $V_y:=B_{\ep(x,y)/4}(y)$. Since $\ep$ is continuous, we may assume that $\ep(x,y)\leq 2\ep(x',y')$ for all $(x',y')\in U(x,y)\times V(x,y)$.

    Let $(x_k,y_k)\in \op{Cut}_M^t$ be a countable family of points such that
    \begin{align*}
        \bigcup_{k\in \N} U(x_k,y_k)\times V(x_k,y_k)\supseteq \op{Cut}_M^t.
    \end{align*}
    Let $\rho_k$ be a locally finite smooth partition of unity for this union, subordinate to this open cover, and define 
    \[
    s:\op{Cut}_M^t\to (0,\infty),\ s(x,y):=\sum_{k\in \N} \rho_k(x,y)s_0(x_k,y_k).
    \]
    We claim that all properties hold for this function. 
    
    Clearly, $s$ is continuous. Fix $(x_0,y_0)\in \op{Cut}_M^t$. There exist two open neighbourhoods $U$ and $V$ of $x_0$ and $y_0$, respectively, and $k_1,...,k_m\in \N$ with
    \begin{align*}
        &(x_0,y_0)\in \op{supp}(\rho_{k_i}) \text{ for all } i=1,...,m
        \text{ and }
        \\
        &(U\cap V)\cap \op{supp}(\rho_{k})=\emptyset \text{ for } k\notin \{k_1,...,k_m\}.
    \end{align*}
    We may assume that $s(x_{k_1},y_{k_1})\geq s(x_{k_i},y_{k_i})$ for all $i=2,...,m$. Then, denoting  $W:=(U\cap U(x_{k_1},y_{k_1}))\times (V\cap V(x_{k_1},y_{k_1}))$, we have $s(x,y)\leq s_0(x_{k_1},y_{k_1})$ on $\op{Cut}_M^t\cap W$. Hence, on $[0,1]\times (\op{Cut}_M^t\cap W)$, we can write
    \begin{align*}
        \bar F(t,x,y)=F_{k_1}(ts(x,y),x,y),
    \end{align*}
    where $F_{k_1}$ denotes the map from the above theorem applied to the point $(x_{k_1},y_{k_1})$.
    Thus, by the preceding theorem, $\bar F$ is well-defined and continuous on $[0,1]\times (\op{Cut}_M^t\cap W)$ and satisfies (i)-(iii) on this set. For (iv), note that, by the theorem and the definition of $V_{y_{k_1}}$, we have $y,\bar F(t,x,y)\in B_{\ep(x_{k_1},y_{k_1})/4}(y_{k_1})$. Hence, $d_h(y,\bar F(t,x,y))\leq \ep(x_{k_1},y_{k_1})/2\leq \ep(x,y)$. Thus, all the properties hold on $[0,1]\times (\op{Cut}_M^t\cap W)$, which is an open neighbourhood of $(t,x_0,y_0)$ in $[0,1]\times \op{Cut}_M^t$. Since $(x_0,y_0)$ was arbitrary, this concludes the proof.
\end{proof}

\subsection{Proof of Theorem \ref{aa}}\label{ppoiouhgs}
To prove the local contractibility of ${\cal NU}$, let $(x_0,y_0)\in {\cal NU}$ be as in Definition \ref{defa}. We consider two cases. The simple case is when $(x_0,y_0)\in I^+$, and the difficult case is when $(x_0,y_0)\in \partial J^+$. Although the proof of the difficult cases also covers the simple case, we will prove both cases separately to emphasize that the first case is considerably easier. We start with the simple case.

\begin{proof}[Proof of Theorem \ref{aa} if $(x_0,y_0)\in I^+$]
Compare to Theorem 3.1 in \cite{Cannarsa/Cheng/Fathi}.
\begin{enumerate}[(a)]
\item 
    Let $(x_0,y_0)\in  {\cal NU}^t\subseteq \op{Cut}_M^t$, and let $W$ be any open neighbourhood of $(x_0,y_0)$. Let $U',V'$ be open neighbourhoods of $x_0$ and $y_0$, respectively, with $U'\times V'\subseteq W$. Pick $U\subseteq U'$, $V$, $s_0$ and $F$ as in Theorem \ref{ab} with $V_{y_0}:=V'$. Without loss of generality, we may assume that $U\times V$ is contractible to $(x_0,y_0)$, i.e.\ there exists a continuous functions $G:[0,1]\times U\times V\to U\times V$ such that $G(0,x,y)=(x,y)$ and $G(1,x,y)=(x_0,y_0)$ for all $(x,y)\in U\times V$.

    Now consider the homotopy 
    \[
    H:[0,s_0+1] \times ({\cal NU}\cap (U\times V))\to {\cal NU}\cap W
    \]
    defined by
    \begin{align*}
        H(s,x,y)=
        \begin{cases}
            (x,F(s,x,y)),& \text{ if } s\leq s_0,
            \\[10pt]
            (p_1\circ G(s-s_0,x,y),F(s_0,G(s-s_0,x,y))),& \text{ if } s> s_0,
        \end{cases}
    \end{align*}
    where $p_1:M\times M\to M$ is the projection onto the first factor. $H$ is well-defined because (1) $F$ maps $ [0,s_0]\times U\times V$ to $V'$ and $G$ maps to $U\times V$, (2) $(x,F(s,x,y))\in  {\cal NU}$  for all $(s,x,y)\in [0,s_0]\times ( {\cal NU}\cap (U\times V))$ thanks to part (iii) of Theorem \ref{ab}, and (3) $(p_1\circ G(s-s_0,x,y),F(s_0,G(s-s_0,x,y)))\in {\cal NU}$ for all $(s,x,y)\in [s_0,1+s_0]\times ( {\cal NU}\cap (U\times V))$ thanks to part (iv) of Theorem \ref{ab}. Thus, $H$ actually maps to $ {\cal NU}\cap W$ and is well-defined. Obviously, $H(0,x,y)=(x,y)$ and $H(s_0+1,x,y)=(x_0,F(s_0,x_0,y_0))$. $H$ is continuous since $G(0,x,y)=(x,y)$, hence $F(s_0,x,y)=F(s_0,G(0,x,y))$. This proves that $H$ satisfies all the required properties.
    \item 
    Let $x_0\in M$, and set ${\cal NU}(x_0):=\{y\in M\mid (x_0,y)\in {\cal NU}\}$. Let $y_0\in {\cal NU}(x_0)\cap I^+(x_0)$. Let $V'$ be any open neighbourhood of $y_0$, and pick $U,V$, $s_0$ and $F$ as in Theorem \ref{ab} with $V_{y_0}=V'$. We may assume that there is a contraction $G:[0,1]\times V\to V$ to $y_0$. We define 
    \[
    H:[0,s_0+1]\times ({\cal NU}(x_0)\cap V)\to {\cal NU}(x_0)\times V'
    \]
    by
    \begin{align*}
        H(y,s)=
        \begin{cases}
            F(s,x_0,y),& \text{ if } s\leq s_0,
            \\[10pt]
            F(s_0,x_0,G(s-s_0,y)),& \text{ if } s> s_0.
        \end{cases}
    \end{align*}
    Then one checks, as above, that $H$ satisfies all the required properties. This concludes the proof.
    \end{enumerate}
\end{proof}

Now we prepare for the proof when $(x_0,y_0)\in {\cal NU}\cap \partial J^+$.

\begin{lemma}
    Let $(x_0,y_0)\in \op{Cut}_M^n$, and let $W$ be any open neighbourhood of $(x_0,y_0)$. Then there exists a smaller open neighbourhood $W'\subseteq W$ of $(x_0,y_0)$ and a continuous homotopy $G:[0,1]\times W'\to W$ satisfying the following properties:
    \begin{enumerate}[(i)]
    \item $G(0,x,y)=(x,y)$ and $G(1,x,y)=(x_0,y_0)$ for all $(x,y)\in W'$.
    \item $G(t,x_0,y_0)=(x_0,y_0)$ for all $t\in [0,1]$.
    \item $G(t,x,y)\in I^+$ for all $t\in (0,1)$ and $(x,y)\in (W'\cap J^+)\backslash \{(x_0,y_0)\}$.
    \item $p_1\circ G(t,x_0,y)=x_0$ for all $t\in [0,1]$ and $y\in M$ with $(x_0,y)\in W'$.
    \end{enumerate}
\end{lemma}
\begin{proof}
    Let $(U,\phi)$ and $(V,\psi)$ be two charts around $x_0$ and $y_0$, respectively, with $U\times V\subseteq W$, which are centered at $x_0$ and $y_0$ (i.e.\ $\phi(x_0)=\psi(y_0)=0$), $\phi(U)=B_1(0)=\psi(V)$ and such that the cone
    \begin{align}
        C:=\{v\in \R^n\mid v_1\geq 0, v_1^2-\sum_{i=2}^n v_i^2\geq 0\}
         \label{eqah} 
    \end{align}
    is contained in $d_x\phi(\op{int}(\C_x)\cup\{0\})$ and $d_y\psi(\op{int}(\C_y)\cup\{0\})$ for all $x\in U$, $y\in V$. It follows that $\phi^{-1}((\phi(x_0)-C)\cap B_1(0))\backslash \{x_0\}\subseteq I^-(x_0)$, and similarly for $y_0$.
    It is easy to see that there exist two smaller open neighbourhoods $U'\subseteq U,V'\subseteq V$ of $x_0$ and $y_0$, respectively, such that, for all $x\in U'$ and $y\in V'$, the rays
    \begin{align*}
        \{\phi(x)-te_1\mid t\geq 0\},\{\psi(y)+te_1\mid t\geq 0\}
    \end{align*}
    intersect $(\phi(x_0)-C)\cap B_1(0)$ and $(\psi(y_0)+C)\cap B_1(0)$. Here $e_1:=(1,0,0,...)$. We set $W':=U'\times V'$. 
    
    We define the continuou map $I_1:U'\to B_1(0)$, assigning to each $x$ the point $\phi(x)-te_1$, where $t\in [0,\infty)$ is minimal with $\phi(x)-te_1\in (\phi(x_0)-C)\cap B_1(0)$. Analogously, we define a continuous map $I_2:V'\to B_1(0)$ for $y$. We then define the homotopy
    \begin{align*}
        &G:[0,1]\times W'\to W,\ G(t,x,y):=
        \\[10pt]
        &\begin{cases}
            \big(\phi^{-1}((1-2t)\phi(x)+2tI_1(x)),\psi^{-1}((1-2t)\psi(y)+2tI_2(y))\big), &\text{ if } t\leq 1/2,
            \\\\
            \big(\phi^{-1}((2-2t)I_1(x)),\psi^{-1}((2-2t)I_2(y))\big), & \text{ if } t\geq 1/2.
        \end{cases}
    \end{align*}
    Note that $G$ is well-defined, i.e.\ for instance that $(1-2t)\phi(x)+2tI_1(x)\in B_1(0)$ for all $t\leq 1/2$.
    This homotopy is obviously continuous and satisfies properties (i), (ii) and (iv). From the definition of $C$, it follows that $\phi^{-1}(\phi(x)-te_1)\in I^-(x)$ as long as $\phi(x)-te_1\in B_1(0)$ and $t>0$. A similar result holds for $y$. Moreover, as noted earlier, $\phi^{-1}((\phi(x_0)-C)\cap B_1(0))\backslash \{x_0\}\subseteq I^-(x_0)\subseteq I^-(y_0)$. A similar result holds for $y_0$. These observations imply that $G(t,x,y)\in I^+$ whenever $(x,y)\in W'\cap J^+$ and $(x,y)\neq (x_0,y_0)$, proving (iii).
\end{proof}

\begin{lemma}\label{lo}
    Let $(x_0,y_0)\in \op{Cut}_M^n$, and let $W$ be any open neighbourhood of $(x_0,y_0)$. Then there exists a smaller open neighbourhood $W'\subseteq W$ of $(x_0,y_0)$ and a continuous homotopy $G:[0,1]\times (W'\cap \op{Cut}_M)\to W\cap \op{Cut}_M$ satisfying the following properties:
    \begin{enumerate}[(i)]
    \item $G(0,x,y)=(x,y)$ and $G(1,x,y)=(x_0,y_0)$ for all $(x,y)\in W'\cap \op{Cut}_M$.
    \item $G(t,x_0,y_0)=(x_0,y_0)$ for all $t\in [0,1]$.
    \item $G(t,x,y)\in \op{Cut}_M^t$ for all $t\in (0,1)$ and $(x,y)\in (W'\cap \op{Cut}_M)\backslash \{(x_0,y_0)\}$. 
    \item $p_1\circ G(t,x_0,y)=x_0$ for all $t\in [0,1]$ and $y\in M$ with $(x_0,y)\in W'\cap \op{Cut}_M$. 
    \end{enumerate}
\end{lemma}
\begin{proof}
   By Corollary \ref{cor1}, there exists $\ep>0$ and an open neighbourhood $W''\subseteq W$ of $(x_0,y_0)$ such that, whenever $\gamma:[0,1]\to M$ is a maximizing causal geodesic with $(\gamma(0),\gamma(1))\in W$, then $\gamma$ can be exended to a geodesic on $[0,1+\ep]$. This extension, however, is not maximizing. 
    
    We define the map $f: W''\cap J^+\to W\cap \op{Cut}_M$ by
    \begin{align*}
        f(x,y):=(x,\exp_x(\alpha(x,v)v)),
    \end{align*}
    where $[0,1]\ni t\mapsto \exp_x(tv)$ is a maximizing geodesic connecting $x$ to $y$.
   
    This map is well-defined since (1) $\alpha(x,v)\leq 1+\ep$ for any such maximizing geodesic, hence $(x,\exp_x(\alpha(x,v)v))$ is defined and belongs to $\op{Cut}_M\cap W$ thanks to the definition of $\alpha$ and the first part of the proof,
    and (2) if $v,w$ both yield maximizing geodesics, then Theorem \ref{thme} guarantees $\alpha(x,v)=\alpha(x,w)=1$, so both possible definitions yield $f(x,y)=(x,y)$.
    Lemma \ref{lq} and Lemma \ref{lm} imply continuity of $f$. 
    Note that, obviously, $f(x,y)\in I^+(x)$ whenever $y\in I^+(x)$ and $f(x,y)=(x,y)$ whenever $(x,y)\in \op{Cut}_M$.

    Now let $G':[0,1]\times W'\to W''$ be the homotopy constructed in the previous lemma (applied with $W:=W''$). Then define
    \begin{align*}
        G:[0,1]\times (W'\cap \op{Cut}_M)\to W\cap \op{Cut}_M,\ 
         G(t,x,y):=f(G'(t,x,y)).
    \end{align*}
    This map is well-defined since $G'(t,x,y)\in W''\cap J^+$ for all $(t,x,y)\in [0,1]\times (W'\cap \op{Cut}_M)$ thanks to properties (i)-(iii) of the previous lemma.  Clearly, $G$ is continuous. Properties (i)-(iv) follow from the corresponding properties of $G'$, the definition and the above mentioned properties of $f$.
\end{proof}

\begin{proof}[Proof of Theorem \ref{aa} if $(x_0,y_0)\in \partial J^+$]
\begin{enumerate}[(a)]
\item 
    Let $(x_0,y_0)\in {\cal NU}\cap \partial J^+$, and let $W$ be any open neighbourhood of $(x_0,y_0)$. Choose open neighbourhoods $U',V'$ of $x_0$ and $y_0$, respectively, such that $U'\times V'\Subset  W$. Let $\ep:\op{Cut}_M^t\to (0,\infty)$ be defined by
    \[
    \ep(x,y):= \min\{d_{h\times h}((x,y),\partial W),d_{h\times h}((x,y),\partial J^+)\}
    \]
    $\text{ for } (x,y)\in (U'\times V')\cap \op{Cut}_M^t$, with an arbitrary continuous extension to $\op{Cut}_M^t$. Let $s$ and $\bar F$ be the maps provided by Corollary \ref{cor4}, and define
    \begin{align*}
         &K:[0,1]\times ((\op{Cut}_M^t\cup {\cal NU})\cap (U'\times V'))\to M,
         \\[10pt]
         &K(t,x,y):=
         \begin{cases}
             \bar F(t,x,y),& \text{ if } (x,y)\in \op{Cut}_M^t\cap (U'\times V'),
             \\\\
             y,& \text{ otherwise, i.e.\ } (x,y)\in {\cal NU}\cap (U'\times V')\cap \partial J^+.
         \end{cases}
    \end{align*}
    We claim that $K$ is continuous. Corollary \ref{cor4} guarantees continuity on the open set $[0,1]\times \op{Cut}_M^t\cap (U'\times V')$, so it remains to consider a sequence $(t_k,x_k,y_k)\to (t,x,y)$, where $(x_k,y_k)\in (\op{Cut}_M^t\cup {\cal NU})\cap (U'\times V')$ and $(x,y)\in {\cal NU}\cap (U'\times V')\cap \partial J^+\subseteq \op{Cut}_M^n$. By definition of $K$, we may assume that $(x_k,y_k)\in \op{Cut}_M^t\cap (U'\times V')$ for all $k$. By Corollary \ref{cor4}(iv), we have
    \[
    d_h(\bar F(t_k,x_k,y_k),y)
    \leq 
    d_h(\bar F(t_k,x_k,y_k),y_k)+d_h(y,y_k)
    \leq
    \ep(x_k,y_k)+d_h(y,y_k)
    \] 
    and the latter expression tends to $0$ as $k\to \infty$ by definition of $\ep$. Thus, $K$ is continuous.
    
    Moreover, Corollary \ref{cor4}(iii),(iv) and the definition of $\ep$ imply $(x,K(t,x,y))\in {\cal NU}\cap W$ for all $t\in [0,1]$ and $(x,y)\in {\cal NU}\cap (U'\times V')$, and $(x,K(1,x,y))\in {\cal NU}\cap W$ for all $(x,y)\in (\op{Cut}_M^t\cup {\cal NU})\cap (U'\times V')$.
    
    By Lemma \ref{lo} (applied with $W:=U'\times V'$), we can find two open neighbourhoods $U\subseteq U',V\subseteq V'$ of $x_0$ and $y_0$, and a homotopy $G:[0,1]\times ((U\times V)\cap \op{Cut}_M)\to U'\times V'$ satisfying the properties listed in the lemma. Finally, we define
    \begin{align*}
        &H:[0,2]\times ({\cal NU}\cap (U\times V))\to {\cal NU}\cap W, 
        \\[10pt]
        &H(t,x,y):=
        \begin{cases}
            (x,K(t,x,y)),\ &\text{ if } t\leq 1,
            \\\\
            (p_1\circ G(t-1,x,y),K(1,G(t-1,x,y))), &\text{ if } t> 1.
        \end{cases}
    \end{align*}
    Note that $H$ is well-defined since (1) ${\cal NU}\cap (U\times V)\subseteq ((\op{Cut}_M^t\cup{\cal NU})\cap (U'\times V')$ and ${\cal NU}\cap (U\times V) \subseteq (U\times V)\cap \op{Cut}_M$, (2) $(x,K(t,x,y))\in {\cal NU}\cap W$ for all $(t,x,y)\in [0,1]\times ({\cal NU}\cap (U\times V))$ as observed above, and (3) $G(t-1,x,y)\in (\op{Cut}_M^t\cup{\cal NU})\cap (U'\times V')$ for all $(t,x,y)\in (1,2]\times ({\cal NU}\cap (U\times V))$ by Lemma \ref{lo}, hence $H(t,x,y)\in {\cal NU}\cap W$ as observed above. Continuity of $H$ follows from continuity of $K$ and $G$. Clearly, $H(0,x,y)=(x,y)$ and $H(2,x,y)=(x_0,y_0)$. 
    \item 
    Fix $x_0\in M$, and let $y_0\in M$ with $(x_0,y_0)\in {\cal NU}\cap \partial J^+$. Let $V'\Subset V''$ be any two open neighbourhoods of $y_0$. Pick two arbitrary open neighbourhoods $U'\Subset U''$ of $x_0$ and set $W:=U''\times V''$. We redo (with the same notation) the proof of part (a), and consider the homotopy 
    \begin{align*}
        &\tilde H:[0,2]\times \{y\in V\mid (x_0,y)\in {\cal NU}\}\to \{y\in V''\mid (x_0,y)\in {\cal NU}\},
        \\[10pt]
        &\tilde H(t,y):=p_2\circ H(t,x_0,y),
    \end{align*}
    where $p_2:M\times M\to M$ denotes the projection onto the second factor. Note that $H(t,x_0,y)=(x_0,\tilde H(t,y))$ thanks to Lemma \ref{lo}(iv). Therefore, all the properties in the definition of local contractibility follow from the corresponding properties of $H$.
\end{enumerate}
\end{proof}

\subsection{Proof of Theorem \ref{b}(b)}\label{ppoiouhgs1}

We will prove Theorem \ref{b} in two steps. In one step, we show that the sets $\op{Cut}_M,\op{Cut}_M^t,\op{Cut}_M(x)$ and $\op{Cut}_M^t(x)$ are strong deformation retracts of the sets $J^+\backslash {\cal A}$, $I^+\backslash {\cal A}$, $J^+(x)\backslash {\cal A}(x)$ and $I^+(x)\backslash {\cal A}(x)$, respectively (Propositions \ref{pro2} and \ref{pro4}). For the versions involving a fixed point $x$, this is particularly intuitive; the point $y$ is moved along the future inextendible geodesic through $x$ and $y$ that is maximizing on the segment between them, until this geodesic intersects $\op{Cut}_M(x)$ (or $\op{Cut}_M^t(x)$). Theorem \ref{thme} and the lemma below ensure that this construction is well-defined, while continuity follows from continuity of the map $\alpha$. A similar result in the compact Riemannian setting can be found in \cite{Klingenberg}, Theorem 2.1.8. The same idea extends to the cases $\op{Cut}_M$ and $\op{Cut}_M^t$, although the proof requires additional refinements. Let us note that neither of these proofs rely on Theorem \ref{main2} or on the results from Sections 3 and 4.

In the second step, we prove that the inclusions from 
\[
{\cal NU},{\cal NU}^t,\{y\in J^+(x)\mid (x,y)\in {\cal NU}\},\{y\in I^+(x)\mid (x,y)\in {\cal NU}\}
\]
to $\op{Cut}_M,\op{Cut}_M^t,\op{Cut}_M(x),\op{Cut}_M^t(x)$ are homotopy equivalences (Proposition \ref{pro33} and \ref{pro5}). This is the point where Theorem \ref{main2} becomes essential.

\begin{lemma}\label{laa}
    Let $x\in M$, $y\in J^+(x)$ and $[0,1]\ni t\mapsto \exp_x(tv)$ be a maximizing geodesic connecting $x$ to $y$. Then the following statements are equivalent:
    \begin{enumerate}[(i)]
        \item $y\notin {\cal A}(x)$
        \item $\alpha(x,v)<\infty$ and $\exp_x(\alpha(x,v)v)$ exists.
    \end{enumerate}  
    Moreover, if either of these conditions holds, then $\exp_x(\alpha(x,v)v)\in \op{Cut}_M(x)$.
\end{lemma}
\begin{proof}
    The implication $(i) \Rightarrow (ii)$, as well as the final statement, follow directly from the definitions. For the converse implication, suppose that $y\in {\cal A}(x)$. Then there is a ray ${\tilde \gamma}:[0,a)\to M$, $a\in (0,\infty]$, with ${\tilde \gamma}(0)=x$ and ${\tilde \gamma}(t_1)=y$ for some $t_1\in (0,a)$. By rescaling, we may assume that $t_1=1$. Write ${\tilde \gamma}=[t\mapsto \exp_x(t\tilde v)]$ for some $\tilde v\in \C_x$. If $v\neq \tilde v$, there would exist two distinct maximizing geodesics connecting $x$ to $y$. By Theorem \ref{thme}, both geodesics must stop being maximizing at $y$, contradicting the fact that ${\tilde \gamma}$ is a ray, i.e.\ maximizing throughout its domain. Thus, $v=\tilde v$. Since ${\tilde \gamma}$ is maximizing and inextendible, it follows that, if $\alpha(x,v)$ is finite, then $\exp_x(\alpha(x,v)v)$ cannot exist. 
\end{proof}

\begin{proposition}\label{pro2}
 For any $x\in M$, $\op{Cut}_M(x)$ is a strong deformation retract of $J^+(x)\backslash {\cal A}(x)$. Moreover, there exists a strong deformation retraction that restricts to a strong deformation retraction from $I^+(x)\backslash {\cal A}(x)$ onto $\op{Cut}^t_M(x)$.
\end{proposition}
\begin{proof} 
    First note that, by the above lemma, $\op{Cut}_M(x)$ and $\op{Cut}_M^t(x)$ are indeed contained in $J^+(x)\backslash {\cal A}(x)$ and $I^+(x)\backslash {\cal A}(x)$, respectively.

    We need to prove existence of a continuous homotopy
    \begin{align*}
        H:[0,1]\times J^+(x)\backslash {\cal A}(x)\to J^+(x)\backslash {\cal A}(x)
    \end{align*} 
    satisfying the properties
    \begin{itemize}
        \item $H(0,y)=y$ for all $y\in J^+(x)\backslash {\cal A}(x)$.
        \item $H(1,y)\in \op{Cut}_M(x)$ for all $y\in J^+(x)\backslash {\cal A}(x)$.
        \item $H(t,y)=y$ for all $y\in \op{Cut}_M(x)$ and $t\in [0,1]$.
        \item $H(t,y)\in I^+(x)$ for all $y\in I^+(x)\backslash {\cal A}(x)$ and $t\in [0,1]$.
    \end{itemize}
    Indeed, the first three properties prove the first statement, and the last one shows the second.
    
     We define 
    \begin{align*}
        H:[0,1]\times J^+(x)\backslash {\cal A}(x)\to J^+(x),\ H(t,y):=\exp_x(((1-t)+t\alpha(x,v))v),
    \end{align*} where $[0,1]\ni t\mapsto \exp_x(tv)$ is any maximizing geodesic connecting $x$ to $y$. 
    
   $H$ is well-defined since (1) $\exp_x(\alpha(x,v)v)$ is defined by Lemma \ref{laa} and $\alpha(x,v)\geq 1$, hence $\exp_x(((1-t)+t\alpha(x,v))v)$ is defined for all $t\in [0,1]$, and (2) if $v,w$ both yield maximizing geodesics, then $\alpha(x,v)=\alpha(x,w)=1$ by Theorem \ref{thme}, hence both possible definitions give $H(t,x,y)=(x,y)$. 

    To prove continuity, suppose $(t_k,y_k),(t,y)\in J^+(x)\backslash {\cal A}(x)$ and let $(t_k,y_k)\to (t,y)$. Denote by $\gamma_k:[0,1]\to M$ maximizing geodesics connecting $x$ to $y_k$, and let $v_k\in \C_{x}$ with $\gamma_k(s)=\exp_{x}(sv_k)$. By Lemma \ref{lm}, after passing to a subsequence, $v_k\to v\in \C_x$, and the curve $[0,1]\ni s\mapsto \exp_x(sv)$ is a maximizing geodesic connecting $x$ to $y$. By definition of ${\cal A}(x)$, $\exp_x(\alpha(x,v)v)$ exists (see Lemma \ref{laa}), implying that $\alpha$ is continuous at $(x,v)$. Thus
    \begin{align*}
        H(t,y)=\exp_x(((1-t)+t\alpha(x,v))v)
        &= \lim_{k\to \infty} \exp_{x}(((1-t_k)+t_k\alpha(x,v_k))v_k)
        \\[10pt]
        &=\lim_{k\to \infty} H(t_k,y_k).
    \end{align*}
    This proves the continuity. 
    
    Clearly, $H(0,y)=y$. By Lemma \ref{laa}, $H(1,y)\in \op{Cut}_M(x)$ for all $y\in J^+(x)\backslash {\cal A}(x)$. Furthermore, $H(t,y)=y$ for all $y\in \op{Cut}_M(x)$ and $t\in [0,1]$, and $H(t,y)\in I^+(x)$ whenever $y\in I^+(x)\backslash {\cal A}(x)$ and $t\in [0,1]$. To conclude the proof, it suffices to prove that $H(t,y)\in J^+(x)\backslash {\cal A}(x)$ for $(t,y)\in [0,1]\times J^+(x)\backslash {\cal A}(x)$: The geodesic $[0,1]\ni s\mapsto \exp_x(s((1-t)+t\alpha(x,v))v)$ is maximizing and connects $x$ to $H(t,y)$. Its maximal future extension has a cut point at $s=\alpha(x,v)((1-t)+t\alpha(x,v))^{-1}$, so Lemma \ref{laa} implies $H(t,y)\notin {\cal A}(x)$.
\end{proof}

\begin{definition}\rm
    Let $X:M\to TM$ be a smooth timelike vector field, whose existence is guaranteed by the time orientability of $M$. We denote by 
    \[
    \varphi:(0,\infty)\times M\supseteq {\cal D}\to M
    \]
    its smooth local flow.
\end{definition}

\begin{lemma}\label{lr}
    Let $x\in M$. The future Aubry set ${\cal A}(x)$ is closed.
\end{lemma}
\begin{proof}
    Since $J^+(x)$ is closed, it suffices to prove that ${\cal A}(x)$ is closed relative to $J^+(x)$, or equivalently, that $J^+(x)\backslash {\cal A}(x)$ is relatively open in $J^+(x)$.

    Let $y\in J^+(x)\backslash {\cal A}(x)$. Suppose, for contradiction, that there exists a sequence $y_k\in {\cal A}(x)$ converging to $y$.  Let $[0,1]\ni t\mapsto \exp_{x}(tv_k)$ be maximizing geodesics connecting $x$ to $y_k$. By Lemma \ref{lm}, after passing to a subsequence, we have $v_k\to v\in \C_x$, and the curve $[0,1]\ni t\mapsto \exp_x(tv)$ is a maximizing geodesic connecting $x$ to $y$. By definition of the future Aubry set, $\alpha(x,v)<\infty$ and $\exp_x(\alpha(x,v)v)$ exists. Hence, $\alpha$ is continuous at $(x,v)$, implying that $\alpha(x,v_k)<\infty$ and that $\exp_x(\alpha(x,v_k)v_k)$ exists for sufficiently large $k$. Thus, by Lemma \ref{laa}, $y_k\notin {\cal A}(x)$ for these $k$, contradicting the assumption. 
\end{proof}

\begin{proposition}\label{pro33}
    For any $x\in M$, the inclusion 
    \[
    \{y\in J^+(x)\mid (x,y)\in {\cal NU}\}\hookrightarrow \op{Cut}_M(x)
    \]
    is a homotopy equivalence, which restricts to a homotopy equivalence 
    \[
    \{y\in I^+(x)\mid (x,y)\in {\cal NU}\}\hookrightarrow \op{Cut}_M^t(x).
    \]
\end{proposition}
\begin{proof}
    Fix $x\in M$. Thanks to the preceding lemma, we can choose a continuous function $T:M\backslash {\cal A}(x)\to (0,\infty)$ such that $\varphi(t,y)\in M\backslash{\cal A}(x)$ for all $y\in M\backslash {\cal A}(x)$ and $t\in [0,T(y)]$. 
    We define the continuous map
    \begin{align*}
        G:[0,1]\times \op{Cut}_M(x)\to J^+(x)\backslash {\cal A}(x),
        \ (t,y)\mapsto \varphi(tT(y),y).
    \end{align*}
    Let $\ep:\op{Cut}_M^t\to (0,\infty),\ \ep(x',y):=d_{h\times h}((x',y),\partial J^+)$, and let $\bar F$ and $s:\op{Cut}_M^t\to (0,\infty)$ denote the maps from Corollary \ref{cor4} associated with $\ep$. Denote by $\tilde H:[0,1]\times (J^+(x)\backslash {\cal A}(x))\to J^+(x)\backslash {\cal A}(x)$ the strong deformation retraction from Proposition \ref{pro2}.
    We define
    \begin{align*}
        &H:[0,1]\times \op{Cut}_M(x)
        \to \op{Cut}_M(x),
        \\[10pt]
        &H(t,y):=
        \begin{cases}
            y,\ &\text{if } t=0,
            \\\\
            \bar F(t,x,\tilde H(1,G(t,y))),\ &\text{ if } t>0.
        \end{cases}
    \end{align*}
    Since the vector field $X$ is timelike, $G(t,y)\in I^+(x)\backslash {\cal A}(x)$ for $t>0$, so Proposition \ref{pro2} implies $\tilde H(1,G(t,y))\in \op{Cut}_M^t(x)$. Hence, $H$ is well-defined and 
    \begin{align}
        (x,H(t,y))\in {\cal NU}^t\subseteq {\cal NU} \text{ for all } t>0. \label{eqaj}
    \end{align}
    
    We claim that $H$ is continuous. As a composition of continuous functions, continuity holds on $(0,1]\times \op{Cut}_M(x)$. For $(t,y)\in [0,1]\times \op{Cut}_M^t(x)$, we have $\tilde H(1,G(t,y))\in \op{Cut}_M^t(x)$ thanks to Proposition \ref{pro2}. Since $\bar F(\cdot,x,\cdot)$ is well-defined and continuous on $[0,1]\times \op{Cut}_M^t(x)$ and $\bar F(0,x,\tilde H(1,G(0,y)))=y$ for $y\in \op{Cut}_M^t(x)$, it is immediate that also $H$ is continuous on the open set $[0,1]\times \op{Cut}_M^t(x)$. Finally, if $t=0$ and $y\in \op{Cut}_M^n(x)$, and $(0,1]\times \op{Cut}_M(x)\ni (t_k,y_k)$ converges to $(0,y)$, then $y_k':=\tilde H(1,G(t_k,y_k))$ converges to $\tilde H(1,y)=y$. Hence, 
    \begin{align*}
    d_{h}(H(t_k,y_k),y)
    &\leq 
    d_h(H(t_k,y_k),y_k')+d_h(y_k',y)
    \\
    &=
    d_{h\times h}\big((x,\bar F(t_k,x,y_k')\big),(x,y_k'))+d_h(y_k',y)
    \\
    &\leq 
    \ep(x,y_k')+d_h(y_k',y)
    \xrightarrow{k\to \infty} 0.
    \end{align*}
    In the last step, we used $(x,y_k')\to (x,y)\in \partial J^+$. This proves the continuity.
    
    We claim that $H(1,\cdot)$ is a homotopy inverse to the inclusion 
    \[
    \iota:\{y\in J^+(x)\mid (x,y)\in {\cal NU}\}\hookrightarrow \op{Cut}_M(x),
    \]
    and that the restriction $H(1,\cdot)_{|\op{Cut}_M^t(x)}$ is a homotopy inverse to the inclusion 
    \[
    \iota^t:\{y\in I^+(x)\mid (x,y)\in {\cal NU}\}\hookrightarrow \op{Cut}^t_M(x).
    \]

    Indeed, by \eqref{eqaj}, $H(1,\cdot)$ maps to $\{y\in J^+(x)\mid (x,y)\in {\cal NU}\}$. The map $\iota\circ H(1,\cdot)=H(1,\cdot)$ is homotopic to $\op{Id}_{\op{Cut}_M(x)}=H(0,\cdot)$ in $\op{Cut}_M(x)$ via the homotopy $H$. Conversely, using \eqref{eqaj}, we have $(x,H(t,y))\in {\cal NU}$ whenever $(x,y)\in {\cal NU}$, so the composition $H(1,\cdot)\circ \iota$ is homotopic to $\op{Id}_{\{y\mid (x,y)\in {\cal NU}\}}$ in $\{y\mid (x,y)\in {\cal NU}\}$ via the homotopy $(t,y)\mapsto H(t,\iota(y))$. Using that $(x,H(t,y))\in {\cal NU}^t$ whenever $t>0$ (see \eqref{eqaj}), the same arguments show that also $H(1,\cdot)_{|\op{Cut}_M^t(x)}$ is a homotopy inverse to $\iota^t$.
\end{proof}

\begin{proof}[Proof of Theorem \ref{b}(b)]
   The inclusion of a strong deformation retract into its ambient space is a homotopy equivalence, and the composition of two homotopy equivalences is again a homotopy equivalence. Hence, the theorem follows from Propositions \ref{pro2} and \ref{pro33}.
\end{proof}

\subsection{Proof of Theorem \ref{b}(a)}\label{ppoiouhgs2}

The strategy in proving (a) is essentially the same as for part (b). However, the analogue of Proposition \ref{pro2}, namely, Proposition \ref{pro4}, requires additional care. As explained in the previous subsection, in the proof of Proposition \ref{pro2}, the key idea was to push a point $y$ along the maximizing geodesic connecting $x$ to $y$ until it intersects $\op{Cut}_M(x)$. This worked since $y\in {\cal A}(x)$. In contrast, if $(x,y)\in {\cal A}$, then the maximal future extension of the maximizing geodesic connecting $x$ to $y$ might be maximizing. More precisely:

\begin{lemma}\label{ls}
    Let $(x,y)\in J^+$, and let $\gamma:[0,1]\to M$ be a maximizing geodesic connecting $x$ to $y$. Then the following are equivalent:
    \begin{enumerate}[(i)]
        \item $(x,y)\notin {\cal A}$
        \item The maximal geodesic extension of $\gamma$ is not maximizing.
    \end{enumerate}
\end{lemma}
\begin{proof}
    The implication $(i)\Rightarrow (ii)$ follows immediately from the definition. For the converse, suppose $(x,y)\in {\cal A}$, and let $\tilde \gamma:I\to M$ be a line through $x$ and $y$. Without loss of generality, $\tilde \gamma(0)=x$ and $\tilde \gamma(1)=y$. Write $\tilde \gamma=\exp_x(t\tilde v)$ for some $\tilde v\in \C_x$. As in the proof of Lemma \ref{laa}, one shows $v=\tilde v$, where $\gamma(t)=\exp_x(tv)$. Therefore, the maximal geodesic extension of $\gamma$ must be $\tilde \gamma$, which is maximizing. 
\end{proof}

\begin{corollary}\label{cor5}
    It holds $\op{Cut}_M\subseteq J^+\backslash {\cal A}$.
\end{corollary}

\begin{lemma}\label{lu}
    There exist two continuous functions 
    \[
    \varphi^+:J^+\backslash {\cal A}\to (1,\infty) 
    \text{ and }
    \varphi^-:J^+\backslash {\cal A}\to (-\infty,0]
    \]
    such that the following holds:

    Whenever $(x,y)\in J^+\backslash {\cal A}$, and $\gamma:[0,1]\to M$ is a maximizing geodesic connecting $x$ to $y$, then the maximal geodesic extension of $\gamma$ is defined but not maximizing on the interval $[\varphi^-(x,y),\varphi^+(x,y)]$.
\end{lemma}
\begin{proof}
    Let $(x,y)\in J^+\backslash {\cal A}$ be arbitrary, and $\gamma:[0,1]\to M$ be a maximizing geodesic connecting $x$ to $y$. We claim that there exist $a\leq 0$ and $b> 1$ and two open neighbourhoods $U$ and $V$ of $x$ and $y$, respectively, such that, whenever $(x',y')\in J^+\cap (U\times V)$ and $[0,1]\to M$ is maximizing geodesic connecting $x'$ to $y'$, then its maximal geodesic extension is defined but not maximizing on $[a,b]$. 
    
    Indeed, if $(x,y)\in \op{Cut}_M$, then the claim follows from Corollary \ref{cor1}.
    
    Otherwise, by Theorem \ref{thme}, $\gamma$ is the unique maximizing geodesic connecting $x$ to $y$. Since $(x,y)\in {\cal A}$, we can pick $a\leq 0$ and $b> 1$ such that the maximal geodesic extension of $\gamma$ (say $\tilde \gamma$) is defined but not maximizing on $[a,b]$. From Lemma \ref{lm}, it follows that, whenever $J^+\ni (x_k,y_k) \to (x,y)$, and $\gamma_k:[0,1]\to M$ is a maximizing geodesic connecting $x_k$ to $y_k$, then its maximal geodesic extension $\tilde \gamma_k$ is defined on $[a,b]$ and converges, in the $C^1([a,b])$-topology, to ${\tilde \gamma}_{|[a,b]}$. In particular, by continuity of the Lorentzian distance, if ${\tilde \gamma_k}$ would be maximizing on $[a,b]$ for infinitely many $k$, so would $\tilde \gamma$. This is a contradiction.

    We do this construction for every $(x,y)\in J^+\backslash {\cal A}$. Obviously, 
    \begin{align*}
        \bigcup_{(x,y)\in J^+\backslash {\cal A}} U_{(x,y)}\times V_{(x,y)} \supseteq J^+\backslash {\cal A}.
    \end{align*}
    Let $\{W'_i\}_{i\in I}$ be a locally finite refinement of this open cover such that each $W'_i$ is compactly embedded in some $U_{(x_i,y_i)}\times V_{(x_i,y_i)}$ (\cite{Lee}, Theorem 1.15). Further, let $\{W_j\}_{j\in J}$ be a locally finite refinement of the open cover $\{W_i'\}$ such that each $W_j$ is compactly embedded in some $W_{i_j}'$ (\cite{Lee}, Theorem 1.15). For each $j\in J$, choose a smooth bump function $\rho_j:M\to [0,1]$ for $\overbar W_j$ in $W_{i_j}'$, i.e.\ $\rho_{|\overbar W_j}\equiv 1$ and $\op{supp}(\rho)\subseteq W_{i_j}'$. We define 
    \begin{align*}
        \varphi^+(x,y):=\max\{\rho_j(x,y)b_{i_j}\mid j\in J\}  \text{ and } 
        \varphi^-(x,y):=\min\{\rho_j(x,y)a_{i_j}\mid j\in J\},
    \end{align*}
    where $b_{i}:=b_{(x_{i},y_{i})}$ and $a_{i}:=a_{(x_{i},y_{i})}$, $i\in I$, are the interval endpoints associated to $(x_i,y_i)$.

    Since the covers $\{W_i'\}_{i\in I}$ and $\{W_j\}_{j\in J}$ are locally finite, these maxima and minima are locally finite maxima and minima of smooth functions, and thus $\varphi^\pm$ is real valued and continuous.
     Obviously, 
     \[
     \max\{b_i\mid (x,y)\in W_i'\}\geq \varphi^+(x,y)\geq b_{i_j}>1
     \]
     and 
     \[
     \min\{a_i\mid (x,y)\in W_i'\}\leq \varphi^-(x,y)\leq a_{i_j}\leq 0
     \]
     on $W_j\subseteq U_{(x_{i_j},y_{i_j})}\times V_{(x_{i_j},y_{i_j})}$. Since the $W_j$ cover $J^+\backslash {\cal A}$, it follows that, whenever $(x,y)\in J^+\backslash {\cal A}$ and $\gamma:[0,1]\to M$ is a maximizing geodesic connecting $x$ to $y$, then its maximal geodesic extension is defined but not maximizing on $[\varphi^-(x,y),\varphi^+(x,y)]$. 
\end{proof}

\begin{definition}\rm
    We define the function $\beta:J^+\backslash {\cal A}\to [0,1)$ by
    \begin{align*}
       \beta(x,y):=\sup\{t\geq 0\mid \gamma \text{ is maximizing on } [t\varphi^-(x,y)),(1-t)+t\varphi^+(x,y)]\},
    \end{align*}
    where $\gamma:[0,1]\to M$ is a maximizing geodesic connecting $x$ to $y$ (and also denotes its maximal extension). We will show below that $\beta$ is well-defined.
\end{definition}

\begin{lemma}
    Let $\gamma:I\to M$ be a causal geodesic defined on an open interval $I$. Suppose that $a,b\in I$ are such that $\gamma$ is maximizing $[a,b]$ but not on any interval $[a-\ep,b+\ep],\ \ep>0$. Then $(\gamma(a),\gamma(b))\in \op{Cut}_M$.
\end{lemma}
\begin{proof}
    Without loss of generality, we may assume that $a=0$ and $b=1$. We must show that $\gamma(1)$ is the cut point of $\gamma(0)$ along $\gamma$. Suppose, by contradiction, that $\gamma(1)$ is not the cut point of $\gamma(0)$ along $\gamma$. Then $\gamma$ is maximizing on $[0,1+\ep]$ for small $\ep>0$. In particular, $\gamma(1)$ is not conjugate to $\gamma(0)$ along $\gamma$ (hence, by Remark \ref{hdsdoa}, $\exp:TM\supseteq \op{dom}(\exp)\to M^2$ is a local diffeomorphism near $(\gamma(0),\dot \gamma(0))$) and $\gamma$ is the unique (up to reparametriaztion) maximizing geodesic connecting these two points. 
    
    Now let $\ep_k\to 0$ be a sequence of positive numbers. By assumption, $\gamma$ is not maximizing on $[-\ep_k,1+\ep_k]$. For each $k$, pick a maximizing geodesic $\gamma_k:[0,1]\to M$ connecting $\gamma(-\ep_k)$ to $\gamma(1+\ep_k)$. By Lemma \ref{lm},  $\gamma_k$ converges in the $C^1$-topology to $\gamma$. Note that $(\gamma_k(0),\dot \gamma_k(0)),(\gamma(-\ep_k),(1+2\ep_k)\dot \gamma(-\ep_k))$ converge to $(\gamma(0),\dot \gamma(0))$. Since
    \begin{align*}
        \exp(\gamma_k(0),\dot \gamma_k(0))=(\gamma(-\ep_k),\gamma(1+\ep_k))=\exp(\gamma(-\ep_k),(1+2\ep_k)\dot \gamma(-\ep_k)),
    \end{align*}
    and since $\exp$ is a local diffeomorphism near $(\gamma(0),\dot \gamma(0))$, it follows that $\dot \gamma_k(0)=(1+2\ep_k)\dot \gamma(0)$ for large $k$, so $\gamma_k$ is a reparametrization of $\gamma_{|[-\ep_k,1+\ep_k]}$. Thus, since $\gamma_k$ is maximizing, also $\gamma$ must be maximizing on $[-\ep_k,1+\ep_k]$, contradicting the assumption. Hence, $\gamma(1)$ must be the cut point of $\gamma(0)$ along $\gamma$.
\end{proof}

The proof of the following lemma is similar to the prove of Proposition 9.33 in \cite{Ehrlich}.

\begin{lemma}
    The map $\beta$ is well-defined and continuous. Moreover, if $\gamma:[0,1]\to M$ is a maximizing geodesic connecting $x$ to $y$, then 
    \begin{align}
        \big(\gamma(\beta(x,y)\varphi^-(x,y)),\gamma((1-\beta(x,y))+\beta(x,y)\varphi^+(x,y))\big)\in \op{Cut}_M. \label{eqvv}
    \end{align}
\end{lemma}
\begin{proof}
    If there exist two distinct maximizing geodesics $[0,1]\to M$ connecting $x$ to $y$, then, by Theorem \ref{thme}, $y$ is the cut point of $x$ along both of them. Hence, since $\varphi^+(x,y)>1$, both possible definitions give $\beta(x,y)=0$. Also, in the defintion of $\beta$, by Lemma \ref{lu}, the condition is violated for $t=1$. Therefore, by continuity of $d$, $\beta$ maps to $[0,1)$. Hence, $\beta$ is well-defined.
    
    To prove continuity, let $(x,y)\in J^+\backslash {\cal A}$, and suppose that $(x_k,y_k)\in J^+\backslash {\cal A}$ is a sequence converging to $(x,y)$. For each $k$, let $[0,1]\ni t \mapsto \exp_{x_k}(tv_k)$ be a maximizing geodesic connecting $x_k$ to $y_k$. Along a subsequence, $(x_k,v_k)\to (x,v)\in \C$ and  $[0,1]\ni t\mapsto \exp_x(tv)$ is a maximizing geodesic connecting $x$ to $y$. 
    
    We must prove that $\beta(x_k,y_k)\to \beta(x,y)$. Denote 
    \[
    \beta:=\beta(x,y),\ \varphi^\pm:=\varphi^\pm(x,y) \text{ and } \varphi_k^\pm:=\varphi^\pm(x_k,y_k).
    \]
    
    \noindent \textbf{Upper semicontinuity:}
    Assume for contradiction that $\limsup_{k\to \infty} \beta(x_k,y_k)\geq \beta+2\ep$ for some $\ep>0$. Without loss of generality, we have $\lim_{k\to \infty} \beta(x_k,y_k)\geq \beta +2\ep$. Choosing $\ep<1-\beta$, we may assume that $\exp_x$ is defined at $-(\beta+\ep)\varphi^-v$ and $((1-(\beta+\ep))+(\beta+\ep)\varphi^+)v$. Then
    \begin{align*}
        &d\big(\exp_x(-(\beta+\ep)\varphi^-v),\exp_x([(1-(\beta+\ep))+(\beta+\ep)\varphi^+]v)\big)
        \\[10pt]
        =
        &\lim_{k\to \infty} 
        d\big(\exp_{x_k}(-(\beta+\ep)\varphi^-_kv_k),\exp_{x_k}([(1-(\beta+\ep))+(\beta+\ep)\varphi_k^+]v_k\big)
        \\[10pt]
        =&\lim_{k\to \infty}\Big[(1-(\beta+\ep))+ (\beta+\ep)(\varphi^+_k-\varphi^-_k)\Big]d(x_k,y_k)
        \\[10pt]
        = &\Big[(1-(\beta+\ep))+(\beta+\ep)(\varphi^+-\varphi^-)\Big]d(x,y)
        \\[10pt]
        = &\ell_g\Big([-(\beta+\ep)\varphi^-,(1-(\beta+\ep))+(\beta+\ep)\varphi^+]\ni t \mapsto \exp_x(tv)\Big)
    \end{align*}
    This contradicts the definition of $\beta$, proving the upper semicontinuity.
    \\
    
    \noindent \textbf{Lower semicontinuity:} Assume, for contradiction, that $\lim_{k\to \infty} \beta(x_k,y_k)\leq \beta-2\ep$ for some $\ep>0$.
    Then, for large $k$, $\exp_{x_k}(tv_k)$ is defined but not maximizing on the interval $[(\beta-\ep)\varphi_k^-,(1-(\beta-\ep))+(\beta-\ep)\varphi_k^+]$. We define the curves $\gamma,\gamma_k:[0,1]\to M$ by
    \begin{align*}
        &\gamma(t):=\exp_x\bigg(\Big[(1-t)\big((\beta-\ep)\varphi^-\big)+t\big( (1-(\beta-\ep))+(\beta-\ep)\varphi^+\big)\Big]v\bigg) \text{ and }
        \\[10pt]
        &\gamma_k(t):=\exp_x\bigg(\Big[(1-t)\big((\beta-\ep)\varphi_k^-\big)+t\big( (1-(\beta-\ep))+(\beta-\ep)\varphi_k^+\big)\Big]v_k\bigg)
    \end{align*}
    By Theorem \ref{thme}, $\gamma$ is the unique maximizing geodesic connecting $\gamma(0)$ to $\gamma(1)$, and $\gamma(1)$ is not conjugate to $\gamma(0)$ along $\gamma$. In particular, $\exp:TM\supseteq \op{dom}(\exp)\to M^2$ is a local diffeomorphism near $(\gamma(0),\dot \gamma(0))$. By assumption, $\gamma_k$ is not maximizing for large $k$. Call $\tilde \gamma_k:[0,1]\to M$ a maximizing geodesic connecting $\gamma_k(0)$ to $\gamma_k(1)$. Lemma \ref{lm} implies that, up to a subsequence, $\tilde \gamma_k$ converges in the $C^1$-topology to a maximizing geodesic $\tilde \gamma:[0,1]\to M$ connecting $\gamma(0)$ to $\gamma(1)$, which must be $\gamma$. Hence, both $(\gamma_k(0),\dot \gamma_k(0))$ and $(\tilde \gamma_k(0),\dot {\tilde \gamma}_k(0))$ converge to  $(\gamma(0),\dot \gamma(0))$. But $\exp$ is a local diffeomorphism around $(\gamma(0),\dot \gamma(0))$ and
    \begin{align*}
        \exp(\gamma_k(0),\dot \gamma_k(0))=\exp(\tilde \gamma_k(0),\dot {\tilde \gamma}_k(0))=(\gamma_k(0),\gamma_k(1)),
    \end{align*}
    implying that $\dot {\tilde \gamma}_k(0)=\dot \gamma_k(0)$ for large $k$, thus $\gamma_k=\tilde \gamma_k$ is maximizing. This is a contradiction and proves continuity.

    Finally, \eqref{eqvv} follows from the previous lemma and the fact that $\gamma$ is maximizing on $[\beta(x,y)\varphi^-(x,y),(1-\beta(x,y))+\beta(x,y)\varphi^+(x,y)]$, thanks to the continuity of $d$.
\end{proof}

\begin{proposition}\label{pro4}
    $\op{Cut}_M$ is a strong deformation retract of $J^+\backslash {\cal A}$. Moreover, there exists a strong deformation retraction that restricts to a strong deformation retraction from $I^+\backslash {\cal A}$ onto $\op{Cut}_M^t$.
\end{proposition}
\begin{proof}
    By Corollary \ref{cor5}, we know that $\op{Cut}_M$ and $\op{Cut}_M^t$ are contained in $J^+\backslash {\cal A}$ and $I^+\backslash {\cal A}$, respectively.
    
    We need to prove existence of a continuous homotopy
    \begin{align*}
        H:[0,1]\times J^+\backslash {\cal A}\to J^+\backslash {\cal A}
    \end{align*} 
    satisfying the properties
    \begin{itemize}
        \item $H(0,x,y)=(x,y)$ for all $(x,y)\in J^+\backslash {\cal A}$.
        \item $H(1,x,y)\in \op{Cut}_M$ for all $(x,y)\in J^+\backslash {\cal A}$.
        \item $H(t,x,y)=(x,y)$ for all $(x,y)\in \op{Cut}_M$ and $t\in [0,1]$.
        \item $H(t,x,y)\in I^+$ for all $(x,y)\in I^+\backslash {\cal A}$ and $t\in [0,1]$.
    \end{itemize}
    
     We define $H:[0,1]\times J^+\backslash {\cal A}\to J^+$ by 
     \begin{align*}
     &H(t,x,y):=
     \\
     &\big(\exp_{x}(-t\beta(x,y)\varphi^-(x,y)v),\exp_x([(1-t\beta(x,y))+t\beta(x,y)\varphi^+(x,y)]v)\big).
    \end{align*} 
    Here, $[0,1]\ni t\mapsto \exp_x(tv)$ is a maximizing geodesic connecting $x$ to $y$.
    
    $H$ is well-defined since $\exp_x$ is defined at 
    \[
    \beta(x,y)\varphi^-(x,y)v \text{ and } [(1-\beta(x,y))+\beta(x,y)\varphi^+(x,y)]v.
    \]
    In the case where multiple maximizing geodesics connect $x$ to $y$, Theorem \ref{thme} implies that $y$ is the cut point of $x$ along both of them. Thus, $\beta(x,y)=0$, hence both possible definitions reduce to $H(t,x,y)=(x,y)$ for all $t\in [0,1]$.

    To prove continuity, suppose that $J^+\backslash {\cal A} \ni (x_k,y_k)\to (x,y)\in J^+\backslash {\cal A}$, and let $[0,1]\ni t\mapsto \exp_{x_k}(tv_k)$ be maximizing geodesics connecting $x_k$ to $y_k$. Then, after passing to a subsequence that we do not relabel, we have $(x_k,v_k)\to (x,v)$ and $[0,1]\ni t\mapsto \exp_x(tv)$ is a maximizing geodesic connecting $x$ and $y$. From the formula above, using the continuity of $\beta$ and $\varphi^\pm$, we conclude
    \begin{align*}
        H(t_k,x_k,y_k)\xrightarrow{k\to \infty} H(t,x,y),
    \end{align*}
    proving continuity.
    
    Clearly, $H(0,x,y)=(x,y)$. By the above lemma, $H(1,x,y)\in \op{Cut}_M$ for all $(x,y)\in J^+\backslash {\cal A}$. Furthermore, $H(t,x,y)=(x,y)$ for all $(x,y)\in \op{Cut}_M$ and $t\in [0,1]$, since in this case $\beta(x,y)=0$. Clearly, $H(t,x,y)\in I^+$ whenever $(x,y)\in I^+\backslash {\cal A}$ and $t\in [0,1]$. To conclude the proof, it suffices to prove that $H(t,x,y)\in J^+\backslash {\cal A}$ for $(t,x,y)\in [0,1]\times J^+\backslash {\cal A}$: Up to reparametrization, the maximal geodesic extension of the maximizing geodesic $\exp_x(tv)$, $t\in [0,1]$, used in the definition of $H(t,x,y)$, is the maximal extension of a maximizing geodesic connecting $p_1\circ H(t,x,y)$ to $p_2\circ H(t,x,y)$. Hence, by Lemma \ref{ls} and $(x,y)\notin {\cal A}$, $H(t,x,y)\in J^+\backslash {\cal A}$.
\end{proof}

\begin{lemma}
    The Aubry set ${\cal A}\subseteq M\times M$ is closed.
\end{lemma}
\begin{proof}
    Let $(x,y)\in J^+\backslash {\cal A}$. Suppose, for contradiction, that there exists a sequence $(x_k,y_k)\in {\cal A}$ converging to $(x,y)$. Let $[0,1]\ni t\mapsto \exp_{x_k}(tv_k)$ be maximizing geodesics connecting $x_k$ to $y_k$. After passing to a subsequence, we have $(x_k,v_k)\to (x,v)$, and the curve $[0,1]\ni t\mapsto \exp_x(tv)$ is a maximizing geodesic connecting $x$ to $y$. 
    
    Since $(x,y)\notin {\cal A}$, the maximal extension $\gamma:I\to M$ of this geodesic is not maximizing. Let $\gamma_k:I_k\to M$ denote the maximal extension of the geodesic connecting $x_k$ to $y_k$. Since $(x_k,y_k)\in {\cal A}$, the geodesics $\gamma_k$ are globally maximizing, and $\gamma_k$ converges to $\gamma$ in the $C^1$-topology on every compact subinterval $[a,b]\subseteq I$. By continuity of $d$, it is easy to conclude that $\gamma$ must be maximizing on each such subinterval, and therefore maximizing on the entire interval $I$. This is a contradiction.
\end{proof}

\begin{proposition}\label{pro5}
    The inclusion ${\cal NU}\hookrightarrow \op{Cut}_M$ is a homotopy equivalence, which restricts to a homotopy equivalence ${\cal NU}^t\hookrightarrow \op{Cut}_M^t$.
\end{proposition}
\begin{proof}
    Thanks to the previous lemma, we can construct a continuous function $T:(M\times M)\backslash {\cal A}\to (0,\infty)$ such that $(x,\varphi(t,y))\in (M\times M)\backslash {\cal A}$ for all $t\in [0,T(x,y)]$. 
    We define the continuous map
    \begin{align*}
        G:[0,1]\times \op{Cut}_M\to (M\times M)\backslash {\cal A},
        \ (t,y)\mapsto (x,\varphi(tT(x,y),y))
    \end{align*}
    Let $\ep:\op{Cut}_M^t\to (0,\infty),\ \ep(x,y):=d_{h\times h}((x,y),\partial J^+)$, and let $\bar F$ and $s:\op{Cut}_M^t\to (0,\infty)$ denote the maps from Corollary \ref{cor4} associated with $\ep$. We denote by $\tilde H:[0,1]\times (J^+\backslash {\cal A})\to J^+\backslash {\cal A}$ the strong deformation retraction from Proposition \ref{pro4}.
    Define
    \begin{align*}
        &H:[0,1]\times \op{Cut}_M
        \to \op{Cut}_M,
        \\[10pt]
        &H(t,x,y):=
        \begin{cases}
            (x,y),\ &\text{if } t=0,
            \\\\
            (\tilde H(1,G(t,x,y)),\bar F(t,\tilde H(1,G(t,x,y)))),\ &\text{ if } t>0.
        \end{cases}
    \end{align*}
    Since $G(t,x,y)\in I^+\backslash {\cal A}$ for $t>0$, we have $\tilde H(1,G(t,x,y))\in \op{Cut}_M^t$ by Propostion \ref{pro4}. Hence, $H$ is well-defined and 
    \begin{align}
        H(t,x,y)\in {\cal NU}^t\subseteq {\cal NU} \text{ for all } t>0. \label{eqak}
    \end{align}
    We claim that $H$ is continuous. As a composition of continuous functions, continuity holds on $(0,1]\times \op{Cut}_M$. For $(t,x,y)\in [0,1]\times \op{Cut}_M^t$, we have $\tilde H(1,G(t,x,y))\in \op{Cut}_M^t$ thanks to Proposition \ref{pro4}. Since $\bar F$ is well-defined and continuous on $[0,1]\times \op{Cut}_M^t$ and $\tilde H(1,G(0,x,y))=(x,y)$ for $(x,y)\in \op{Cut}_M$, it is immediate that also $H$ is continuous on the open set $[0,1]\times \op{Cut}_M^t$. Finally, if $t=0$ and $(x,y)\in \op{Cut}_M^n$, and $(0,1]\times \op{Cut}_M\ni (t_k,x_k,y_k)$ converges to $(0,x,y)$, then $(x_k',y_k'):=\tilde H(1,G(t_k,x_k,y_k))$ converges to $\tilde H(1,x,y)=(x,y)$. Hence,
    \begin{align*}
        d_{h\times h}(H(t_k,x_k,y_k),(x,y))
        &\leq 
        d_{h\times h}(H(t_k,x_k,y_k),(x_k',y_k'))+d_{h\times h}((x_k',y_k'),(x,y))
        \\
        &\leq 
        \ep(x_k',y_k')+d_{h\times h}((x_k',y_k'),(x,y))
        \xrightarrow{k\to \infty} 0.
    \end{align*}
    In the last step, we used $(x_k',y_k')\to (x,y)\in \partial J^+$. This proves the continuity.
    
    We claim that $H(1,\cdot)$ is a homotopy inverse to the inclusion $\iota: {\cal NU}\hookrightarrow \op{Cut}_M$, and that the restriction, $H(1,\cdot)_{|\op{Cut}_M^t}$, is a homotopy inverse to the inclusion $\iota^t:{\cal NU}^t\hookrightarrow \op{Cut}^t_M$.

    Indeed, by \eqref{eqak}, $H(1,\cdot)$ maps to ${\cal NU}$. The map $\iota\circ H(1,\cdot)=H(1,\cdot)$ is homotopic to $\op{Id}_{\op{Cut}_M}=H(0,\cdot)$ in $\op{Cut}_M$ via the homotopy $H$. Conversely, using \eqref{eqak}, we have $H(t,x,y)\in {\cal NU}$ whenever $(x,y)\in {\cal NU}$, so the composition $H(1,\cdot)\circ \iota$ is homotopic to $\op{Id}_{\cal NU}$ via the homotopy $(t,x,y)\mapsto H(t,\iota(x,y))$. Using that $H(t,x,y)\in {\cal NU}^t$ whenever $t>0$ (see \eqref{eqak}), the same arguments show that also $H(1,\cdot)_{|\op{Cut}_M^t}$ is a homotopy inverse to $\iota^t$.
\end{proof}

\begin{proof}[Proof of Theorem \ref{b}(a)]
   As in the proof of part (b), the result now follows from the preceding proposition in conjunction with Proposition \ref{pro4}
\end{proof}

\section{Appendix}
We want to prove the following lemma:
\begin{lemma}
\begin{enumerate}[(a)]
    \item The function 
    \begin{align*}
        \C:(0,\infty)\times M\times M\to \R\cup\{+\infty\},\ (t,x,y)\mapsto c_t(x,y),
    \end{align*}
    is real-valued, continuous on $(0,\infty)\times J^+$, and locally semiconcave on $(0,\infty)\times I^+$.
    \item If $x\in M$ and $y\in I^+(x)$, then the set of super-differentials of $\C$ at the point $(t,x,y)$ is given by
    \begin{align*}
        \partial^+ \C(t,x,y)
        =
        \op{conv}\left(\bigg\{\left(\partial_t c_t(x,y),-\frac{\partial L}{\partial v}(x,\dot \gamma(0)), \frac{\partial L}{\partial v}(y,\dot \gamma(t))\right)\bigg\}\right),
    \end{align*}
    where the set runs over all maximizing geodesics $\gamma:[0,t]\to M$ connecting $x$ to $y$.

    In particular, $\C$ is differentiable at $(t,x,y)$ if and only if there is a unique maximizing geodesic connecting $x$ to $y$ in time $t$ (equivalently, in time $1$).
\end{enumerate}
\end{lemma}
\begin{proof}
    The continuity in part (a) is trivial. For the rest, note that it suffices to check that $c_1$ is locally semiconcave on $I^+$ and that 
    \begin{align}
        \partial^+ c_1(x,y)
        =
        \op{conv}\left(\bigg\{\left(-\frac{\partial L}{\partial v}(x,\dot \gamma(0)), \frac{\partial L}{\partial v}(y,\dot \gamma(t))\right)\bigg\}\right), \label{eqadd}
    \end{align}
    where the set runs over all maximizing geodesics $\gamma:[0,1]\to M$ connecting $x$ to $y$.

    The proof is oriented towards \cite{Fathi/Figalli}, Theorem B19. We use a similar strategy and notation.

Let $(x_0,y_0)\in I^+$. Let $(U,\phi_1)$ and $(V,\phi_2)$ be two charts around $x_0$ and $y_0$ respestively with $\phi_1(x_0)=0$, $\phi_1(U)=\R^n$ and $\phi_2(y_0)=0$, $\phi_2(V)=\R^n$ and $\overbar U\times \overbar V\subseteq I^+$.
Set $U_1:=\phi_1^{-1}(B_1(0))$ and $V_1:=\phi_2^{-1}(B_1(0))$. 

From Lemma \ref{lm} (whose proof only requires the well-known properties), we know that the set $\Gamma$ of all maximizing geodesics $\gamma:[0,1]\to M$ with $\gamma(0)\in \overbar U$ and $\gamma(1)\in \overbar V$ is compact. It follows that there exists $\ep\in (0,1)$ and a compact set $K\subseteq \op{int}(\C)$ such that, for all  $\gamma\in \Gamma$, it holds
\begin{enumerate}[(1)]
\item
 $\gamma([0,\ep])\subseteq  \phi_1^{-1}(B_2(0)) \text{ and } \gamma([1-\ep,1])\subseteq \phi_2^{-1}(B_2(0))$. \label{dgaidwuasasasafioas}
\item $(\gamma(t),\dot \gamma(t))\in K$ for all $t\in [0,1]$.
\end{enumerate}
 In particular, there exist two compact sets $K_1\subseteq T\phi_1(TU\cap \op{int}(\C))\subseteq \R^n\times \R^n$ and $K_2\subseteq T\phi_2(TV\cap \op{int}(\C))\subseteq \R^n\times \R^n$ such that $T\phi_1(\gamma(t),\dot \gamma(t))\in K_1$ for all $\gamma$ and all $t\in [0,\ep]$ and $T\phi_2(\gamma(t),\dot \gamma(t))\in K_2$ for all $\gamma\in \Gamma$ and all $t\in [1-\ep,1]$.

Then let $\delta>0$ such that
\begin{align*}
   B_{\frac{2\delta}\ep}(K_1)\Subset T\phi_1(TU\cap \op{int}(\C))
   \text{ and } 
   B_{\frac{2\delta}\ep}(K_2)\Subset  T\phi_2(TV\cap \op{int}(\C)).  
\end{align*}
Now, if $\gamma\in \Gamma$, and $h:=(h_1,h_2)\in \R^n\times \R^n$ with $|h_1|,|h_2|\leq \delta$, let us set
\begin{align}
 \gamma_h:[0,1]\to ,\  \gamma_h(t):=
\begin{cases}
 \phi_1^{-1}\(\frac{\ep-t}{\ep}h_1+\phi_1(\gamma(t))\),\ &t\leq \ep,
\\[10pt]
\gamma(t),\ &t\in [\ep,1-\ep],
\\[10pt]
\phi_2^{-1}\(\frac{t-(1-\ep)}{\ep}h_2+\phi_2(\gamma(t))\),\ &t\geq 1-\ep.
\end{cases} \label{hdaiujoaiopd}
\end{align}
Note that
\begin{align*}
    &T\phi_1(  \gamma_h(t), \dot \gamma_h(t))\in B_{\frac{2\delta}\ep}(K_1),\ t\in [0,\ep], \text{ and } 
    \\[10pt]
    &T\phi_2( \gamma_h(t), \dot \gamma_h(t))\in B_{\frac{2\delta}\ep}(K_1),\  t\in [1-\ep,1].
\end{align*}

Given $(x_1,y_1),(x_2,y_2)\in B_{\frac{\delta}{2}}(0)\times B_{\frac{\delta}{2}}(0)$ set $h_1:=x_2-x_1$ and $h_2:=y_2-y_1$. Let $\gamma:[0,1]\to M$ be a maximizing geodesic connecting $\phi_1^{-1}(x_1)$ to $\phi_2^{-1}(y_1)$, let $h:=(h_1,h_2)$ and let $\gamma_h$ be the piecewise smooth curve as in \eqref{hdaiujoaiopd}. We can then estimate 
\begin{align*}
&c_1(\phi_1^{-1}(x_2),\phi_2^{-1}(y_2))-c_1(\phi_1^{-1}(x_1),\phi_2^{-1}(y_1))
\\[10pt]
\leq \ &\int_0^1 L({\gamma_h}(t),\dot {\gamma}_h(t))\, dt-\int_0^1 L(\gamma(t),\dot \gamma(t))\, dt
\\[10pt]
=\ &\int_0^\ep L({\gamma_h}(t),\dot {\gamma}_h(t))-L(\gamma(t),\dot \gamma(t))\, dt+\int_{1-\ep}^{1} L({ \gamma_h}(t),\dot { \gamma}_h(t))-L(\gamma(t),\dot \gamma(t))\, dt.
\end{align*}
We deal with the first integral $I_1$ exclusively since the second can be treated analogously. As in \cite{Fathi/Figalli}, we define the new Lagrangian
\begin{align*}
 L_{1}:\R^{n}\times \R^{n}\to \R,\ L_{1}(x,v):=L(\phi_1^{-1}(x),d_x\phi_1^{-1}(v))
\end{align*}
and the new curves
\begin{align*}
    \gamma_1:=\phi_1\circ \gamma \text{ and } \gamma_{h1}=\phi_1 \circ \gamma_h,
\end{align*}
so that 
\begin{align*}
    I_1=\int_0^\ep  L_{1}( \gamma_{h1}(t),\dot {\gamma}_{h1}(t))- L_{1}(\gamma_{h1}(t),\dot \gamma_{h1}(t))\, dt.
\end{align*}
 The Lagrangian $L_1$ is smooth on $T\phi_x(\op{int}(\C)\cap TU)$, so we can find a Lipschitz constant $C_1$ for the derivative $DL_1$ restricted to $B_{\frac {2\delta}\ep}(T\phi_1(K_1))$. Using the mean value theorem 
\begin{align*}
I_1
\leq
\int_0^\ep 
D L_{1}(\gamma_1(t),\dot { \gamma}_1(t))\Big[\frac{\ep-t}{\ep}h_1,-\frac{h_1}{\ep}\Big]\, dt +\frac{C_1|h_1|^2}{\ep}.
\end{align*}
We do the same computation for the second integral. With obvious notations, we obtain
\begin{align*}
&c_1(\phi_1^{-1}(x_2),\phi_2^{-1}(y_2))-c_1(\phi_1^{-1}(x_1),\phi_2^{-1}(y_1))
\\[10pt]
\leq
&\int_0^\ep 
DL_1(\gamma_1(t),\dot { \gamma}_1(t))\Big[\frac{\ep-t}{\ep}h_1,-\frac{h_1}{\ep}\Big]\, dt +
\int_{1-\ep}^1 
DL_2(\gamma_2(t),\dot { \gamma}_2(t))\Big[\frac{t-(1-\ep)}{\ep}h_2,\frac{h_2}{\ep}\Big]\, dt
\\[10pt]
+&\frac{C_1|h_1|^2}{\ep}+\frac{C_2|h_2|^2}{\ep}.
\end{align*}
Since $C_1$ and $C_2$ are independent of $(x_1,y_1),(x_2,y_2)\in B_{\frac{\delta}{2}}(0)\times B_{\frac{\delta}{2}}(0)$, and since the two integrals are linear in $h_1$ and $h_2$, respectively, this proves the local semiconcavity.

Moreover, using the Euler-Lagrange equation for timelike $L$-minimizers (recall that $L$ is smooth on $\op{int}(\C)$) and integrating the above integrals by parts, this proof also shows that a super-differential of $c_1$ at some point $(x,y)\in I^+$ is given by
\begin{align*}
\left(-\frac{\partial L}{\partial v}(x,\dot \gamma(0)),
\frac{\partial L}{\partial v}(y,\dot \gamma(1))\right),
\end{align*}          
where $\gamma:[0,1]\to M$ is a maximizing geodesic connecting $x$ to $y$ (see also Corollary B20 in \cite{Fathi/Figalli}).

In particular, since $\partial^+c_1(x,y)$ is a convex set, we proved $\supseteq$ in \eqref{eqadd}. However, it is well-known (\cite{Philippis}, Proposition A.3) that $\partial^+c_1(x,y)$ is given by the convex hull of reaching gradients, that is, it is the convex hull of all covectors $(p,q)\in T_x^*M\times T_y^*M$ of the form
\[
(p,q)=\lim_{k\to \infty} d_{(x_k,y_k)}c_1
\]
where $(x_k,y_k)\in I^+$ is a sequence of differentiability points of $c_1$ converging to $(x,y)$. However, if $\gamma_k:[0,1]\to M$ are maximizing geodesics connecting $x_k$ to $y_k$, by the above we must have
\begin{align*}
d_{(x_k,y_k)}c_1
= 
\left(-\frac{\partial L}{\partial v}(x_k,\dot \gamma_k(0)),
\frac{\partial L}{\partial v}(y_k,\dot \gamma_k(1))\right).
\end{align*}
Thanks to Lemma \ref{lm} we get that, along a subsequence, $\gamma_k$ converges in the $C^1$-topology to some maximizing geodesic $\gamma:[0,1]\to M$ connecting $x$ to $y$. Thus, we get
\begin{align*}
(p,q)=
\left(-\frac{\partial L}{\partial v}(x,\dot \gamma(0)),
\frac{\partial L}{\partial v}(y,\dot \gamma(1))\right).
\end{align*}          
This proves $\subseteq$ in \eqref{eqadd}, concluding the proof.
\end{proof}

\section*{Acknowledgements}
I’m especially grateful to Stefan Suhr, who first guided me toward the field of Lorentzian weak KAM theory. The topic and results of this paper are based on his original idea, and I truly appreciate his ongoing advice and the support he provided throughout the writing process.
I’d also like to thank my supervisor, Markus Kunze, who has supported me from the very beginning of my PhD.

\bibliography{maximalgeodesics}
\end{document}